\documentclass[11pt,reqno,twoside]{amsart}

\usepackage{multicol}
\usepackage{graphicx}
\usepackage{amsmath}
\usepackage{amssymb}
\usepackage{xcolor}
\usepackage{fancyhdr}
\usepackage{mathrsfs}
\usepackage{mathabx}
\usepackage{cancel}
\usepackage{placeins}
\usepackage{afterpage}
\usepackage{flafter}

\usepackage{float}
\usepackage{algorithm}
\usepackage{algpseudocode}

\makeatletter
\algrenewcommand\alglinenumber[1]{\makebox[2.2em][l]{\footnotesize #1}}
\makeatother

\usepackage{booktabs}
\usepackage{tikz}
\usepackage{comment}
\usepackage{subcaption}
\usepackage[normalem]{ulem}
\usepackage[colorlinks=true,
 linkcolor=blue,
 citecolor=blue,
 urlcolor=blue]{hyperref}

\newcounter{algorithmfourb}
\renewcommand{\thealgorithmfourb}{4B}

\usetikzlibrary{arrows.meta, positioning, calc}

\newcommand{\sech}{\operatorname{sech}}

\theoremstyle{remark}
\newtheorem{remark}{Remark}

\theoremstyle{remark}

\def\R{{\mathbb{R}}}

\title[Ground states with PINNS \& Deep Ritz]
{Soliton Profiles in 1D: Enhanced Neural-Network Solvers\\ 
and Constrained Deep Ritz Methods 
}

\author[C. Haight]{Chandler Haight}
\address{Department of Mathematics \& Statistics, Florida International University, Miami, FL 33199, USA}
\email{chaig004@fiu.edu}

\author[A. D. Rodriguez]{Alex David Rodriguez}\address{Department of Mathematics \& Statistics, Florida International University, Miami, FL 33199, USA}
\email{arodr1128@fiu.edu}

\author[S. Roudenko]{Svetlana Roudenko}
\address{Department of Mathematics \& Statistics, Florida International University, Miami, FL 33199, USA}
\email{sroudenko@fiu.edu}

\date{}



\begin{document}

\begin{abstract} 
We develop and compare {\it enhanced} neural-network solvers for one-dimensional ground-state profiles, which are solutions of 
$$
Q^{\prime\prime}=bQ-\gamma Q^{p},\quad Q>0,\quad p>1, \quad b, \gamma >0,
$$
the stationary profile equation for the nonlinear Schr\"odinger, generalized Korteweg--de Vries, and nonlinear Klein--Gordon equations. Its explicit $\sech$-type solutions serve as exact reference for error computation { and benchmark hyperparameter selection}. We improve on the standard Adam-trained physics-informed neural network (PINN) baseline by introducing higher-order optimization via Gauss--Newton/Levenberg--Marquardt and L--BFGS refinements, as well as Fourier-feature PINNs, and constrained Deep Ritz (DR) formulations. In DR approaches we introduce two new neural-network formulations, {which implement two classical variational characterizations:} 
a {\it Nehari-constrained} Deep Ritz method, in which the height is learned by projecting onto the Nehari/Pokhozhaev manifold, and a {\it Weinstein} Deep Ritz method, in which { a normalized} profile minimizes the Weinstein quotient.  
The Deep Ritz formulations perform best overall, in particular, the Weinstein variant gives the best error--runtime balance. 
The neural-network based solvers remain less accurate and efficient than classical Petviashvili iteration in 1D, {however, the current study provides a proof of concept: the mesh-free formulation, non-reliance on exact profiles or the related features}  make them promising in higher dimensions and for more complicated nonlinearities and dispersive operators.

\end{abstract}

\maketitle

\hypersetup{linkcolor=black}
\tableofcontents
\hypersetup{linkcolor=magenta}

\section{Introduction}\label{S:intro}

Ground-state and solitary-wave solutions play a profound role in the study of nonlinear dispersive, or wave-type, partial differential equations: 
their profiles and their various norms serve as key elements for understanding wave propagation,   
expressing thresholds in long-time solutions dynamics (to a given PDE), they are the main players in understanding 
stability or instability mechanisms, and fundamental in the description of singularity formation such as self-similar collapses and blow-ups. 
The ground state profile equation, describing the stationary ground state or solitary wave profiles (see derivation for several nonlinear dispersive models in Section \ref{S:models}), reduces to a nonlinear elliptic-type differential equation, whose analytical solutions are often difficult to obtain or inaccessible with current analytical methods, and thus, accurate numerical solution is highly desirable. 

In the one-dimensional setting with a single-power nonlinearity, we are interested in ground state solutions of the second-order ODE 
\begin{equation}\label{E:intro-ODE}
\qquad \qquad Q^{\prime\prime} = bQ - \gamma |Q|^{p-1}Q, \quad p>1, \quad b >0, \quad \gamma>0, 
\end{equation}
which are positive, smooth, decaying at infinity solutions (e.g., in $H^1(\mathbb R)$), thus, the previous equation can be simplified to 
\begin{equation}\label{E:intro-ODE2}
\qquad \qquad Q^{\prime\prime} = bQ - \gamma Q^p, \qquad Q>0, \quad Q'(0)=0, \quad Q(\infty)=0.
\end{equation}

This equation arises, for example, from standing waves in the nonlinear Schr\"odinger (NLS) equation, traveling solitary waves in the generalized Korteweg--de Vries (gKdV) equation and in the nonlinear Klein--Gordon (NLKG) equation. A key feature of the one-dimensional case is that it admits an explicit exact
solution
$$
Q(x)=A\,\sech^{\frac2{p-1}}(\beta x),
$$
with an amplitude $A$ and inverse width $\beta$ depending on given $b$, $\gamma$, and $p$, see \eqref{E:coeff}. This explicit ground truth makes the 1D case an ideal testing and verification setting for
numerical benchmarking.  

For the purpose of this paper, we consider integer exponents $p>1$, typically studying $p=3$ (subcritical), $p=5$ (critical) and $p=7$ (supercritical) cases.

In our previous work \cite{HRW2026}, we compared neural-network solvers 
-- PINNs and operator learning approaches such as DeepONets and Fourier Neural Operators -- with classical methods such as finite-difference with Newton's iterations and Petviashvili method, for computing these 1D profiles. That study showed that while neural-network methods can approximate solitary-wave profiles, classical Petviashvili-based solvers remained substantially more accurate and efficient.

The present paper aims to shrink this gap by introducing, 
developing and evaluating {\it enhanced} neural-network-based algorithms for the ground-state equation \eqref{E:intro-ODE}. 
{ This explicitly solvable 1D problem serves as a proof of concept and a controlled benchmark for neural solvers, rather than a setting in which they are expected to replace the Petviashvili iteration. Before advancing to more complicated problems, we, therefore, concisely survey the main methods, introduce new features, and compare them in this setting, in the hope that it will allow future development and comparison in more complicated cases and models.}
Specifically, we propose three directions of improvement:
\begin{itemize}
\item {\it Higher-order optimization}: replacing or augmenting standard Adam training with Gauss--Newton (GN) and Levenberg--Marquardt (LM) refinement strategies, applied to the PINN residual after an initial Adam warm-up phase;
\item {\it Fourier-feature PINNs (FF--PINNs)}: embedding the scalar coordinate input into a space with Fourier basis 
before passing it to the MLP (multilayer perceptron), in order to reduce the spectral bias of standard coordinate-based networks;
\item {\it Constrained Deep Ritz formulations}: several variational methods that minimize the action functional or related quotients over constrained neural-network trial classes (instead of the profile equation itself), including an action-based constrained Deep Ritz formulation, a Nehari-constrained Deep Ritz formulation, and a Weinstein-quotient-based Deep Ritz formulation.
\end{itemize}

The last point is worth commenting on further as we { develop neural implementations via  Deep Ritz (DR) formulations of two classical variational characterizations} for obtaining ground state solutions in focusing dispersive PDE: the {\it Nehari-constrained} (NDR) and the {\it Weinstein} Deep Ritz (WDR).

The WDR formulation is based on minimizing the Weinstein functional quotient (sometimes also called Lyapunov functional or nonlinear Rayleigh, e.g., see \cite{KR2026}) defined as 
\begin{equation}\label{E:intro-Weinstein}
R_b = \inf_{u\in H^1(\R)\setminus\{0\}}
\frac{\|u'\|_{L^2}^2 + b\|u\|_{L^2}^2}{\|u\|_{L^{p+1}}^2},
\end{equation}
over a normalized, hard-constrained neural-network trial class. One of the main advantages of the WDR formulation, compared to the action-based constrained Deep Ritz method, is that WDR method does not require prescribing the peak amplitude $Q(0)$ from the exact solution {(so its training objective can be formed without an explicit profile)}, thus, making it particularly suitable for problems where no closed-form solution is available. The first method, the {\it constrained Deep Ritz}, in a general class of solutions may not have a converging result, as it relies on the initial maximum value $A=Q(0)$, which in its turn helps selecting the correct {\it branch} of positive ground state solutions. We first improve that in the {\it Nehari-constrained Deep Ritz} (NDR) formulation, by incorporating the first Pokhozhaev (also called Nehari) identity \eqref{E:P1} into the formulation and the second Pokhozhaev identity \eqref{E:P2} into the loss. We then observe that changing the action functional optimization to the Weinstein quotient \eqref{E:intro-Weinstein} optimization further improves the results ({ e.g., gives a more favorable error-runtime balance}).

Our results show that these enhanced methods substantially improve upon the standard Adam-trained PINN baseline. Among the methods tested, the { complete} Deep Ritz {procedures, including their final residual refinement},  give the strongest overall performance, { with the new Weinstein Deep Ritz formulation} achieving the best {(among DR variants)} error--runtime for $p=3$ (cubic) and $p=5$ (quintic). For the sharper $p=7$ (septic) profile, the Weinstein and Nehari-constrained formulations produce errors on the same order, and better compared with the other formulations. 
The PINN and FF--PINN variants also benefit from higher-order optimization and improved feature representation, but they remain more sensitive to loss formulation, optimizer choice, initialization, and computational device. Overall, although the enhanced neural-network solvers do not yet match the accuracy of Petviashvili-type iterations reported in~\cite{HRW2026} for the one-dimensional setting, they represent a significant step forward in neural-network-based computation of nonlinear ground states. We expect that in higher-dimensional problems or with more intricate nonlinearities or dispersion, the neural-network-based methods will be an essential addition to the field. 

\smallskip

The paper is organized as follows. In Section~\ref{S:groundstate} we set up the ground-state equation and its variational structure, describe the explicit solutions, boundary conditions, domain truncation, the Pokhozhaev identities, and introduce the Weinstein variational characterization used in the Deep Ritz formulations. Section~\ref{S:PINNs} covers PINNs and FF--PINNs, including the use of second-order optimization and comparisons with standard Adam training.
In Section~\ref{S:DEEP_RITZ} we study Deep Ritz formulations, including constrained Deep Ritz, Nehari-constrained Deep Ritz, 
and Weinstein-type Deep Ritz methods.
Section~\ref{S:Conclusion} concludes with a discussion of the main findings.
\smallskip

{\bf Acknowledgments.} We thank Javier G\'omez-Serrano for helpful discussions and for suggesting the use of a Fourier-feature-type PINN architecture and higher-order optimization strategies, such as Gauss--Newton, quasi--Newton, and Adam followed by higher-order refinement, as possible ways to improve upon the standard Adam-trained PINN baseline motivated by the strong Petviashvili performance observed in our previous work.

The training on basic neural-network-based approaches for PDE came out of the STEM summer ML seminar that S.R. co-organized in Summer 2025 (with Zhongming Wang), which was partially supported by the NSF grant DUE/EDU-2221491 (PI: S. Roudenko); A.D.R. and C.H. participated in that seminar; 
A.D.R., C.H. and S.R. were partially supported by the NSF grants DUE/EDU-2221491 and DMS-2452782. 

{\bf  Declaration of usage of Large Language Models.} During the preparation of this work the authors used GPT-5.5 and Claude Opus 4.7 in order to assist with code development and debugging. GPT-5.6 Sol and Claude Opus 4.7 were later used for grammar 
and language editing of the manuscript. 
While using these tools, the authors reviewed and decided whether the content is relevant, and if it was, it was re-worked and edited before incorporating into the paper; the authors take full responsibility for the content of the article.


\section{Ground-state equation and problem setup}\label{S:groundstate}

In this section, we describe the profile equation, its origin from the NLS, gKdV (we write this equation with an absolute value in nonlinearity for consistency of notation), and NLKG models, the explicit solutions, the finite-domain formulation used in the numerical experiments, and the variational structure underlying the Deep Ritz methods.

\subsection{Dispersive models and profile reduction}\label{S:models}
The nonlinear dispersive models considered in this paper are the focusing, power-type one-dimensional equations:
\begin{align}
    i u_t + u_{xx} + \gamma |u|^{p-1}u &= 0,
    \qquad &&\text{(NLS)} \label{NLS}\\
    u_t + \big( u_{xx} + \gamma |u|^{p-1}u \big)_x &= 0,
    \qquad &&\text{(gKdV)} \label{gKdV}\\
    u_{tt} - u_{xx} + b u - \gamma |u|^{p-1} u &= 0,
    \qquad &&\text{(NLKG)} \label{NLKG}
\end{align}
with parameters $b>0$, $\gamma>0$ and integer nonlinearity exponent $p>1$.

The sign convention in front of the nonlinear term is chosen so that these
equations are focusing, and hence, possess solitary-wave or ground-state solutions (the defocusing case with  the opposite sign does not admit ground states). These models arise in nonlinear optics, plasma physics, and water waves; see, for example, \cite{Agrawal2001,KA2003,Z1972,ZakRub1973,SS1999,KS,50years} and references therein. The solitary waves or ground states are spatially localized solutions that propagate without changing shape and often play a key role in characterizing the long-time dynamics of the flow, forming the main components of the soliton resolution conjecture~\cite{Soffer2006,Tao2009}.

Each of the models \eqref{NLS}--\eqref{NLKG} enjoys conservation laws on its
maximal interval of existence (see, e.g., the textbooks 
Cazenave~\cite{Cazenave2003}, Tao~\cite{Tao2006}, Linares--Ponce~\cite{LinaresPonce2015}). For the NLS equation \eqref{NLS}, the conserved mass and Hamiltonian (energy) are given by
\begin{equation}\label{E:mass}
    M[u(t)] = \int_{\mathbb{R}} |u(t,x)|^2 \, dx,
    \qquad
    E[u(t)] = \int_{\mathbb{R}} \Big( \tfrac{1}{2}|u_x|^2
    - \tfrac{\gamma}{p+1}|u|^{p+1} \Big) dx,
\end{equation}
while for the gKdV equation \eqref{gKdV} one has analogous conservation of the $L^2$-norm and
a suitable Hamiltonian functional. (There is also a momentum in NLS or an integral of the solution in gKdV conservation, but we omit those in this paper.)
In the NLKG equation \eqref{NLKG}, the conserved energy in the phase space $(u,u_t)$ is
\begin{equation}\label{E:energy}
    E[u(t),u_t(t)] = \int_{\mathbb{R}} \Big(
        \tfrac{1}{2}|u_t|^2 + \tfrac{1}{2}|u_x|^2
        + \tfrac{b}{2}|u|^2 - \tfrac{\gamma}{p+1}|u|^{p+1}
    \Big)\,dx.
\end{equation}

Solitary-wave or ground-state solutions are obtained by substituting special ansatzes that reduce the PDE to an ODE for a spatial profile $Q$. For the NLS equation \eqref{NLS}, the standing-wave ansatz $u(x,t) = e^{ibt}Q(x)$ yields the profile equation
\begin{equation}\label{E:1}
    Q^{\prime\prime} - b Q + \gamma |Q|^{p-1} Q = 0.
\end{equation}
For the gKdV equation \eqref{gKdV}, the traveling-wave ansatz $u(x,t)=Q(x-bt)$ with speed $b>0$ similarly gives
$Q^{\prime\prime} = bQ - \gamma |Q|^{p-1}Q$
after integration and use of the localization condition $Q(\xi)\to 0$ as $|\xi|\to\infty$.

For the NLKG equation \eqref{NLKG}, the traveling-wave ansatz $u(x,t)=Q(x-ct)$ with speed $c$ leads, after rescaling parameters, to the same equation \eqref{E:1}.

In all three models, we look for {\it positive} $H^1$ solutions of \eqref{E:1}, which are known to exist and are unique (e.g., see \cite{BerestyckiLions1983} and \cite{Kwong1989}), therefore, ground-state profiles $Q>0$ satisfy the second-order profile ODE
\begin{equation}\label{2ndODE}
    Q^{\prime\prime} = b Q - \gamma Q^p.
\end{equation}
Under suitable assumptions on $p$, these ground states can be constructed variationally as minimizers of some functionals such as the {\it action} functionals (with or without constraints) or the energy functional at fixed mass and are orbitally stable (in 1D for $p<5$) in the sense of Weinstein~\cite{Wein85modu,Wein1986Lya} and Grillakis--Shatah--Strauss~\cite{GSS87,GSS90}, for a review see \cite[\S 2.2]{KR2026}. For gKdV, well-posedness and scattering are discussed in the work of Kenig--Ponce--Vega~\cite{kenig1993gkdv}, and the solitary waves and their (in)stability, for example, in Bona--Souganidis--Strauss~\cite{BSS} and references therein.

In this paper we focus on the accurate numerical solutions of the profile equation~\eqref{2ndODE} via neural-network approaches enhanced with various optimizers and modified with different formulations to achieve more efficient and higher precision results compared to 
our previous study~\cite{HRW2026}.

\subsection{Problem formulation and explicit solutions}\label{S:problem_formulation}

We consider localized positive solutions of the one-dimensional nonlinear elliptic equation
\begin{align}
Q^{\prime\prime} &= b Q - \gamma Q^p, \label{2ndODE-again}\\ 
Q^{\prime}(0) &= 0,\quad Q(\infty)=0, \label{2ndODEBC}
\end{align}
with parameters $b>0$, $\gamma>0$ and exponent $p>1$. As we mentioned in the introduction, for the purpose of this paper, we consider integer exponents, typically, $p=3,5,7$.  The condition $Q^\prime(0)=0$ corresponds to the even symmetry of the full-line ground state, while $Q(\infty)=0$ imposes localization, which is exponential (in the case of a standard Laplacian or second order derivative operator).

Since $Q(x)\to 0$ as $x\to\infty$ and $p>1$, the nonlinear term in \eqref{2ndODE-again} is asymptotically lower order in the tail. Thus, for large $x$, the decay is governed by the linearized equation
$$
 -Q''+b Q=0,
$$
whose localized solution is proportional to 
$e^{-\sqrt{b}x}$, and thus, satisfies
$$
Q'(x)+\sqrt{b}\,Q(x)=0.
$$
Therefore, on the finite half-interval $[0,L]$, it is typical in numerical approaches to replace the second condition in \eqref{2ndODEBC} by the mixed asymptotic (sometimes called Robin) boundary condition
\begin{equation}\label{E:mixed-BC}
Q'(L)+\sqrt{b}\,Q(L)=0.
\end{equation}
We take that into consideration and discuss the boundary conditions in the next subsection.

As we mentioned in the introduction, the ground state profile equation \eqref{2ndODE-again} admits an explicit nontrivial solution of the form
\begin{equation}\label{Qexact}
Q(x) = A \,\sech^{\frac{2}{p-1}}(\beta x) \equiv Q_{\mathrm{exact}}(x),
\end{equation}
where
\begin{align}\label{E:coeff}
A &= A(b,\gamma,p) = 
    \left( \frac{b(p+1)}{2\gamma} \right)^{\frac{1}{p-1}},
\quad \beta = \beta(b,\gamma,p) = \frac{\sqrt{b}(p-1)}{2}.
\end{align}
This explicit profile is used both as a reference solution for computing errors and, in some neural-network formulations, as a way to prescribe the correct peak amplitude $Q(0)=A$. {(Specifically, in what follows, we use $Q_0=A$ as input in Dirichlet PINN, FF-PINN, and action-based constrained DR baselines, whereas in NDR and WDR we determine the amplitude variationally. We also point out that apart from the prescribed boundary value in those baselines, no interior solution samples enter the training losses. The present hyperparameter selection uses errors against the exact cubic profile, which we discuss next.)} 

The cubic case 
$$
Q^{\prime\prime}=Q-Q^3,
$$
corresponding to $b=\gamma=1$ and $p=3$, is used as the main benchmark for hyperparameter tuning. After suitable configurations are selected from the cubic benchmark, the optimized methods are applied typically to the cubic, quintic ($p=5$), and septic ($p=7$) cases in the main comparisons. It should be mentioned that for higher nonlinearities it might also be possible to optimize the parameters for each nonlinearity independently to improve errors and avoid anomalies that we observe later (e.g., see Figure ~\ref{fig:GN_converge_bad}).

\subsection{Domain truncation and boundary conditions}\label{S:BCs}

In our simulations, the infinite-domain problem is approximated on $[0,L]$. 
We use $L=30$, which is sufficiently large for the benchmark cases considered here; for example, when $b=\gamma=1$ and $p=3$, the exact value $Q(30)$ is on the order of $10^{-13}$. 
We note that although all our computations are performed on the interval $[0,\,30]$, the solution profiles are plotted only on $[0,\,10]$ in order to enlarge the details from the 
peak-amplitude region, which is typically the most challenging for the neural-network approximation to resolve accurately.

Consider the following finite-domain boundary conditions:
\begin{align}
Q(0)&=A,\qquad Q(L)=0, \label{E:Dirichlet}\\
Q'(0)&=0,\qquad Q(L)=0, \label{E:mixed}\\
Q(0)&=A,\qquad Q'(L)+\sqrt{b}\,Q(L)=0, \label{E:Robin_D}\\
Q'(0)&=0,\qquad Q'(L)+\sqrt{b}\,Q(L)=0. \label{E:Robin_m}
\end{align}
We refer to \eqref{E:Dirichlet} as the Dirichlet formulation,
\eqref{E:mixed} as the Neumann--Dirichlet formulation, \eqref{E:Robin_D}
as the Robin--Dirichlet formulation, and \eqref{E:Robin_m} as the
Neumann--Robin formulation.

On the truncated interval $[0,L]$, the decay condition at infinity may
be approximated either by the Dirichlet condition $Q(L)=0$ or by the asymptotic Robin condition $Q'(L)+\sqrt{b}\,Q(L)=0$.
Since $Q(x)=\mathcal{O}(e^{-\sqrt{b}\,x})$ as $x\to\infty$, the Robin
condition incorporates the leading-order behavior of the linearized tail and
is generally more accurate on short computational domains. In our computations, $b=1$ and $L\geq 15$ (typically $L=30$), so that $Q(L)$ is already
negligibly small relative to the observed neural-network approximation
errors. We therefore retain the
Dirichlet condition for simplicity (and use either \eqref{E:Dirichlet} or \eqref{E:mixed}), while noting that the Robin condition may
be advantageous when a shorter interval is used (e.g., $L=8$) or when $b$ is very small. We note that in higher-dimensional settings or different nonlinearities, such boundary conditions would need to be taken into account. 

\subsection{Pokhozhaev identities}\label{S:Pokhozhaev}

For later use in the constrained Deep Ritz formulations in Section~\ref{S:DEEP_RITZ}, we record the two integral identities satisfied by ground state solutions to \eqref{2ndODE-again}. In 1D, multiplying the equation by $Q$ and integrating by
parts gives
\begin{equation}\label{E:P1}
\mathcal{P}_1(Q)
:= b\|Q\|_{L^2}^2+\|\partial_x Q\|_{L^2}^2
-\gamma\|Q\|_{L^{p+1}}^{p+1}=0.
\end{equation}
The corresponding scaling identity (obtained by multiplying by $xQ'$ and integrating) is
\begin{equation}\label{E:P2}
\mathcal{P}_2(Q)
:= b\|Q\|_{L^2}^2-\|\partial_x Q\|_{L^2}^2
-\frac{2\gamma}{p+1}\|Q\|_{L^{p+1}}^{p+1}=0.
\end{equation}
We use these identities either { as penalties in the action-based and Nehari-constrained DR formulations, and as diagnostics for Weinstein DR formulation} in Section~\ref{S:DEEP_RITZ}.

{
For the purpose of the computations on the finite interval $[0,L]$, we denote the following quantities by 
$M=\int_0^L Q^2\,dx$, $K=\int_0^L|Q'|^2\,dx$, and $N=\int_0^L|Q|^{p+1}\,dx$, and rewrite the above Pokhozhaev identities, which now contain the boundary terms. Assuming a sufficiently smooth solution of the ODE, integration by parts gives the exact finite-interval identities, namely,} 
$$
P_{1,L}=K+bM-\gamma N+Q(0)Q'(0)-Q(L)Q'(L)=0,
$$
and
$$
P_{2,L}=bM-K-\frac{2\gamma}{p+1}N+L|Q'(L)|^2-b L Q(L)^2+\frac{2\gamma L}{p+1}|Q(L)|^{p+1}=0.
$$

\subsection{Variational characterization via the Weinstein quotient}\label{S:Weinstein}

The ground-state equation \eqref{2ndODE-again} admits a variational characterization that plays a key role in the Weinstein Deep Ritz method developed in Section~\ref{S:DEEP_RITZ}. Recalling the equation \eqref{2ndODE-again}, for a given $b>0$ we define the Weinstein functional quotient 
\begin{equation}\label{E:Weinstein-Rb}
R_b := \inf_{u\in H^1(\R)\setminus\{0\}}
\frac{\|u'\|_{L^2(\R)}^2 + b\|u\|_{L^2(\R)}^2}{\|u\|_{L^{p+1}(\R)}^2}.
\end{equation}
If $v$ is a positive optimizer of this quotient, normalized by $\|v\|_{L^{p+1}(\R)}=1$, then $q_b(v)=R_b$, where $q_b(u):=\|u'\|_{L^2}^2+b\|u\|_{L^2}^2$, and the Euler--Lagrange equation is 
\begin{equation}\label{E:Weinstein-EL}
-v^{\prime\prime}+bv = R_b v^p.
\end{equation}
The ground-state profile is then recovered by the rescaling
\begin{equation}\label{E:Weinstein-recovery}
Q = \left(\frac{R_b}{\gamma}\right)^{\frac{1}{p-1}} v.
\end{equation}
This {classical characterization is related (after spatial rescaling) to the sharp interpolation variational principle} of Weinstein~\cite{Weinstein1983}, and is especially convenient for a Deep Ritz implementation. The quotient \eqref{E:Weinstein-Rb} is invariant under multiplication of $u$ by a nonzero constant (and thus, does not select the trivial solution), and the normalization constraint $\|v\|_{L^{p+1}}=1$ removes the remaining scale invariance. For fixed $b$ and $p$, the variational training is independent of $\gamma$, which enters only through the final amplitude rescaling \eqref{E:Weinstein-recovery}. Crucially, this method does not require prescribing the peak amplitude $Q(0)=A$ from the exact solution \eqref{Qexact}, making it suitable for computing ground states when no closed-form solution is available.

\subsection{Related work}\label{S:related_work}

Numerical solutions to ODE have been extensively studied in the literature using classical approaches (see discussion in our previous paper \cite{HRW2026} on that and references therein).
Classical numerical methods, as we showed in \cite{HRW2026}, typically provide reliable accuracy and convergence for ground state solutions of~\eqref{2ndODE}, although usually require mesh design and problem-specific tuning. In contrast, neural-network-based solvers---including physics-informed neural networks (PINNs)~\cite{RPE19}---approximate the profile through trainable function classes and can be implemented in a mesh-free manner.

Recent work on PINNs has shown that their accuracy depends strongly on the optimization method and that first-order training followed by higher-order refinement can improve convergence \cite{RathoreLeiFrangellaLuUdell2024}. These findings are directly related to our comparison of Adam \cite{KingmaBa2015}, L--BFGS \cite{LiuNocedal1989}, and Adam--L--BFGS.
PINN training has also been treated as a nonlinear least-squares problem using Gauss-Newton/Levenberg Marquardt (GN/LM) methods, in which  damping is used to improve stability
\cite{Levenberg1944,Marquardt1963,ShahabMukhlashSusanto2026LM}. The successful use of first-order training followed by high-precision GN/LM refinement for challenging nonlinear PDE profiles is closely related to, and motivates, our Adam--GN/LM solver \cite{WangBennaniMartens2025}. 
Related work on Fourier features 
has shown that they help neural networks represent high-frequency and multiscale functions and can reduce spectral bias in PINNs 
\cite{Tancik2020,WangWangPerdikaris2021},
motivating our Fourier-feature PINN (FF--PINN) formulation.
For the variational methods,
the Deep Ritz approach trains a neural network by minimizing an energy or an action functional  rather than a differential-equation residual \cite{EYu2018}.
However, direct minimization may not select the desired ground-state branch, unless more conditions or constrains are added.
Recent neural approaches that incorporate 
normalization or projection constraints include \cite{BaoChangZhao2025,ChangWenZhao2025}, and together with recent Deep Ritz work for stationary Schr\"odinger problems in \cite{YaoLuLerman2026DeepRitzSchrodinger}, provide additional methods for the constrained neural variational methods considered here. 
Variational characterizations based, on the other hand, on sharp Gagliardo--Nirenberg inequalities and normalized gradient flows, and  yield approaches that can properly identify 
the desired ground states
\cite{Wein85modu,Wang2022,LiuYuanZhao2023,GomesPastor2021}. 
We consider some of these and suggest novel approaches to {construct and compare constrained neural solvers} to obtain the ground state solutions in the focusing dispersive PDE, for example, by including Nehari or Pokhozhaev manifolds or the Weinstein quotient. Neural-network implementations of these constrained formulations have recently been developed in \cite{ChangWenZhao2025}, who enforce the $L^{p+1}$-normalization and Nehari manifold constraints as hard constraints through normalization and projection layers. Their setting is the {\it defocusing} action ground state, for which the action is {\it bounded below} and an unconstrained formulation is also available. In the {\it focusing} problem, which we consider, the action is unbounded below (see Section~\ref{S:NehariDeepRitz}), so a constraint is not one option among several but a structural necessity. Our contribution is therefore not the variational characterizations themselves, which are classical, but their realization {(including the combined implementation and numerical assessment)} for the truncated profile equation \eqref{2ndODE}: hard-constrained trial classes that encode symmetry, positivity and the endpoint condition simultaneously, a boundary-corrected finite-interval Pokhozhaev penalty, as well as a staged Adam followed by Gauss--Newton/Levenberg--Marquardt residual refinement. A recent work \cite{LiuWangZhao2026} analyzes the $L^{p+1}$-normalized gradient flow with asymptotic Lagrange multiplier for computing the action ground state of the focusing NLS, {using a gradient-flow rather than a neural-network approach.
}

\subsection{Computational setup}\label{S:hardware}

All reported runtimes in the tables are recorded as elapsed wall-clock times and measured using \texttt{time.perf\_counter()}. Therefore, the reported CPU times should be interpreted as practical implementation-dependent runtimes on the stated hardware, rather than hardware-independent performance measurements. The computations were performed on a device equipped with a 64-bit Windows laptop with an AMD Ryzen 5 7535HS CPU at 3.30 GHz, 6 cores / 12 threads, 16~GB RAM (15.3~GB usable), integrated Radeon graphics, and a dedicated NVIDIA GeForce RTX 4050 Laptop GPU with 6~GB VRAM\footnote{The purpose is not to use the fastest possible (and costly) computational facilities, but rather make it sufficiently fast and individually accessible.}. In Sections \ref{S:secondaryPINN}, \ref{S:secondaryFFPINN}, and \ref{S:secondaryDeepRitz}, for comparison purposes we also use a ``secondary'' device equipped with an AMD Ryzen 7 9700X (3.8~GHz base, 5.5~GHz boost, 8 cores / 16 threads), 32~GB of RAM (31.1~GB usable), and a dedicated NVIDIA GeForce RTX 4070 Super GPU with 12~GB of VRAM. 

For CUDA-enabled experiments, the reported elapsed times reflect the end-to-end runtime of the implemented \texttt{PyTorch} computation on this device, including CPU-side launching and synchronization overheads as well as GPU computation. 
We therefore label all secondary-device timings as \emph{GPU+CPU} rather than GPU-only: the \texttt{CUDA} device carries the tensor and neural-network computations, but the optimizer control logic, backtracking decisions, checkpointing, logging, data transfers, and any calls to CPU-based linear algebra or \texttt{NumPy} routines remain on the CPU. The reported wall-clock time is thus the combined execution of GPU-accelerated tensor operations and CPU-side back-end work, and should not be read as a measure of processor performance alone. Accordingly, the primary- and secondary-device results should be interpreted as a comparison between the two complete workstation configurations rather than as an isolated comparison of CPU and GPU performance. This convention applies to Sections~\ref{S:secondaryPINN}, \ref{S:secondaryFFPINN}, \ref{S:secondaryDeepRitz}, \ref{S:secondary_nehari}, and \ref{S:secondary_weinstein}.

All reported computations use double precision (\texttt{float64}) and fixed \texttt{np.random.seed(42)} for reproducibility. The tabulated values correspond to the individual runs reported in the respective experiments.

\FloatBarrier
\section{PINNs Variations}\label{S:PINNs}

We first study residual-based neural-network solvers for
\eqref{2ndODE-again}. The common benchmark problem, boundary conditions,
finite-domain truncation, and timing protocol were specified in
Section~\ref{S:groundstate}. In this section, we focus on
Physics-Informed Neural Networks (PINNs) and Fourier-Feature PINNs
(FF--PINNs), together with the optimization strategies used to train them.

For each method, we first describe the neural-network ansatz and loss
function. We then perform a hyperparameter study on the {\it cubic} trial problem
$$
Q''=Q-Q^3.
$$
Then we use the selected configuration for analysis and comparison of the 
cubic, quintic and septic cases. The purpose of this section is to assess the
accuracy, robustness, and computational cost of residual-based neural-network solvers, including with the enhancements, before comparing them with the variational Deep Ritz
formulations in Section~\ref{S:DEEP_RITZ}.

Building on our previous study of neural-network-based solvers for solitary wave computation \cite{HRW2026}, which considered PINNs, DeepONets, and FNOs, we extend the numerical comparison in this section by incorporating higher-order optimization strategies and Fourier-feature PINNs. 

The goal of this section is to evaluate 
the numerical performance in terms of accuracy and computation efficiency, where computation cost is measured using wall-clock elapsed time \texttt{time.perf\_counter()}. These baselines results later compared against the results obtained from the Deep Ritz formulations of Section~\ref{S:DEEP_RITZ}. 

\subsection{Physics-Informed Neural Network (PINN)}\label{S:PINN}
In the PINN framework \cite{RPE19}, we approximate the solitary-wave profile $Q(x)$ by the standard feedforward multilayer perceptron (MLP) $Q_\theta(x)$ with $m$ hidden layers and $n$ neurons in each layer, namely,
$$
Q(x)\approx Q_\theta(x),
$$
where $\theta$ denotes the collection of trainable network parameters. All spatial derivatives appearing in the loss are
evaluated using the reverse-mode automatic-differentiation engine \texttt{autograd} provided by \texttt{PyTorch} \cite[Sec.~4.3]{PaszkeEtAl2019}. Rather than approximating derivatives by finite differences, \texttt{autograd} records the elementary tensor operations executed during the forward evaluation of $Q_\theta(x)$ and applies the chain rule through the resulting computational graph, so that the computed derivatives are exact up to floating-point error \cite[Secs.~1 and~3]{BaydinEtAl2018}. The spatial derivatives $Q^\prime_\theta(x)$ and $Q^{\prime\prime}_\theta(x)$ are obtained by differentiating the network output with respect to the input $x$; computing $Q^{\prime\prime}_\theta(x)$ requires retaining the graph of the first derivative so that $Q^\prime_\theta(x)$ can itself be differentiated. The same construction subsequently allows the PINN loss, which depends on $Q_\theta$ and its spatial derivatives, to be differentiated with respect to the trainable parameters $\theta$.

In the numerical experiments, we test several {\it smooth activation functions}, including $\tanh$, SiLU, Gaussian, super-Gaussian, $\sech$, and the nonstandard activation
\begin{equation}\label{E:SiLUPLUS}
\sigma_{\mathrm{SiLUPLUS}}(x) =
\frac{x}{1+e^{-x^2}}.
\end{equation}
No single activation function is assumed to be optimal across all neural-network formulations. Instead, the activation is treated as a hyperparameter and selected separately for each method, boundary-condition formulation, and optimizer configuration. The activation chosen for each final setup is reported together with the remaining selected hyperparameters in the corresponding configuration tables.

For the single-power stationary problem, the network minimizes the residual associated with the finite-domain profile ODE \eqref{2ndODE}. Specifically, we define
\begin{equation}\label{E:residual}
\mathcal{R}(x;\theta)=Q^{\prime\prime}_\theta(x)-b Q_\theta(x)+\gamma Q_\theta^p(x).
\end{equation}
We evaluate this residual at collocation points $\{x_i\}_{i=1}^{N_f} \subset [0,L]$, where the computational domain is $[0,L]$. 
The boundary points $x=0$ and $x=L$ are treated separately via the Dirichlet boundary loss
\begin{equation}\label{E:pinn-loss}
    \mathcal{L}_{\mathrm{BC}}(\theta)=|Q_\theta(0)-Q_0|^2+|Q_\theta(L)|^2.
\end{equation}
Here, the residuals are scalar, so $|\cdot|$ denotes the usual absolute value. We write the total PINN loss as
\begin{equation}\label{E:L-1}
\mathcal{L}_{\mathrm{PINN}}(\theta)=\frac{1}{N_f}\sum^{N_f}_{i=1}|\mathcal{R}(x_i;\theta)|^2 + \lambda_{\mathrm{BC}}\left( |Q_\theta(0)-Q_0|^2 + |Q_\theta(L)|^2\right),
\end{equation}
where $\lambda_{\mathrm{BC}}$ is the boundary-weight parameter. The network parameters are then determined by solving
\begin{equation}\label{E:PINN-objective}
\theta^*=\arg\min_\theta\mathcal{L}_{\mathrm{PINN}}(\theta).
\end{equation}

Building on the PINN results that we obtained in \cite{HRW2026}, our present work focuses specifically on {\it higher-order} optimization strategies for training the PINN loss. In particular, we investigate the Gauss--Newton/ Levenberg--Marquardt (GN/LM) method and the quasi-Newton L--BFGS method. We also consider hybrid procedures, in which Adam is used only as an initialization or warm-up stage before refinement with GN/LM or L--BFGS. The general PINN algorithm with second-order optimizers is summarized in Algorithm~\ref{alg:PINN} below.

\subsubsection{Higher-order optimizer refinements for PINNs}\label{S:higher-order-PINNs}
The PINN loss in \eqref{E:L-1} has a nonlinear least-squares structure, that is, we can rewrite it as
\begin{equation}\label{eq:least-square-loss}
\mathcal{L}(\theta)=\|F(\theta)\|_{l^2}^2,
\end{equation}
where $F(\theta)$ is the finite-dimensional vector of rescaled residuals,
\begin{equation}\label{eq:residual-vector}
F(\theta)=
\begin{bmatrix}
N_f^{-1/2}\mathcal{R}(x_1;\theta)\\
\vdots\\
N_f^{-1/2}\mathcal{R}(x_{N_f};\theta)\\
\sqrt{\lambda_{\mathrm{BC}}}\bigl(Q_\theta(0)-Q_0\bigr)\\
\sqrt{\lambda_{\mathrm{BC}}}\,Q_\theta(L)
\end{bmatrix}.
\end{equation}
Thus, the first $N_f$ components of $F(\theta)$ are the interior ODE residuals evaluated at the collocation points and scaled by $N_f^{-1/2}$, while the final two components are the scaled Dirichlet boundary residuals. With this notation, $\|F(\theta)\|_{l^2}^2=\sum_{j=1}^{N_f+2}|F_j(\theta)|^2$, so \eqref{E:L-1} and \eqref{eq:least-square-loss}--\eqref{eq:residual-vector} describe the same loss.

Since $F(\theta)$ depends on $Q_\theta$ and its spatial derivatives, the parameter gradient of the loss must be propagated through the operations that construct $Q^\prime_\theta(x)$ and $Q^{\prime\prime}_\theta(x)$:
$$
\theta\longmapsto Q_\theta(x)\longmapsto(Q_\theta^\prime(x),\,Q^{\prime\prime}_\theta(x))\longmapsto F(\theta)\longmapsto\mathcal{L}(\theta).
$$
The gradient therefore carries the dependence of both the network values and their spatial derivatives on the weights and biases; this nested differentiation is handled by the same \texttt{PyTorch} automatic-differentiation graph described above \cite{BaydinEtAl2018,PaszkeEtAl2019}. It is this propagation through the spatial derivatives that makes the residual Jacobian $J_F$ below more expensive to form and, in general, more poorly conditioned than in an ordinary regression problem (a point we return to when introducing Levenberg--Marquardt damping in \eqref{eq:lm-update}).

Writing $\theta\in\R^{N_\theta}$, where $N_\theta$ is the total number of trainable weights and biases, let
$$
J_F(\theta)=\frac{\partial F(\theta)}{\partial\theta}\in\R^{(N_f+2)\times N_\theta}
$$
denote the Jacobian of the residual vector with respect to the trainable parameters. The chain rule applied to \eqref{eq:least-square-loss} gives
\begin{equation}\label{eq:least-squares-gradient}
\nabla_\theta\mathcal{L}(\theta)=2J_F(\theta)^T F(\theta)
=2\sum_{i=1}^{N_f+2}F_i(\theta)\,\nabla_\theta F_i(\theta),
\end{equation}
where the sum over residual components is exactly the matrix product $J_F(\theta)^T F(\theta)$.

This least-squares formulation motivates comparing optimizers that use different levels of derivative information. The three families considered here use it differently: Adam is a first-order adaptive-gradient method that uses $\nabla_\theta\mathcal{L}$ together with component-wise running estimates formed from present and past gradients; L--BFGS is a quasi-Newton method that infers curvature of the scalar loss from changes in successive parameter gradients; and GN/LM is a nonlinear least-squares method that uses the residual vector $F(\theta)$ directly to build a problem-specific curvature approximation. Throughout, ``higher-order'' is used broadly for methods that incorporate curvature information, or an approximation to it, rather than in the sense of forming exact second derivatives: L--BFGS estimates curvature from gradient differences, and GN/LM builds its curvature approximation from the residual Jacobian.

We use Adam as the baseline warm-up method, since it is a first-order adaptive-gradient method based on running estimates of the first and second moments of the gradient \cite{KingmaBa2015}. (Here, the ``second moment'' is a running average of squared gradient components, not a second derivative of the loss.) Adam is useful during the initial stage of training because it is relatively robust from random initialization and tolerates poor gradient scaling; however, as a first-order method, it may converge slowly once the iterate is already near a local minimizer. We therefore test two curvature-informed refinement strategies: L--BFGS and damped Gauss--Newton / Levenberg--Marquardt (GN/LM).
\smallskip

\paragraph{$\bullet$ \bf L--BFGS refinement.}
L--BFGS is a limited-memory quasi-Newton method for large-scale optimization \cite[Sec.~2]{LiuNocedal1989}, \cite[Secs.~6.1 and~7.2]{NocedalWright2006}. Rather than forming the $N_\theta\times N_\theta$ Hessian $\nabla_\theta^2\mathcal{L}$, which is impractical to store for large networks, it maintains an implicit approximation $H_k$ to the inverse Hessian built from a limited number of recent parameter and gradient differences, and produces search directions
\begin{equation}\label{eq:search-direction}
    p_k=-H_k\,\nabla_\theta\mathcal{L}(\theta_k).
\end{equation}
The parameters are then updated by $\theta_{k+1}=\theta_k+\alpha_k p_k$, where the step length $\alpha_k>0$ is selected by a strong Wolfe line search \cite{Wolfe1969,Wolfe1971,MoreThuente1994}. The strong Wolfe conditions enforce both a sufficient decrease in the loss and a curvature condition that keeps the BFGS inverse-Hessian approximation positive definite, and hence the search direction a descent direction; we use the strong Wolfe line search provided by the \texttt{PyTorch} implementation of L--BFGS. The derivation of \eqref{eq:search-direction} from the Newton direction and the precise Wolfe conditions are standard and are collected in Appendix~\ref{S:appendix-optimizers}. The stability of the strong Wolfe search is particularly helpful for the FF--PINNs of Section~\ref{S:FF--PINN}, where the Fourier-feature embedding makes the optimization landscape more sensitive to step size. Unlike GN/LM, however, L--BFGS is a general-purpose method for smooth scalar objectives and does not exploit the decomposition of the loss into the individual residuals in $F(\theta)$.
\smallskip

\paragraph{$\bullet$ \bf Gauss--Newton/Levenberg--Marquardt refinement.}
The GN/LM method is especially natural here because the PINN loss has the nonlinear least-squares form \eqref{eq:least-square-loss}. Linearizing $F$ near the current iterate $\theta_k$,
\begin{equation}\label{eq:F-linearize}
F(\theta_k+\delta)\approx F(\theta_k)+J_F(\theta_k)\delta,
\end{equation}
which follows by stacking the first-order component expansions $F_i(\theta_k+\delta)\approx F_i(\theta_k)+\nabla_\theta F_i(\theta_k)^T\delta$, whose gradient rows $\nabla_\theta F_i(\theta_k)^T$ are exactly the rows of $J_F(\theta_k)$. Gauss--Newton chooses the increment $\delta_k$ by minimizing the linearized residual norm $\|F(\theta_k)+J_F(\theta_k)\delta\|_{l^2}^2$, giving the normal equations
\begin{equation}\label{eq:gn-normal}
J_F(\theta_k)^T J_F(\theta_k)\,\delta_k
= -J_F(\theta_k)^T F(\theta_k),
\qquad \theta_{k+1}=\theta_k+\delta_k.
\end{equation}
Equivalently, $J_F(\theta_k)^T\bigl(F(\theta_k)+J_F(\theta_k)\delta_k\bigr)=0$, which states that the remaining linearized residual is orthogonal to the column space of $J_F(\theta_k)$. 
Geometrically, $J_F(\theta_k)\delta_k$ is chosen to approximate $-F(\theta_k)$ as closely as the first-order model allows: if exact cancellation is attainable the linearized residual is zero, and otherwise the step yields the smallest attainable $l^2$-norm. Since the components of $F(\theta)$ are the ODE and boundary residuals, $F(\theta)=0$ would mean the network satisfies the equation and boundary conditions exactly at all sampled points; thus zero is the ideal, and Gauss--Newton seeks the best locally attainable reduction toward it.

The connection to the Hessian follows by differentiating \eqref{eq:least-squares-gradient}:
 \begin{equation}\label{eq:least-square-exact-hessian}
     \nabla^2_\theta\mathcal{L}(\theta)=2J_F(\theta)^TJ_F(\theta)+2\sum^{N_f+2}_{i=1}F_i(\theta)\,\nabla^2_\theta F_i(\theta).
 \end{equation}
Gauss--Newton neglects the second term and uses $2J_F(\theta)^TJ_F(\theta)$ as a curvature approximation, which is especially appropriate near a solution where the residual components $F_i(\theta)$ are small. Unlike L--BFGS, which infers curvature from successive loss gradients, Gauss--Newton builds its approximation directly from the residual Jacobian.

In practice we use the damped Levenberg--Marquardt modification
\begin{equation}\label{eq:lm-update}
\left(J_F(\theta_k)^T J_F(\theta_k)+\lambda I\right)\delta_k
= -J_F(\theta_k)^T F(\theta_k),
\qquad \lambda>0,
\end{equation}
in which $\lambda I$ regularizes the linear system and improves stability when $J_F^T J_F$ is ill-conditioned or the linearization \eqref{eq:F-linearize} is inaccurate. For large $\lambda$ the increment approaches $\delta_k\approx-\tfrac{1}{\lambda}J_F(\theta_k)^T F(\theta_k)$, a scaled negative-gradient step by \eqref{eq:least-squares-gradient}, while as $\lambda\to0$ the update approaches the Gauss--Newton step; the undamped method is recovered at $\lambda=0$ \cite{Levenberg1944,Marquardt1963}. This damping is particularly relevant for PINN training because, as noted above, the residual depends on derivatives of the network output and the residual Jacobian can be poorly conditioned.

{\scriptsize
\begin{algorithm}[!htbp]
\footnotesize
\makeatletter
\algrenewcommand{\alglinenumber}[1]{\makebox[1.8em][r]{\scriptsize #1}\hspace{0.6em}}
\makeatother
\caption{PINN Method with Second-Order Optimizers for $Q''=b Q-\gamma Q^{p}$}
\begin{algorithmic}[1]
\State Choose domain $[0,L]$, collocation points $\{x_j\}_{j=1}^{N_f}$, and boundary condition type.
\State Define a neural network $q_\theta(x)$.
\State Define the PINN approximation by
$$
Q_\theta(x)=q_\theta(x).
$$
\State Define the ODE residual
$$
\mathcal R_\theta(x)
=
Q_\theta''(x)
-b Q_\theta(x)
+\gamma Q_\theta(x)^{p}.
$$
\State Define
$$
\mathcal{L}_{\mathrm{ODE}}(\theta)
=
\frac{1}{N_f}
\sum_{j=1}^{N_f}
\left|\mathcal R_\theta(x_j)\right|^2.
$$
\State Define the boundary loss
$$
\mathcal{L}_{\mathrm{BC}}(\theta) =\left|Q_\theta(0)-Q_0\right|^2+\left|Q_\theta(L)\right|^2,
$$
\State Define the total loss
$$
\mathcal{L}_{\mathrm{PINN}}(\theta)
=
\mathcal{L}_{\mathrm{ODE}}(\theta)
+
\lambda_{\mathrm{BC}}
\mathcal{L}_{\mathrm{BC}}(\theta).
$$
\State Choose one of the optimizer modes
$$
\mathsf{Opt}
\in
\{\text{Adam}, \text{ L--BFGS}, \text{ GN/LM}, 
\text{ Adam}\to\text{L--BFGS}, 
\text{ Adam}\to\text{GN/LM}\}.
$$
\For{$n=1,2,\dots,N_{\mathrm{steps}}$}
    \State Use automatic differentiation to compute $Q_\theta'$ and $Q_\theta''$.
    \State Compute $\mathcal{L}_{\mathrm{ODE}}(\theta)$, $\mathcal{L}_{\mathrm{BC}}(\theta)$, and $\mathcal{L}_{\mathrm{PINN}}(\theta)$.
    \State Update the trainable parameters $\theta$ according to the selected optimizer mode $\mathsf{Opt}$.
    \State Periodically resample the random subset of the collocation points and store checkpoint diagnostics.
\EndFor
\State Return the trained approximation $Q_\theta(x)$.
\end{algorithmic}
\label{alg:PINN}
\end{algorithm}
}

In summary, Adam uses first-order gradients with component-wise adaptive scaling but does not model curvature coupling between parameter directions; L--BFGS constructs an implicit inverse-Hessian approximation from parameter and gradient differences; and GN/LM exploits the nonlinear least-squares structure to form the problem-specific curvature approximation $J_F(\theta)^T J_F(\theta)$ directly from the residual Jacobian.

In addition to applying L--BFGS and GN/LM directly from the network initialization, we consider {\it hybrid} strategies in which an Adam warm-up precedes L--BFGS or GN/LM refinement. Such hybrid procedures are common in PINN practice and have been reported to improve accuracy and/or training efficiency relative to a single optimizer; see, e.g., \cite{LouMengKarniadakis2021,SahinVonDanwitzPopp2024,RathoreLeiFrangellaLuUdell2024}. Recent work indicates that stronger second-order optimizers can substantially outperform the classical Adam--L--BFGS procedure in demanding high-accuracy regimes \cite{JniniVellaZeinhofer2024,WangBennaniMartens2025}. When L--BFGS or GN/LM is applied directly from the network initialization, without an Adam warm-up, we refer to the resulting procedure as \emph{pure L--BFGS} or \emph{pure GN/LM}, respectively.

\subsubsection{Implementation of PINN with second order optimizers} 
In our implementation in the following sections, we record specific parameters and errors, mostly described in tables, such as Table \ref{tab:lbfgs_optimizer_schedule_sweep}. Throughout the numerical results, we report the standard \textit{relative} $L^2$ error and $L^\infty$ error, namely, 
$$
\frac{\|Q_\theta-Q_{\mathrm{exact}}\|_{L^2(0,L)}}{\|Q_{\mathrm{exact}}\|_{L^2(0,L)}} \quad \mbox{and} \quad \|Q_\theta-Q_{\mathrm{exact}}\|_{L^\infty(0,L)}.
$$

For the higher-order optimizer PINN comparison, the selected L--BFGS, Adam--L--BFGS, GN/LM, and Adam--GN/LM configurations from Tables~\ref{tab:selected_lbfgs_optimizer_configs} and \ref{tab:selected_gn_optimizer_configs} were evaluated with the Dirichlet boundary condition \eqref{E:Dirichlet} for the ODE $Q^{\prime\prime}=Q-Q^p$, with $p=3,\,5,$ and $7$. The results in Tables~\ref{tab:pinn_LBFGS_compared} and \ref{tab:pinn_GN_compare_D_explicit} show that all four refinement strategies are capable of reducing the training loss $\mathcal{L}(\theta)$ 
to small values, but the resulting profile accuracy depends strongly on the optimizer family.  
We use the term \emph{training loss} here because $\mathcal{L}(\theta)$ is the residual-based loss minimized by the PINN, whether the network is trained directly with Adam or further refinement using L--BFGS or GN/LM. More generally, \emph{training objective} is used when referring to both the PINN losses and the variational objectives minimized in the Deep Ritz methods as seen in Section~\ref{S:DEEP_RITZ}.

{\footnotesize
\begin{table}[!htbp]
\centering
\begin{tabular}{ll}
\toprule
\textbf{Parameter} & \textbf{Values tested / description} \\
\midrule
Domain length $L$ & $30$ \\

Pure L--BFGS runs
& L--BFGS initialized directly from the network initialization \\

Adam--L--BFGS runs
& Adam warm-up followed by L--BFGS refinement \\

Activation functions 
& $\tanh(x)$, $\mathrm{SiLU}(x)$, $\mathrm{SiLUPLUS}(x)$,
$e^{-x^2}$, $e^{-x^4}$, and $\sech(x)$ \\

Number of hidden layers & $3$, $4$, or $5$ \\

Neurons per layer & $32$, $64$, or $96$ \\

Number of collocation points $N_f$ & $2^8$, $2^9,\,2^{10},$ or $2^{11}$ \\

Adam warm-up length
& $2{,}500$, $5{,}000$, $7{,}500$, $10{,}000$, $12{,}500$, and $15{,}000$ epochs \\
Adam learning rate $lr$ &  $10^{-3}$ initial; StepLR decay \\

Boundary condition weight $\lambda_{\mathrm{BC}}$ & $0.1,
\, 1,\,5$, or $25$ \\

L--BFGS refinement steps & $25$, $50$, $75$, $100$, $125$, $150$, or $200$ \\

L--BFGS history size & $20,\,50$,  or $100$\\

L--BFGS max iterations per step & $1,\,5$, $10$, $20$, or $30$\\

L--BFGS learning rate & 0.5 and 1.0 \\

Selection criterion
& Smallest $L^2$ error, with runtime used to break near-ties \\

\bottomrule
\end{tabular}
\caption{Schedule quantities varied in the L--BFGS and Adam--L--BFGS optimizer tests.}
\label{tab:lbfgs_optimizer_schedule_sweep}
\end{table}
}

For these higher-order PINN refinement methods, hyperparameter selection was performed using the representative cubic problem $Q^{\prime\prime}=Q-Q^3$. We fixed $p=3$ and swept through the PINN, L--BFGS, and GN/LM optimizer parameters summarized in Tables~\ref{tab:lbfgs_optimizer_schedule_sweep} and \ref{tab:gn_optimizer_schedule_sweep}. 
After selecting the best configurations according to the smallest $L^2$ error, with runtime used to break near-ties, the same optimizer settings were then held fixed and evaluated on the cases $p=3,\,5,$ and $7.$ 
Although the cubic case was used during hyperparameter selection, the final selected configuration was rerun for $p=3$ on the primary computational device so that all three reported benchmark cases were obtained using the same finalized setup and timing environment. 
{We remark that when no exact profile is available, hyperparameters could instead be selected using the ODE and boundary residuals on an independent, finer validation grid, together with Pokhozhaev identities (or one of them, which is not algebraically enforced by the trial function). Then properties such as positivity, localization, centering of the computed profile under domain and quadrature refinement could also be checked. In this paper, we do not address this, and will consider this elsewhere, e.g., see an example of equations with combined potentials in \cite{HRR-cNLS}.}

The internal L--BFGS settings were fixed throughout our computations: strong-Wolfe line search, gradient tolerance $10^{-10}$, and change tolerance $10^{-12}$. These lower-level internal settings were not swept. The optimizer-schedule comparison varied the optimizer type, activation function, network architecture, number of collocation points, Adam warm-up length, and L--BFGS refinement schedule. For the checkpoint-based Adam--L--BFGS tests, Adam checkpoints were compared at different warm-up lengths, followed by L--BFGS refinement checkpoints. The selected configurations are therefore optimizer-dependent rather than common to all L--BFGS-based runs. For the pure L--BFGS case, the selected Dirichlet configuration used $50$ L--BFGS steps, history size $100$, a maximum of $10$ inner iterations per step, L--BFGS learning rate $1.0$, and $\lambda_{\mathrm{BC}}=5$. For the Adam--L--BFGS case, the selected Dirichlet configuration used a $15{,}000$-step Adam warm-up followed by $25$ L--BFGS refinement steps, history size $100$, and a maximum of $5$ inner iteration per step, L--BFGS learning rate $1.0$, and $\lambda_{\mathrm{BC}}=1$. In both cases, the strong-Wolfe line search was used for the L--BFGS updates. Table~\ref{tab:lbfgs_optimizer_schedule_sweep} summarizes the L--BFGS schedule quantities considered in the optimizer tests, while Table~\ref{tab:selected_lbfgs_optimizer_configs} reports the selected Dirichlet configurations for each L--BFGS-based optimizer.

\begin{table}[!htbp]
\centering
\resizebox{\textwidth}{!}{ 
\begin{tabular}{ccccccccccc}
\toprule
 \textbf{Optimizer}
& \textbf{Activation}
& \shortstack[c]{\textbf{Hidden}\\ \textbf{Layer}}
& \shortstack[c]{\textbf{Neurons}\\ \textbf{per Layer}}
& $N_f$
& \shortstack[c]{\textbf{Adam}\\ \textbf{Warm-up}}
& \shortstack[c]{\textbf{L--BFGS}\\ \textbf{Steps}}
& \shortstack[c]{\textbf{History}\\ \textbf{size}}
&  \shortstack[c]{\textbf{Max Iters}\\ \textbf{per Step}}
& \shortstack[c]{\textbf{L--BFGS}\\ \textbf{LR}}
& $\lambda_{\mathrm{BC}}$

\\
\midrule
 L--BFGS
& $\tanh(x)$
& $4$
& $32$
& $2^{10}$
& --
& $50$
& $100$
& $10$
& $1.0$
& $5.0$
\\

Adam--L--BFGS
& $\tanh(x)$
& $4$
& $64$
& $2^9$
& $15{,}000$
& $25$
& $100$ 
& $5$
& $1.0$
& $1.0$
\\

\bottomrule
\end{tabular}
}
\caption{Selected L--BFGS and Adam--L--BFGS for Dirichlet boundary conditions.}
\label{tab:selected_lbfgs_optimizer_configs}
\end{table}

\begin{table}[!htbp]
\centering
\setlength{\tabcolsep}{3pt}
\renewcommand{\arraystretch}{1.05}
\resizebox{\textwidth}{!}{ 
\begin{tabular}{|c||c|c|c|c||c|c|c|c|}
\hline
 & \multicolumn{4}{c||}{PINN with L--BFGS optimization} & \multicolumn{4}{c|}{PINN with Adam--L--BFGS optimization} \\
\hline\hline
$p$ & $\mathcal{L}_{\mathrm{PINN}}(\theta)$ & CPU time & $L^\infty$ error & $L^2$ error &  $\mathcal{L}_{\mathrm{PINN}}(\theta)$ & CPU time & $L^\infty$ error & $L^2$ error \\
\hline\hline
3 & 1.314525e-08 & 15.328056 & 7.162973e-03 & 5.864162e-03 & 5.917457e-08 & 353.818283 & 7.896724e-03 & 6.463914e-03
\\ \hline
5 & 1.679169e-07 & 14.411380 & 8.026411e-03 & 7.032712e-03 & 2.705050e-07 & 338.015580 & 5.668667e-03 & 4.915231e-03
\\ \hline
7 & 9.080274e-08 & 10.815290 & 2.407888e-03 & 2.187603e-03  & 2.613929e-07 & 326.965276 & 2.374574e-03 & 2.225127e-03
\\ \hline
\end{tabular} 
}
\caption{Comparison of PINN performance with L--BFGS and Adam--L--BFGS optimization and Dirichlet boundary conditions \eqref{E:Dirichlet} for $Q^{\prime\prime}=Q-Q^{p}$.}
\label{tab:pinn_LBFGS_compared}
\end{table}
In Table~\ref{tab:pinn_LBFGS_compared}, the left group of columns reports the pure L--BFGS results, while the right group reports the Adam--L--BFGS results.
For the L--BFGS-based methods, the pure L--BFGS runs are substantially cheaper than Adam--L--BFGS, and the added Adam warm-up does not uniformly improve the final accuracy. 
{ For $p=3$, the pure L--BFGS run gives the smaller $L^\infty$ and $L^2$ errors despite requiring much less CPU time. }

For $p=5$, Adam--L--BFGS gives the smaller $L^\infty$ and $L^2$ errors.
For $p=7$, Adam--L--BFGS gives a slightly smaller $L^\infty$ error, whereas pure L--BFGS gives a slightly smaller $L^2$ error, with both methods producing errors of the same order.
We therefore find that the Adam warm-up does not uniformly improve the L--BFGS refinement stage, even though it can help for the $p=5$ case. 

\begin{figure}[!htbp]
\centering
\includegraphics[width=0.49\textwidth, height=0.32\textwidth]{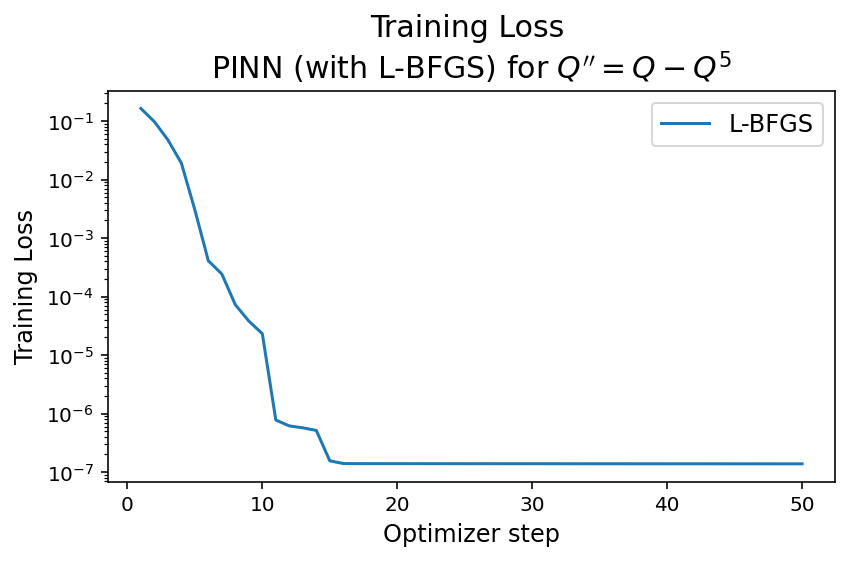}
\includegraphics[width=0.49\textwidth, height=0.32\textwidth]{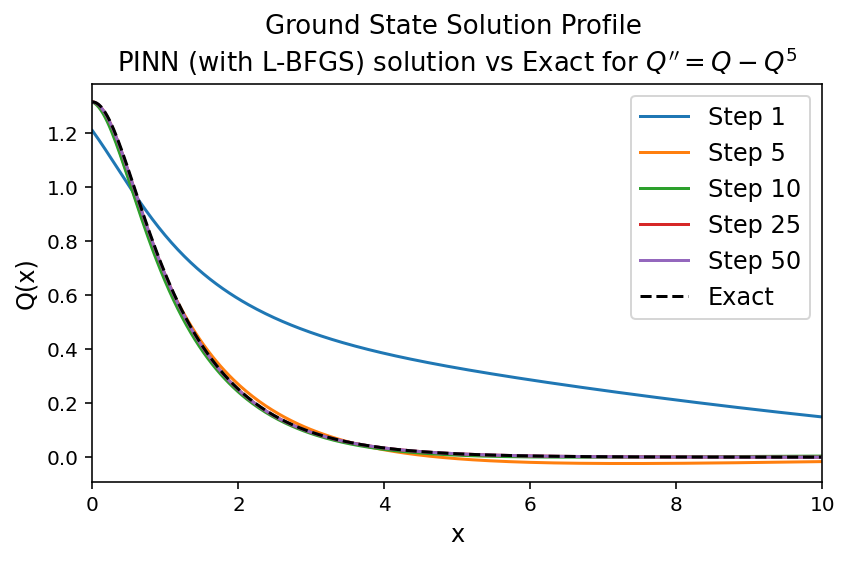}
\caption{Training loss (left) and ground state solution profile convergence (right) for L--BFGS PINN (with parameters as seen in Table~\ref{tab:selected_lbfgs_optimizer_configs}) for $Q^{\prime\prime}=Q-Q^5$. 
The (black dashed) curve labeled \emph{Exact} in the right panel is the explicit solution given by \eqref{Qexact}.}
\label{fig:LBFGS_converge}
\end{figure}
Figure~\ref{fig:LBFGS_converge} shows in the left panel that pure L--BFGS rapidly decrease the PINN training loss during the initial refinement steps before reaching a plateau. The profile plot in the right panel shows the 
learned profiles converging toward the exact ground state given by \eqref{Qexact} as the number of L--BFGS steps increases. We observe that the final profile after $50$ L--BFGS steps agrees closely with the exact profile, up to the errors reported in the left group of columns in Table~\ref{tab:pinn_LBFGS_compared}. For $p=5$, the final $L^\infty$ and $L^2$ errors are  $8.026411\times10^{-3}$ and $7.032712\times10^{-3}$, respectively, and are therefore both of order $10^{-3}$.
We show the computations for the representative quintic case in Figure \ref{fig:LBFGS_converge}, as it is the $L^2$-critical (borderline) case in one dimension, with similar convergence behavior observed for the cubic and septic cases.

\begin{figure}[!htbp]
\centering
\includegraphics[width=0.49\textwidth, height=0.32\textwidth]{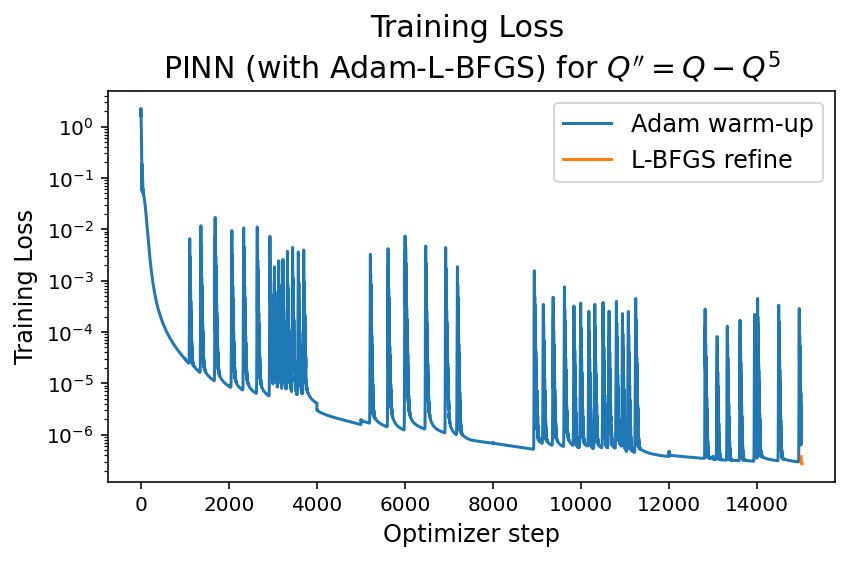}
\includegraphics[width=0.49\textwidth, height=0.32\textwidth]{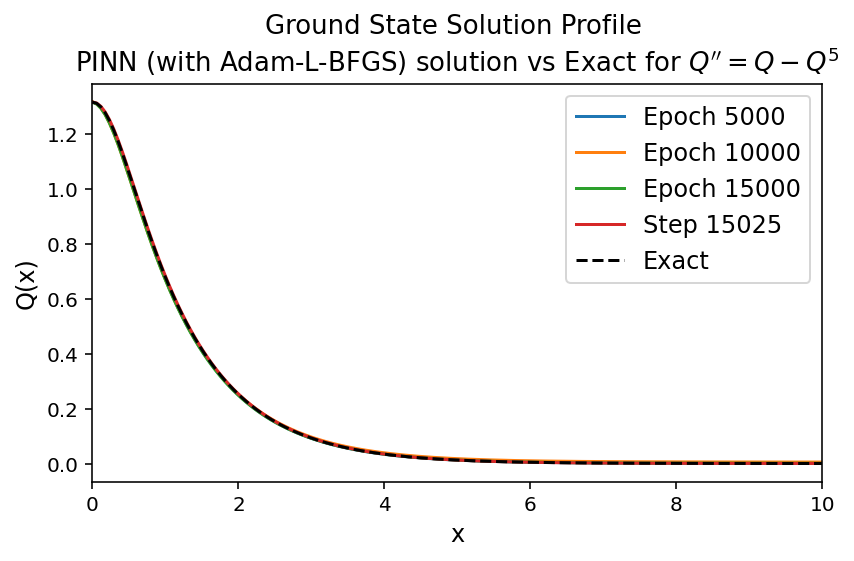}
\caption{Training loss (left) and ground state solution profile convergence (right) for Adam--L--BFGS PINN (with parameters as seen in Table~\ref{tab:selected_lbfgs_optimizer_configs}) for $Q^{\prime\prime}=Q-Q^5$. 
The (black dashed) curve labeled \emph{Exact} in the right panel is the explicit solution given by \eqref{Qexact}. The final learned profile corresponds to $15{,}000$ Adam steps followed by $25$ L--BFGS steps.}
\label{fig:adam_LBFGS_converge}
\end{figure}
Figure~\ref{fig:adam_LBFGS_converge} shows in the left panel that the Adam warm-up gradually reduces the training loss before the L--BFGS refinement stage produces the final decrease. In the profile plot in the right panel of the same figure, we observe that the learned profiles improve across the recorded Adam and L--BFGS checkpoints and approach the exact ground state given by \eqref{Qexact}. The labels of the form $15{,}000+j$ denote the total optimizer-step count, with the first $15{,}000$ epochs corresponding to Adam and the subsequent $j$ steps corresponding to the L--BFGS refinement. The final profile after the $15{,}000$-epoch Adam warm-up and $25$ L--BFGS steps agrees closely with the exact profile, up to the errors reported in the right group of columns of Table~\ref{tab:pinn_LBFGS_compared}. For $p=5$, the final $L^\infty$ and $L^2$ errors are $5.668667\times10^{-3}$ and $4.915231\times10^{-3}$, respectively, and are therefore both of order $10^{-3}$. Similar qualitative behavior was also observed for $p=3$ and $p=7$.

For the GN/LM PINN refinements, the optimizer comparison varied the optimizer type, activation function, network architecture, number of collocation points, Adam warm-up length, GN/LM refinement length, and damping parameters. For the pure GN/LM case, the method was initialized directly from the random network initialization, while in the Adam--GN/LM case, the GN/LM refinement was applied after an Adam warm-up. The damping-related quantities, including the initial damping, maximum number of backtracking steps, damping shrink factor, damping growth factor, and minimum damping, were swept as part of the optimizer tuning.

{\footnotesize
\begin{table}[!htbp]
\centering
\begin{tabular}{ll}
\toprule
\textbf{Parameter} & \textbf{Values tested / description} \\
\midrule
Domain length $L$ & $30$ \\
Pure GN/LM runs 
& GN/LM initialized directly from the network initialization \\

Adam--GN/LM runs 
& Adam warm-up followed by GN/LM refinement \\

Activation functions 
& $\tanh(x)$, $\mathrm{SiLU}(x)$, $\mathrm{SiLUPLUS}(x)$, 
$e^{-x^2}$, $e^{-x^4}$, and $\sech(x)$ \\

Number of hidden layers & $3$, $4$, or $5$ \\

Neurons per layer & $32$, $64$, or $96$ \\

Number of collocation points $N_f$ & $2^8$, $2^9,\,2^{10},$ or $2^{11}$ \\

Adam warm-up checkpoints 
& $2{,}500$, $5{,}000$, $7{,}500$, $10{,}000$, $12{,}500$, and $15{,}000$ epochs \\

Adam learning rate $lr$ &  $10^{-3}$ initial; StepLR decay \\

Boundary condition weight $\lambda_{\mathrm{BC}}$ & $0.1,
\, 1.\,5$ or $25$ \\

GN/LM refinement steps & $25$, $50$, $75$, $100$, $125$, $150$, or $200$ \\

Initial damping 
& $5\times 10^{-2}$, $10^{-3}$, and $5\times 10^{-4}$ \\

Maximum backtracking steps 
& $2$, $4$, 6, 8, and 10 \\

Damping shrink factor 
& Values in $[0.4,1.1]$ \\

Damping growth factor 
& $5$, $10$, and $15$ \\

Minimum damping 
& $10^{-12}$, $10^{-14}$, and $10^{-16}$ \\

Selection criterion 
& Smallest $L^2$ error, with runtime time used to break near-ties \\

\bottomrule
\end{tabular}
\caption{Schedule and damping quantities varied in the GN/LM and Adam--GN/LM optimizer tests.}
\label{tab:gn_optimizer_schedule_sweep}
\end{table}
}

The selected GN/LM-based Dirichlet configurations shared the collocation size, refinement length, damping parameters, and boundary-condition weight, { but not the activation function}. In both cases, the selected runs used $N_f=2^9$, $150$ GN/LM refinement steps, initial damping $10^{-3}$, $2$ maximum backtracking steps, damping shrink factor $0.6$, damping growth factor $10$, and $\lambda_{\mathrm{BC}}=5$. 
{Pure GN/LM was started directly from the random network initialization and used $\sech(x)$ activation with  $4$ hidden layers with $64$ neurons per layer. 
Adam--GN/LM used $\tanh(x)$ activation with with $5$ hidden layers and $32$ neurons per layer, following a $15{,}000$-epoch Adam warm-up.} 
Table~\ref{tab:gn_optimizer_schedule_sweep} summarizes { the tuning ranges, and 
Table~\ref{tab:selected_gn_optimizer_configs} reports the selected configurations.
}

{\footnotesize
\begin{table}[!htbp]
\centering
\resizebox{\textwidth}{!}{ 
\begin{tabular}{cccccccccccc}
\toprule
\textbf{Optimizer} 
& \textbf{Activation}
& \shortstack[c]{\textbf{Hidden}\\ \textbf{Layer}}
& \shortstack[c]{\textbf{Neurons}\\ \textbf{per Layer}}
& $N_f$
& \shortstack[c]{\textbf{Adam}\\ \textbf{Warm-up}}
& \shortstack[c]{\textbf{GN/LM}\\ \textbf{Steps}}
& \shortstack[c]{\textbf{Max}\\ \textbf{Backtracks}}
& \shortstack[c]{\textbf{Damping}\\ \textbf{Shrink}}
& \shortstack[c]{\textbf{Damping}\\ \textbf{Growth}}
& \shortstack[c]{\textbf{Initial}\\ \textbf{Damping}}
& \textbf{$\lambda_{\mathrm{BC}}$}
\\
\midrule

 GN/LM 
& $\sech(x)$ 
& $4$
& $64$
& $2^9$
& -- 
& $150$ 
& $2$
& $0.6$
& $10$
& $10^{-3}$
& $5$
\\

Adam--GN/LM 
& $\tanh(x)$ 
& $5$
& $32$
& $2^9$
& $15{,}000$ 
& $150$ 
& $2$
& $0.6$
& $10$
& $10^{-3}$
& $5$
\\

\bottomrule
\end{tabular}
}
\caption{Selected GN/LM and Adam--GN/LM configurations for Dirichlet boundary conditions.}
\label{tab:selected_gn_optimizer_configs}
\end{table}
}

\begin{figure}[!htbp]
\centering
\includegraphics[width=0.49\textwidth, height=0.32\textwidth]{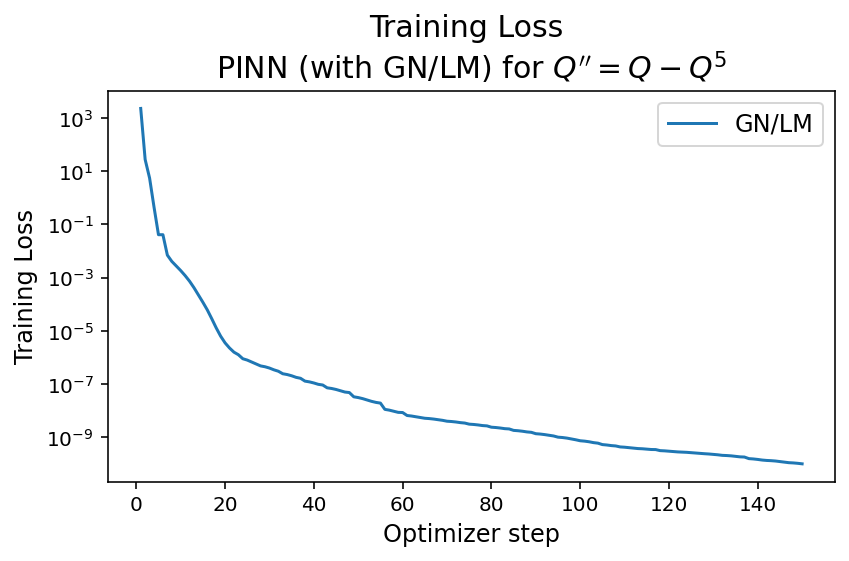}
\includegraphics[width=0.49\textwidth, height=0.32\textwidth]{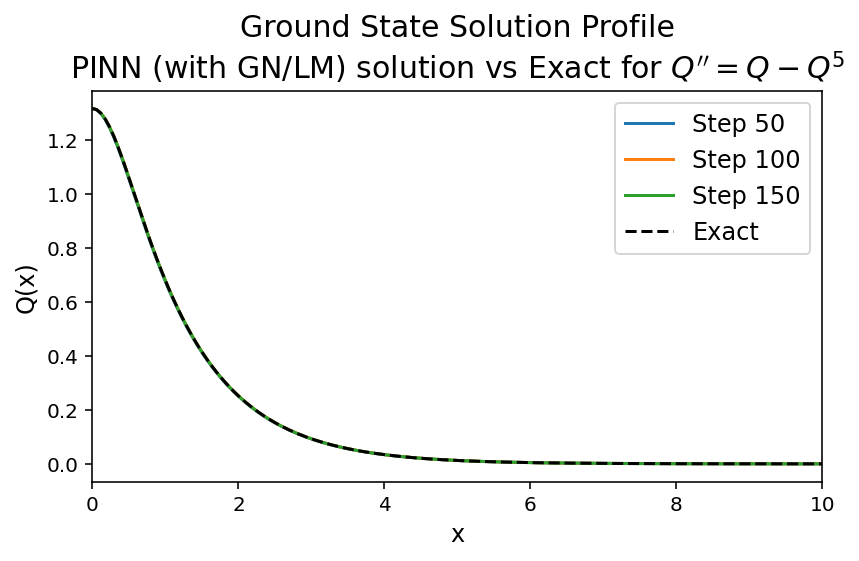}
\caption{Training loss (left) and ground state solution profile convergence (right) for GN/LM PINN (with parameters as seen in Table~\ref{tab:selected_gn_optimizer_configs}) for $Q^{\prime\prime}=Q-Q^5$. 
The (black dashed) curve labeled \emph{Exact} in the right panel is the explicit solution given by \eqref{Qexact}.}
\label{fig:GN_converge}
\end{figure}

\begin{figure}[!htbp]
\centering
\includegraphics[width=0.49\textwidth, height=0.32\textwidth]{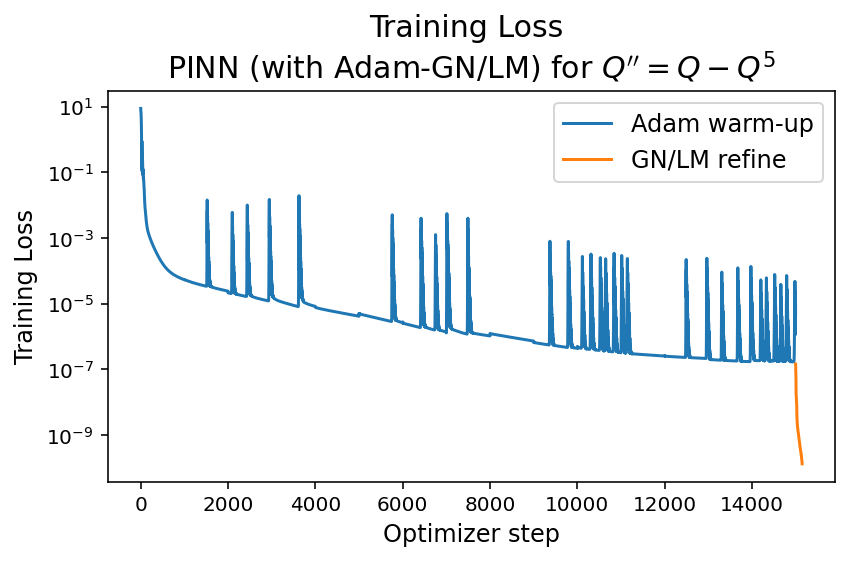}
\includegraphics[width=0.49\textwidth, height=0.32\textwidth]{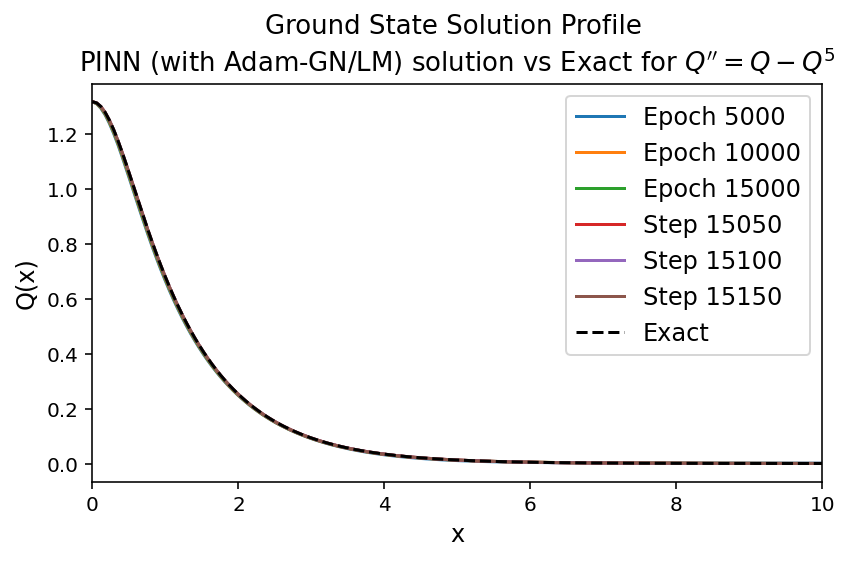}
\caption{Training loss (left) and ground state solution profile convergence (right) for Adam--GN/LM PINN (with parameters as seen in Table~\ref{tab:selected_gn_optimizer_configs}) for $Q^{\prime\prime}=Q-Q^5$. 
The (black dashed) curve labeled \emph{Exact} in the right panel is the explicit solution given by \eqref{Qexact}. The final learned profile corresponds to $15{,}000$ Adam steps followed by $150$ GN/LM steps.} 
\label{fig:adam_GN_converge}
\end{figure}

The GN/LM-based methods give smaller errors than the L--BFGS-based methods in most cases, but at a much higher computational cost when run on the primary CPU device. In particular, pure GN/LM obtains errors near $10^{-4}$--$10^{-5}$ for $p=3$ and $p=5$, while Adam--GN/LM gives the best result among the PINN optimizers for $p=7$. The corresponding convergence plots in Figures~\ref{fig:GN_converge} and \ref{fig:adam_GN_converge} show a rapid reduction of the training objective  under GN/LM refinement and close agreement between the computed profile, that the neural network learned, and the exact ground state given by \eqref{Qexact}. Overall, we find that these results suggest GN/LM refinement is the most accurate PINN refinement mechanism tested here, while the Adam warm-up can improve robustness for the more difficult higher-power case.

\begin{table}[!htbp]
\centering
\setlength{\tabcolsep}{3pt}
\renewcommand{\arraystretch}{1.05}
\resizebox{\textwidth}{!}{%
\begin{tabular}{|c||c|c|c|c||c|c|c|c|}
\hline
 & \multicolumn{4}{c||}{PINN with GN/LM optimization} & \multicolumn{4}{c|}{PINN with Adam--GN/LM optimization} \\
\hline\hline
$p$ & $\mathcal{L}_{\mathrm{PINN}}(\theta)$ & CPU time & $L^\infty$ error & $L^2$ error &  $\mathcal{L}_{\mathrm{PINN}}(\theta)$ & CPU time & $L^\infty$ error & $L^2$ error \\
\hline\hline
3 & 2.763789e-11 & 4532.584975 & 1.193728e-04 & 9.808560e-05 & 6.763931e-13 & 1639.024311 & 1.417711e-04 & 1.158408e-04
\\ \hline
5 & 1.015142e-10 & 4940.055658 & 6.839973e-05 & 5.982481e-05 & 1.328934e-10 & 1872.362646 & 5.831789e-04 & 5.065382e-04
\\ \hline
7 & 3.184690e-08 & 5205.814115 & 2.681668e-03 & 2.388390e-03 & 4.723684e-10 & 1812.308855 & 9.558302e-04 & 8.583382e-04
\\ \hline
\end{tabular}%
}

\caption{Comparison of PINN performance with GN/LM and Adam--GN/LM optimization and Dirichlet boundary conditions \eqref{E:Dirichlet} for $Q^{\prime\prime}=Q-Q^{p}$.}
\label{tab:pinn_GN_compare_D_explicit}
\end{table}

In Table~\ref{tab:pinn_GN_compare_D_explicit}, the left group of columns reports the pure GN/LM results, while the right group reports the Adam--GN/LM results. We observe that the two methods give errors of comparable order for $p=3$, while pure GN/LM gives the more accurate profile for $p=5$. In particular, the pure GN/LM errors for the quintic case are of order $10^{-5}$, compared with order $10^{-4}$ for Adam--GN/LM. For $p=7$, however, the Adam warm-up produces the smaller errors. Thus, we find that, although pure GN/LM attains the smallest error for the cubic and quintic cases, Adam--GN/LM is more robust across the increasing nonlinear test cases.

Figure~\ref{fig:GN_converge} shows in the left panel that pure GN/LM produces a sharp and sustained decrease in the PINN training loss over the refinement steps. In the right panel, we observe that the successive learned profiles approach the exact ground state given by \eqref{Qexact}. By step $50$, the learned and exact curves are already visually indistinguishable at the plotting scale. We nevertheless continue the optimization to the selected final value of $150$ GN/LM steps  and report the errors at that final step in the left group of columns of Table~\ref{tab:pinn_GN_compare_D_explicit} because the order of error at step $50$ is only $10^{-3}$, compared with $10^{-5}$ at step $150$. For $p=5$, the final $L^\infty$ and $L^2$ errors are $6.839973\times10^{-5}$ and $5.982481\times10^{-5}$, respectively, and are therefore both of order $10^{-5}$. The error at step $50$ is not reported because this curve is an intermediate visualization checkpoint rather than the final trained configuration used in the comparison table, although the secondary-device results in Section~\ref{S:secondaryPINN} show that pure GN/LM can be sensitive to damping and branch selection.

Figure~\ref{fig:adam_GN_converge} shows that Adam--GN/LM combines the gradual stabilization provided by Adam 
with a rapid final decrease in the training loss during the
GN/LM refinement. In the right panel, we compare the profiles at the recorded Adam and GN/LM checkpoints with the exact ground state given by \eqref{Qexact}. The labels of the form $15000+j$ denote the total optimizer-step count, with the first $15{,}000$ steps corresponding to the Adam warm-up and the subsequent $j$ steps corresponding to GN/LM refinement. We retain only a subset of the available checkpoints in the legend, with no more than eight displayed curves, so that the plot remains readable. Thus, the displayed refinement-step values need not be equally spaced. We observe that the final profile after $150$ GN/LM refinement steps agrees closely with the exact profile, up to the errors reported in the right group of columns of Table~\ref{tab:pinn_GN_compare_D_explicit}. For $p=5$, the final $L^\infty$ and $L^2$ errors are $5.831789\times10^{-4}$ and $5.065382\times10^{-4}$, respectively, thus, both are of order $10^{-4}$, supporting Adam--GN/LM as the most robust PINN refinement strategy (among the tested methods). Although pure GN/LM gives the smaller error for the quintic case, we regard Adam--GN/LM as more robust among the cases tested, since Adam warm-up provides a more reliable initialization for the GN/LM refinement stage.

\subsubsection{Comparison with the ``secondary” computational device for PINNs with higher-order optimizer refinements}\label{S:secondaryPINN}
As a secondary-device check, we repeated the selected higher-order PINN refinement experiments on the separate GPU-capable workstation described in Section~\ref{S:hardware}. The purpose of this comparison is not to introduce a new hyperparameter sweep, but to evaluate the same selected PINN configurations under a different hardware and linear algebra environment. Therefore, the architecture, collocation size, boundary-condition weight, optimizer schedule, and stopping criteria were kept fixed, while only the computational device was changed. The larger-damping GN/LM experiments considered below are reported separately as a sensitivity test and are not part of this fixed-configuration device comparison. As in the preceding subsection, we may refer to L--BFGS or GN/LM applied directly from the network initialization, without an Adam warm-up, as \textit{pure L--BFGS} or \textit{pure GN/LM}, respectively.

As explained in Section~\ref{S:hardware}, these timings are reported as GPU+CPU rather than GPU-only, and are not directly comparable to CPU-only timings as a measure of processor performance alone.

\begin{table}[!htbp]
\centering
\setlength{\tabcolsep}{3pt}
\renewcommand{\arraystretch}{1.05}
\resizebox{\textwidth}{!}{%
\begin{tabular}{|c||c|c|c|c||c|c|c|c|}
\hline
 & \multicolumn{4}{c||}{PINN with L--BFGS optimization} & \multicolumn{4}{c|}{PINN with Adam--L--BFGS optimization} \\
\hline\hline
$p$ & $\mathcal{L}_{\mathrm{PINN}}(\theta)$ & GPU+CPU time & $L^\infty$ error & $L^2$ error &  $\mathcal{L}_{\mathrm{PINN}}(\theta)$ & GPU+CPU time & $L^\infty$ error & $L^2$ error \\ \hline
\hline
3 & 5.329101e-08 & 7.505031 & 6.583645e-03 & 5.421850e-03 & 3.593178e-08 & 62.361742 & 3.780430e-03 & 3.079454e-03 \\ \hline

5 & 7.202707e-08 & 6.157369 & 5.845897e-03 & 5.090068e-03 & 1.164522e-07 & 63.238091 & 2.992726e-03 & 2.691600e-03 \\ \hline

7 & 4.067500e-08 & 6.788750 & 1.609846e-03 & 1.493517e-03 & 2.353134e-07 & 61.959399 & 9.146505e-04 & 1.052700e-03 \\ 
\hline

\end{tabular}%
}
\caption{Secondary-device GPU+CPU comparison of PINN performance with L--BFGS and Adam--L--BFGS optimization and Dirichlet boundary conditions \eqref{E:Dirichlet} for $Q^{\prime\prime}=Q-Q^{p}$.}
\label{tab:pinn_LBFGS_compare_D_explicit_secondary}
\end{table}

Tables~\ref{tab:pinn_LBFGS_compare_D_explicit_secondary} and \ref{tab:pinn_GN_compare_D_explicit_secondary-b} show that the secondary GPU+CPU runs substantially reduce the reported wall-clock time for the higher-order refinement experiments. This reduction is especially clear for the Adam--L--BFGS and GN/LM-based methods, where the secondary device gives much shorter run times than the corresponding primary-device CPU runs. The L--BFGS-based methods remain broadly comparable to the primary-device CPU runs, with errors of the same general order. The Adam--GN/LM runs are also stable across the secondary-device comparison and remain close to the exact ground state given by \eqref{Qexact} for $p=3,\,5,$ and $7$.

In Table~\ref{tab:pinn_LBFGS_compare_D_explicit_secondary}, the left group of columns reports the pure L--BFGS results, while the right group reports the Adam--L--BFGS results. We observe that all of the profile errors are of order $10^{-3}$. On the secondary device, the Adam warm-up reduces the profile errors for each value of $p$, but it also increases the GPU+CPU time from $6-8$ seconds for pure L--BFGS to approximately $62-63$ seconds for Adam--L--BFGS. Thus, we find the same general trade-off observed on the primary device: the Adam warm-up can improve the final profile, but at a substantially greater computational cost.

The main discrepancy in the secondary-device PINN runs occurs for the pure GN/LM optimizer. With the selected initial damping used in Table~\ref{tab:selected_gn_optimizer_configs}, the septic case $p=7$ does not converge to the centered ground-state branch. Instead (although the profile preserves the prescribed value $Q(0)$), the GN/LM iterates develop an additional bump further to the right. This produces the large $L^\infty$ and $L^2$ errors reported for the $p=7$ entry, even though the training loss is still reduced. Hence, we have to distinguish between reduction of the sampled residual-based loss and convergence to the desired ground-state profile.

\begin{table}[!htbp]
\centering
\setlength{\tabcolsep}{3pt}
\renewcommand{\arraystretch}{1.05}
\resizebox{\textwidth}{!}{%
\begin{tabular}{|c||c|c|c|c||c|c|c|c|}
\hline
 & \multicolumn{4}{c||}{PINN with GN/LM optimization} & \multicolumn{4}{c|}{PINN with Adam--GN/LM optimization} \\
\hline\hline
$p$ & $\mathcal{L}_{\mathrm{PINN}}(\theta)$ & GPU+CPU time & $L^\infty$ error & $L^2$ error &  $\mathcal{L}_{\mathrm{PINN}}(\theta)$ & GPU+CPU time & $L^\infty$ error & $L^2$ error \\ \hline
\hline
3 & 4.025347e-11 & 538.109242 & 9.699502e-05 & 7.902533e-05 & 1.294946e-12 & 506.161066 & 5.594483e-05 & 4.555836e-05 \\ \hline
5 & 4.330355e-11 & 534.144825 & 2.525764e-04 & 2.192186e-04 & 1.701413e-10 & 511.624151 & 5.432745e-04 & 4.717845e-04 
\\ \hline
7 & 5.115845e-06 & 538.628306 & 1.256511 & 1.411610 & 1.485399e-10 & 531.554408 & 6.169997e-04 & 5.539982e-04
\\ \hline

\end{tabular}%
}
\caption{Secondary-device GPU+CPU comparison of PINN performance with GN/LM and Adam--GN/LM optimization and Dirichlet boundary conditions \eqref{E:Dirichlet} for $Q^{\prime\prime}=Q-Q^{p}$.}
\label{tab:pinn_GN_compare_D_explicit_secondary-b}
\end{table}

In Table~\ref{tab:pinn_GN_compare_D_explicit_secondary-b}, the left group of columns reports the pure GN/LM results, while the right group reports the Adam--GN/LM results. For $p=3$, both methods give errors between orders of $10^{-4}$ to $10^{-5}$, while for $p=5$ both give errors of order $10^{-4}$. The most significant difference occurs for $p=7$: pure GN/LM profile has errors of order one because it converges to the incorrect branch, whereas Adam--GN/LM gives $L^\infty$ and $L^2$ errors of $6.169997\times 10^{-4}$ and $5.539982\times 10^{-4}$, respectively, both of order $10^{-4}$. Thus, one can note that the Adam warm-up prevents the branch-selection failure in the septic case without increasing the secondary-device run time relative to pure GN/LM.

\begin{table}[!htbp]
\centering
\setlength{\tabcolsep}{3pt}
\renewcommand{\arraystretch}{1.05}
\begin{tabular}{|c||c|c|c|c|}
\hline
 & \multicolumn{4}{c|}{PINN with GN/LM optimization} \\
\hline\hline
$p$ & $\mathcal{L}_{\mathrm{PINN}}(\theta)$ & GPU+CPU time & $L^\infty$ error & $L^2$ error \\ \hline
\hline
3 & 2.483710e-08 & 591.697312 & 1.413589e+00 & 1.414096e+00 \\ \hline
5 & 1.058034e-06 & 586.388381 & 1.318001e+00 & 1.417144e+00 \\ \hline
7 & 1.221674e-09 & 574.305437 & 6.810766e-04 & 6.118432e-04 \\ \hline

\end{tabular}%
\caption{Secondary-device GPU+CPU comparison of PINN performance with GN/LM with larger initial damping ($10^{-1}$) and Dirichlet boundary conditions \eqref{E:Dirichlet} for $Q^{\prime\prime}=Q-Q^{p}$.}
\label{tab:pinn_GN_compare_D_explicit_secondary2}
\end{table}

\begin{figure}[!htbp]
\centering
\includegraphics[width=0.49\textwidth, height=0.32\textwidth]{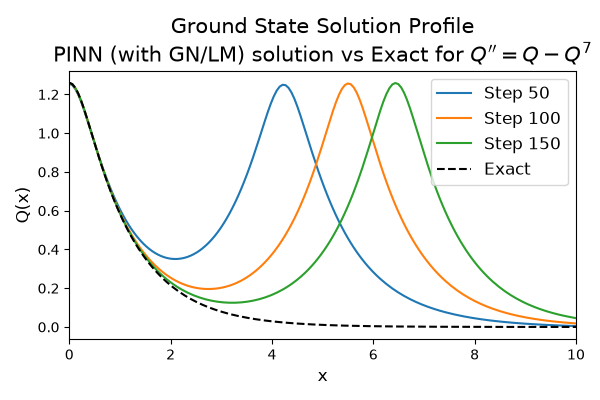}
\includegraphics[width=0.49\textwidth, height=0.32\textwidth]{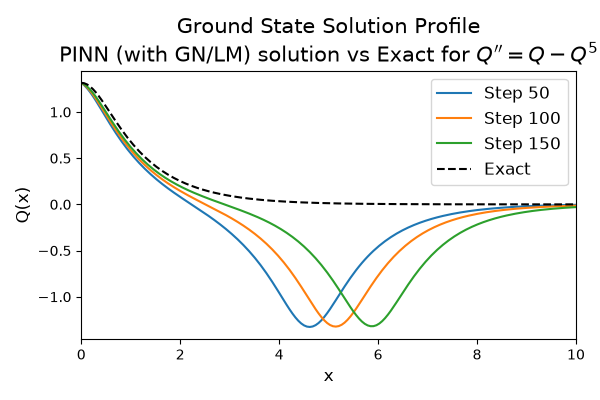}
\caption{Secondary device ground state solution profile convergence for GN/LM PINN with $p=7$ with parameters as seen in Table~\ref{tab:selected_gn_optimizer_configs} (left) and $p=5$ with GN/LM initial damping $10^{-1}$ (right) for $Q^{\prime\prime}=Q-Q^p$. The (black dashed) curve labeled \emph{Exact} in the right panel is the explicit solution given by \eqref{Qexact}.}
\label{fig:GN_converge_bad}
\end{figure}

To test whether this failure was caused by insufficient damping, we repeated the pure GN/LM runs with larger initial damping  values $10^{-2}$ and $10^{-1}$. Table~\ref{tab:pinn_GN_compare_D_explicit_secondary2} reports the results obtained with initial damping $10^{-1}$. This change corrects the septic case and yields an accurate $p=7$ profile. We note, however, that the same damping choice destabilizes the lower-power cases. For $p=3$ the learned solution develops a secondary bump away from the origin, while for $p=5$ the optimizer converges to a negative, off-centered profile rather than the positive, centered at $x=0$ ground state. We show representative non-ground state profiles in Figure~\ref{fig:GN_converge_bad}. We note that increasing the damping resolves the branch-selection issue for the septic (supercritical) case but introduces the same type of non-ground state convergence for the cubic and quintic cases. 

From Table~\ref{tab:pinn_GN_compare_D_explicit_secondary2}, we observe that the larger damping reverses the behavior seen with the selected damping. The cubic and quintic profiles now have $L^\infty$ and $L^2$ errors of order one, while the septic profile has errors of order $10^{-4}$. Thus, no single one of these two initial damping choices reliably selects the desired branch for all three powers when pure GN/LM is initialized directly from the random network parameters.

A plausible explanation is that the pure GN/LM PINN loss does not uniquely select the ground-state branch. We interpret the observed failure as evidence of the branch-selection ambiguity. In the Dirichlet setting, the loss enforces the residual together with $Q(0)=Q_0$ and $Q(L)=0$, but it does not directly impose positivity, monotone decay, or $Q^\prime(0)=0$. As a result, the network can satisfy the boundary conditions and maintain a small residual at the sampled collocation points while developing a secondary interior bump, producing a non-ground-state profile.

This issue is more pronounced for pure GN/LM because it starts from the random initialization and takes large linearization-based correct steps in a nonconvex least-squares landscape. Small numerical differences between the devices or linear algebra back-ends may therefore alter the accepted GN/LM steps and steer the iteration toward a different local minimizer. Thus, the secondary device double-bump profile does not necessarily indicate a change in the equation or code; rather, we interpret it as showing that pure GN/LM is branch-sensitive unless the desired ground-state structure is enforced or reached through a more stable path such as Adam--GN/LM.
In the left panel of Figure~\ref{fig:GN_converge_bad}, we compare the pure GN/LM profile for $p=7$, obtained using the selected damping from Table~\ref{tab:gn_optimizer_schedule_sweep}, with the exact ground state given by \eqref{Qexact}. The final learned profile develops an additional positive bump and does not align with the exact profile; the corresponding $L^\infty$ and $L^2$ errors in the left group of columns of Table~\ref{tab:pinn_GN_compare_D_explicit_secondary-b} are 1.256511 and 1.411610, respectively, and are therefore of order one. In the right panel, we compare the $p=5$ profile obtained with the larger initial damping $10^{-1}$ with the exact solution from \eqref{Qexact}. The final learned solution profile is negative and off-centered rather than positive and centered at the origin, and its $L^\infty$ and $L^2$ errors from Table~\ref{tab:pinn_GN_compare_D_explicit_secondary2} are 1.318001 and 1.417144, respectively, again of order one.

In conclusion, we find that pure GN/LM is highly sensitive to damping, initialization, and device-dependent numerical details. 
We observe that this is indeed the case when testing higher order powers in the nonlinear term, where changing the initial damping yields results similar to those in Figure \ref{fig:GN_converge_bad}. 
In contrast, Adam--GN/LM is more robust because the Adam warm-up moves the network closer to the desired ground-state branch before the GN/LM refinement begins. We emphasize that our secondary-device comparison supports the main conclusion that Adam--GN/LM is the more reliable PINN refinement strategy, while pure GN/LM can achieve high accuracy only when the damping and initialization place the iteration on the correct solution branch. 

\subsection{Fourier-Feature Physics-Informed Neural Network (FF--PINN)}\label{S:FF--PINN}
We also test a Fourier-feature physics-informed neural network (FF--PINN), in which the scalar coordinate input is first embedded into a trigonometric feature space before being passed to the multilayer perceptron. In our implementation, the physical coordinate is scaled by the computational length $L$. When the raw coordinate is included, the embedding is
$$
x\mapsto \left[ \frac{x}{L},\,
    \sin\!\left(2\pi \frac{x}{L}B\right),\,
    \cos\!\left(2\pi \frac{x}{L}B\right)
  \right],
$$
where $B\in\R^{d\times N_{\mathrm{FF}}}$ is a {\it frequency matrix} in $d$ dimensions. In our 1D case, $B$ should be referred to as a {\it frequency vector} as we are studying the one-dimensional ground state problem. %

We also test the corresponding embedding without the raw coordinate, 
$$
x \mapsto   \left[  \sin\! \left(2\pi \frac{x}{L}B\right),\, \cos\!\left(2\pi \frac{x}{L}B\right)  \right].
$$

The entries of $B$ are initialized from a Gaussian distribution whose scale is determined by $\sigma_{\mathrm{FF}}$. We consider both fixed-frequency embeddings, where $B$ is kept fixed after initialization and trainable-frequency embeddings, where $B$ is optimized together with the neural-network weights and biases.
This construction is related to random Fourier features for kernel approximation \cite{RahimiRecht2007} and to Fourier-feature embeddings used to reduce the spectral bias of coordinate-based neural networks \cite{Tancik2020}. Fourier and multiscale feature embeddings have also been studied in the PINN setting for problems where standard PINNs may struggle to represent high-frequency or multiscale components \cite{WangWangPerdikaris2021}.
Thus, inclusion of the raw coordinate $x/L$ is treated as a hyperparameter in the FF--PINN study.

The role of this embedding is to provide the network with oscillatory building blocks before training the MLP. With fixed random frequencies, one obtains a richer input representation without introducing additional trainable frequency parameters. More generally, one could instead treat $B$ as trainable, so that the optimization adjusts the frequencies themselves rather than only the network weights. This may be advantageous when the relevant frequency scale is not known in advance, since the model can adapt its embedding to the problem. However, it also makes the optimization more delicate, because the network must then learn both the shape of the solution and the frequencies used to represent it. In particular, trainable frequencies may drift toward unnecessarily large values, which can produce oscillatory artifacts and reduce stability.

This issue is especially important because the residual contains the second derivative $Q^{\prime\prime}$. Indeed, a high-frequency component of the form $\sin(kx)$ satisfies
$$
\frac{d^2}{dx^2}\sin(kx)=-k^2\sin(kx),
$$
so high frequencies are amplified in the differential operator through the factor $k^2$. Consequently, overly large Fourier frequencies can have a disproportionate effect on the residual, making the optimization more sensitive and potentially favoring oscillatory low-residual profiles that do not correspond to the desired ground state. For this reason, the Fourier-feature scale and the number of features must be chosen with care.

The motivation for testing this architecture is that our target equation is a stationary profile problem, i.e., $Q^{\prime\prime}=b Q-\gamma Q^{p}$,
and Fourier-based classical methods, such as Petviashvili iterations combined with FFT discretizations, perform very effectively for localized solitary-wave profiles. The FF--PINN is therefore used to test whether enriching the input representation with sinusoidal features improves the PINN approximation relative to the standard coordinate-input architecture. At the same time, this advantage is not automatic in the present setting: the desired ground state is smooth, positive, centered, and localized, rather than highly oscillatory. Thus, Fourier features increase the expressive power of the network, but they can also introduce
unnecessary oscillations.

For this reason, we pair the FF--PINN experiments with  ground-state-oriented constraints incorporated directly into the neural network trial function. To target the nontrivial ground-state branch, we use the structured trial function

\begin{equation}\label{E:constrained_ansatz}
{
Q_\theta(x)= Q_0 \, \frac{\left(1-\left(\frac{x}{L}\right)^2\right)^3}{1+x^2\operatorname{softplus}(q_\theta(x))}, \quad Q_0 = A(b,\gamma,p) >0.}
\end{equation}
This type of hard-constraint trial function is consistent with PINN approaches that encode boundary behavior directly in the network ansatz rather than enforcing all constraints only through penalty terms \cite{SukumarSrivastava2022,Berrone2023Dirichlet}.
This ansatz \eqref{E:constrained_ansatz} has three useful features. First, the factor
$$
\left(1-\left(\frac{x}{L}\right)^2\right)^n
$$
enforces $Q_\theta(L)=0$, also its derivative vanishes at $x=0$. In particular, we take $n=3$ because near $x=L$ this factor behaves like a constant multiple of $(L-x)^3$, so the envelope and its first two derivatives vanish at the right boundary. Since the residual involves the second derivative $Q_\theta''$, this stronger flattening reduces boundary-induced curvature and leads to more stable and accurate training. In our experiments, the choice $n=2$ also satisfied the boundary conditions, but produced less favorable numerical performance.
{ Second, $Q_0$ is the fixed peak amplitude from \eqref{E:coeff}, not a trainable parameter. Thus, the reported FF--PINN experiments use exact peak information, a freely trained positive amplitude would require a separate nonzero normalization to rule out collapse toward the trivial solution.}
Third, the denominator
$$
1+x^2\operatorname{softplus}(q_\theta(x))
$$
is positive and has zero derivative at $x=0$, so it preserves the condition $Q_\theta'(0)=0$. Additionally, the use of $\operatorname{softplus}$ lets the network reshape the profile in the interior regions of the domain, so the ansatz is not just one fixed polynomial envelope. We note that since the hard-constrained ansatz satisfies $Q_\theta(0)=Q_0$ and $Q_\theta(L)=0$ identically, no separate boundary loss or boundary-penalty parameter is required. All FF--PINN results reported below use this hard-constrained formulation.

In this way, the ansatz \eqref{E:constrained_ansatz} biases the optimization toward a positive, nontrivial ground-state profile rather than the trivial solution. In particular, the residual and boundary losses alone may not select the positive ground-state branch, especially when the Fourier embedding gives the network additional oscillatory flexibility. Thus, we use a hard ground-state trial function to bias the admissible profiles toward positive, localized, boundary-compatible solutions. This is consistent with the broader PINN literature on enforcing constraints through the architecture or trial function rather than relying only on soft penalty terms, see, e.g., \cite{LagarisLikasFotiadis1998,SukumarSrivastava2022}. Recent hard-constrained Fourier-feature PINN formulations provide further motivation for combining Fourier-type embeddings with built-in boundary or structural constraints, as in 
\cite{LiDengWuZhangLiWang2024,StraubBrendelMedvedevRosskopf2025}.
We also tested the simpler linear Dirichlet hard ansatz 
$$
Q_\theta(x)=Q_0 \,\left(1-\frac{x}{L}\right)+\frac{x}{L}\left(1-\frac{x}{L}\right)q_\theta\left(\frac{x}{L}\right);
$$
however, the hard-constrained decay ansatz \eqref{E:constrained_ansatz} performed better { for FF--PINN and was used for its final comparisons. The action-based DR baseline instead uses the linear trial \eqref{E:constrained_ansatz_linear}, as specified in Section~\ref{S:actionDeepRitz}.}

We note that the FF--PINN should be viewed not as a replacement for the classical Fourier spectral method, but as an architectural modification of the PINN: the differential-equation residual and boundary conditions remain the training objective, while the neural network is supplied with a richer trigonometric input representation. Its performance is therefore assessed by the same diagnostics used for the standard PINN, including the $L^2$ error and the $L^\infty$ error.

We summarize the implementation of the FF--PINN used in our numerical experiments in Algorithm~\ref{alg:FFPINN}. 

{\scriptsize
\begin{algorithm}[!htbp]
\caption{Fourier-feature PINN with Second-Order Optimizer for $Q^{\prime\prime}=b Q-\gamma Q^{p}$}
\begin{algorithmic}[1]
\State Choose domain $[0,L]$, collocation points $\{x_j\}_{j=1}^{N_f}$, Fourier-feature dimension $N_{\mathrm{FF}}$, Fourier scale $\sigma_{\mathrm{FF}}$, and whether to include the raw coordinate
\State Initialize a fixed frequency matrix $B\in\mathbb R^{1\times N_{\mathrm{FF}}}$ with entries 
drawn from a Gaussian distribution with scale $\sigma_{FF}$.
\State Choose whether $B$ is kept fixed after initialization or trained together with the neural-network parameters.
\State Define the scaled coordinate
$$
\xi=\frac{x}{L}.
$$
\State Define the Fourier-feature mapping by
$$
\Phi(x) =\begin{cases}
  \left[\xi,\, \sin(2\pi\xi B),\, \cos(2\pi\xi B)
  \right], 
  & \text{if the raw coordinate is included}, \\[0.8em]
  \left[\sin(2\pi\xi B),\,\cos(2\pi\xi B)\right], & \text{otherwise}.
  \end{cases}
$$
\State Define a neural network 
$$
q_\theta(x)=\mathrm{MLP}_\theta(\Phi(x)).
$$
\State Define the FF--PINN approximation by
$$
Q_\theta(x)
=\displaystyle
{Q_0 } 
\frac{\left(1-\left(\frac{x}{L}\right)^2\right)^3}{1+x^2\operatorname{softplus}(q_\theta(x))},
$$
where 
{ $Q_0=A(b,\gamma,p)$ from \eqref{E:coeff} is fixed throughout training.}
\State Define the ODE residual 
$$
\mathcal R_\theta(x)
=
Q_\theta''(x)
-b Q_\theta(x)
+\gamma Q_\theta(x)^{p}.
$$
\State Define the residual loss
$$
\mathcal{L}_{\mathrm{ODE}}(\theta)
=
\frac{1}{N_f}
\sum_{j=1}^{N_f}
\left|\mathcal R_\theta(x_j)\right|^2.
$$
\State Define the total loss
$$
\mathcal{L}_{\mathrm{FFPINN}}(\theta)=\mathcal{L}_{\mathrm{ODE}}(\theta).
$$
\State Choose one of the optimizer modes
$$
\mathsf{Opt}\in\{\text{Adam}, \text{L--BFGS}, \text{GN/LM},\text{Adam}\to\text{L--BFGS},\text{Adam}\to\text{GN/LM }\}.
$$
\For{$n=1,2,\dots,N_{\mathrm{steps}}$}
  \State Use automatic differentiation to compute $Q_\theta'$ and $Q_\theta''$.
    \State Compute $\mathcal{L}_{\mathrm{FFPINN}}(\theta)=\mathcal{L}_{\mathrm{ODE}}(\theta)$.
  \State Update the parameters using the selected optimizer mode $\mathrm{Opt}$.
  \State Periodically resample collocation points and store checkpoint diagnostics.
\EndFor
\State Return the trained approximation $Q_\theta(x)$.
\end{algorithmic}
\label{alg:FFPINN}
\end{algorithm}
}

To optimize the FF--PINN implementation, we carried out the parameter sweeps summarized in Table~\ref{tab:ffpinn_hyperparams}, while fixing $p=3$. These tests varied the Fourier-feature dimension, Fourier scale, raw-coordinate inclusion, fixed versus trainable frequency matrices, optimizer schedule, collocation resolution, and refinement settings. To keep FF--PINN comparison consistent with the standard PINN baseline in section~\ref{S:PINN}, while also using the best-performing FF--PINN settings from the sweeps, the working FF--PINN parameters were chosen with $N_f=2^{9}$ total residual points, depth $=4$, and width $=32$.

{\scriptsize
\begin{table}[!htbp]
\centering
\begin{tabular}{ll}
\toprule
\textbf{Parameter} & \textbf{Values tested / Description} \\
\midrule
Domain length $L$ & $30$ \\
Boundary condition & Dirichlet\\

Activation functions 
& $\tanh(x)$, $\mathrm{SiLU}(x)$, $\mathrm{SiLUPLUS}(x)$, 
$e^{-x^2}$, $e^{-x^4}$, and $\sech(x)$ \\

Number of hidden layers & $3$, $4$, or $5$ \\

Neurons per layer & $32$, $64$, or $96$ \\

Optimizer & Adam, L--BFGS, GN/LM, Adam--L--BFGS, or Adam--GN/LM \\

Fourier feature count $N_{\mathrm{FF}}$ & $4$, $8$, $16$, or $32$ \\

Fourier scale $\sigma_{\mathrm{FF}}$ & $0.05$, $0.1$, $0.25$, or $0.5$ \\

Raw coordinate input 
& Tested both with and without the raw coordinate $x/L$ included \\

Fourier matrix $B$ & Fixed or trainable \\

Number of residual points $N_f$ & $2^8$, $2^9,\,2^{10},$ or $2^{11}$ \\

Adam warm-up length 
& $2{,}500$, $5{,}000$, $7{,}500$, $10{,}000$, $12{,}500$, and $15{,}000$ epochs \\
Adam learning rate $lr$ & $10^{-3}$ initial; StepLR decay \\

L--BFGS refinement steps & $25$, $50$, $75$, $100$, $125$, $150$, or $200$ \\

L--BFGS history size & $20$, $50$, or $100$ \\

L--BFGS max iterations per step & $1$, $5$, $10$, $20$, or $30$ \\

L--BFGS learning rate & $0.5$ or $1.0 $\\

GN/LM refinement steps & $25$, $50$, $75$, $100$, $125$, $150$, or $200$ \\

GN/LM initial damping & $5\times 10^{-4}$, $10^{-3}$, or $2\times 10^{-3}$ \\

GN/LM Maximum backtracking steps 
& $2$, $4$, 6, 8, and 10 \\

GN/LM Damping shrink factor 
& Values in $[0.4,1.1]$ \\

GN/LM Damping growth factor 
& $5$, $10$, and $15$ \\

Selection criterion & Smallest $L^2$ error, with runtime used to break near-ties \\
\bottomrule
\end{tabular}
\caption{FF--PINN architecture and optimization hyperparameters tested for solving \eqref{2ndODE} on $[0,30]$.}
\label{tab:ffpinn_hyperparams}
\end{table}
}

The Dirichlet boundary condition formulation favors $N_{\mathrm{FF}}=32$ and the $\mathrm{SiLU}$ activation. For this formulation, Adam--L--BFGS with $\sigma_{\mathrm{FF}}=0.25$ gave the best balance between accuracy and runtime. Although some more expensive refinements produced slightly smaller $L^2$ errors, the improvement was small compared with the additional CPU time. We therefore use the faster near-tie Adam--L--BFGS configuration as the current FF--PINN baseline. For the optimizer internals, the best Dirichlet results were obtained using an L--BFGS history size of $20$ and a maximum of $30$ inner iterations per outer step. The selected FF--PINN baseline configurations for our investigations are found in Table~\ref{tab:ffpinn_baselines}.

\begin{table}[!htbp]
\centering
\setlength{\tabcolsep}{3pt}
\renewcommand{\arraystretch}{1.05}
\resizebox{\textwidth}{!}{%
\begin{tabular}{cccccccccccccc}
\hline
 \textbf{Optimizer} 
& \textbf{Activation} 
& \shortstack[c]{\textbf{Hidden}\\\textbf{Layers}}
& \shortstack[c]{\textbf{Neurons}\\\textbf{per Layer}}
 & $N_f$
& \shortstack[c]{\textbf{Adam}\\ \textbf{Warm-up}}
& \textbf{$N_{\mathrm{FF}}$} 
& \textbf{$\sigma_{\mathrm{FF}}$} 
& \shortstack[c]{\textbf{Train}\\
\textbf{$B$}}
& \shortstack[c]{\textbf{Include} \\ \textbf{Raw}}
& \shortstack[c]{\textbf{L--BFGS}\\ \textbf{Steps}}
& \shortstack[c]{\textbf{History}\\ \textbf{Size}}
&  \shortstack[c]{\textbf{Internal}\\ \textbf{Iteration}}
& \shortstack[c]{\textbf{L--BFGS}\\ \textbf{LR}}
\\
\hline
Adam--L--BFGS 
& $\mathrm{SiLU}(x)$ 
& $4$
& $32$
& $2^9$
& $2{,}500$ 
& $32$ 
& $0.25$ 
& yes 
& no
& $75$ 
& $20$
& $30$
& $1.0$
\\
\hline
\end{tabular}
}
\caption{Selected FF--PINN configurations for Dirichlet boundary conditions.}
\label{tab:ffpinn_baselines}
\end{table}

Additionally, we found that including the raw scaled coordinate $x/L$ alongside the Fourier features did not give better accuracy than using only the sinusoidal embedding, so the following FF--PINN results do not include the raw input together with the trigonometric features.

\begin{table}[!htbp]
\centering
\setlength{\tabcolsep}{3pt}
\renewcommand{\arraystretch}{1.05}
\begin{tabular}{|c||c|c|c|c|}
\hline
 & \multicolumn{4}{c|}{FF--PINN with Dirichlet boundary conditions} \\
\hline\hline
$p$ & $\mathcal{L}_{\mathrm{FFPINN}}(\theta)$ & CPU time & $L^\infty$ error & $L^2$ error  \\
\hline\hline
3 & 1.756630e-07 & 219.333194 & 5.093381e-04 & 6.852868e-04 \\ \hline
5 & 1.004573e-07 & 200.046171 & 4.623447e-04 & 7.603941e-04 \\ \hline
7 & 5.792016e-08 & 234.282703 & 3.989951e-04 & 7.267693e-04 \\ \hline
\end{tabular}%
\caption{Comparison of FF--PINN performance with Dirichlet \ref{E:Dirichlet} boundary conditions (with parameters as in Table~\ref{tab:ffpinn_baselines}) for $Q^{\prime\prime}=Q-Q^{p}$.}
\label{tab:FFPINN_psweep}
\end{table}

We find that the FF--PINN gives accurate approximations for all tested powers $p=3,\,5,$ and $7$. As shown in Table~\ref{tab:FFPINN_psweep}, the final losses remain on the order of $10^{-8} \div 10^{-7}$, with $L^\infty$ errors below $5.1\times10^{-4}$ and $L^2$ errors below $7.7\times10^{-4}$. The CPU times are 
comparable across the three cases, ranging from approximately $200$ to $234$ seconds.
We obtain the smallest $L^2$ error for $p=3$,
while the smallest $L^\infty$ error is obtained for $p=7$. The errors remain of the same order across all three tested powers.

\begin{figure}[!htbp]
\centering
\includegraphics[width=0.49\textwidth, height=0.32\textwidth]{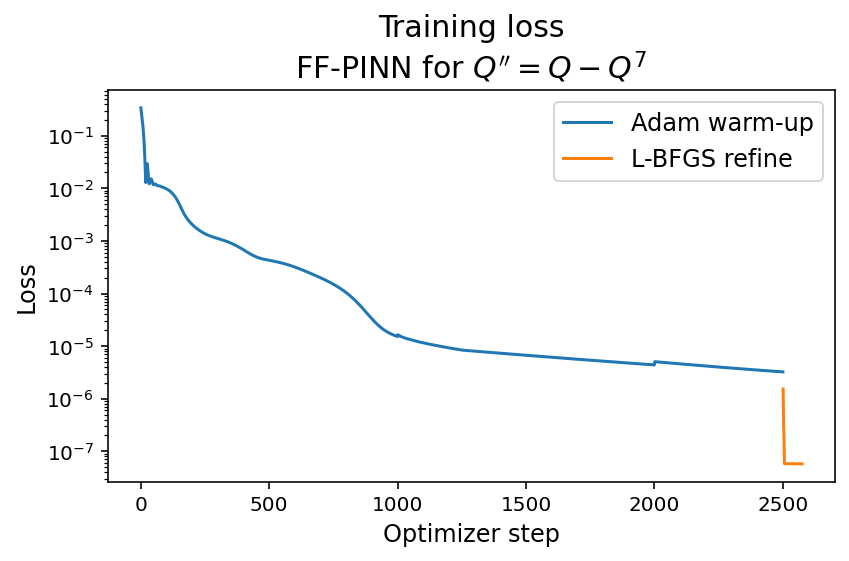}
\includegraphics[width=0.49\textwidth, height=0.32\textwidth]{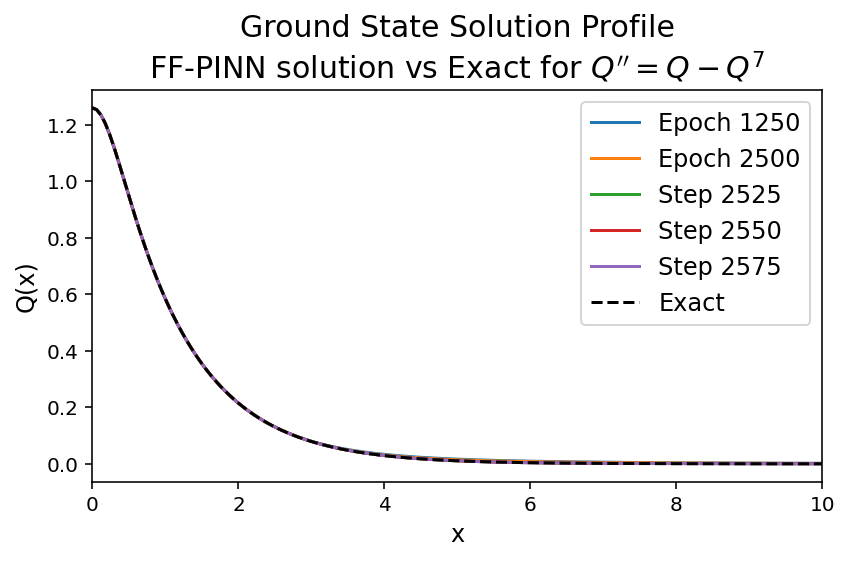}
\caption{Training loss (left) and ground state solution profile convergence (right) for FF--PINN (with parameters as seen in Table~\ref{tab:ffpinn_baselines}) with $Q^{\prime\prime}=Q-Q^5$. The (black dashed) curve labeled \emph{Exact} in the right panel is the explicit solution given by \eqref{Qexact}.  The final learned profile corresponds to $2{,}500$ Adam steps followed by $75$ L--BFGS steps.}
\label{fig:FFP_converge}
\end{figure}

We mention that we use the term \textit{training loss} in Figure~\ref{fig:FFP_converge} because $\mathcal{L}_{\mathrm{FFPINN}}(\theta)$ is the residual-based loss minimized by the FF--PINN. As described above, \textit{training objective} is the more general term that we use when referring collectively to residual-based PINN and FF--PINN losses and the variational objectives used for the Deep Ritz methods.

In the left panel of Figure~\ref{fig:FFP_converge}, we observe that the Adam warm-up stage steadily reduces the loss before the L--BFGS refinement stage further decreases it. In the right panel, we observe that the FF--PINN profiles rapidly converge to the exact ground state give by \eqref{Qexact}. By the final iterate, the numerical and exact curves are visibly indistinguishable on the plotted scale. The remaining discrepancy is consistent with the $p=5$ errors reported in Table~\ref{tab:FFPINN_psweep}, which is orders $10^{-4}$.
The representative convergence plot in Figure~\ref{fig:FFP_converge} illustrates this behavior for $p=5$: we see the gradual decrease during the Adam stage followed by the final L--BFGS refinement in the left panel, while the final profile in the right panel aligns closely with the exact solution throughout the plotted domain.

\subsubsection{Comparison with the ``secondary” computational device for FF--PINN}\label{S:secondaryFFPINN}
The FF--PINN experiments in this section use the same baseline hyperparameters selected in Table~\ref{tab:ffpinn_baselines}, but are run on the secondary GPU-capable workstation. No additional tuning was performed for these runs. In particular, the Fourier-feature dimension, Fourier-feature scaling, raw-coordinate inclusion, activation function, network architecture, collocation size, boundary-condition weight, Adam warm-up, and refinement parameters were all kept fixed to the previously selected FF--PINN baseline. The purpose of this test was to check whether the selected FF--PINN baseline gives comparable accuracy when moved to the GPU+CPU workstation, rather than to search for a new optimal configuration on the secondary device.

Table~\ref{tab:FFPINN_psweep_secondary} reports the FF--PINN results on the secondary GPU+CPU device. We observed that the method produces accurate ground-state approximations for all three power cases $p=3,\,5,\,7$, with the learned profiles remaining positive, centered at $x=0$, and visually consistent with the exact solution given by \eqref{Qexact}. The errors decrease as $p$ increases, with the smallest $L^\infty$ and $L^2$ errors obtained for $p=7$.
\begin{table}[!htbp]
\centering
\setlength{\tabcolsep}{3pt}
\renewcommand{\arraystretch}{1.05}
\begin{tabular}{|c||c|c|c|c|}
\hline
 & \multicolumn{4}{c|}{FF--PINN with Dirichlet boundary conditions} \\
\hline\hline
$p$ & $\mathcal{L}_{\mathrm{FFPINN}}(\theta)$ & GPU+CPU time & $L^\infty$ error & $L^2$ error  \\\hline
3 & 1.284711e-07 & 67.670952 & 7.727327e-04 & 1.107474e-03 \\ \hline
5 & 1.017865e-07 & 67.329122 & 4.419572e-04 & 7.254083e-04 \\ \hline
7 & 6.116914e-08 & 72.062521 & 2.975830e-04 & 5.681980e-04 \\ \hline
\end{tabular}%
\caption{Secondary-device GPU+CPU comparison of FF--PINN performance with Dirichlet boundary conditions \ref{E:Dirichlet} (with parameters as in Table~\ref{tab:ffpinn_baselines}) for $Q^{\prime\prime}=Q-Q^{p}$.}
\label{tab:FFPINN_psweep_secondary}
\end{table}

From Table~\ref{tab:FFPINN_psweep_secondary}, we observe that the $L^\infty$ errors remain of order $10^{-4}$ for all three powers. The $L^2$ errors are also of order $10^{-4}$ for $p=5$ and $p=7$, while the $p=3$ value is of order $10^{-3}$. Thus,
the errors decrease with $p$, and the fixed FF--PINN configuration continues to approximate the exact ground state accurately across the three tested powers. We also observe that GPU+CPU time 
remains comparable across the three cases, ranging from approximately $67$ to $72$ seconds.

Compared with the corresponding primary-device results in Table~\ref{tab:FFPINN_psweep}, the secondary device runs have
substantially have shorter the wall-clock times. In particular, the secondary-device times of approximately $67$--$72$ seconds are considerably smaller than the approximately $200$--$234$-second primary-device times. 
For $p=3$, the primary-device run gives smaller $L^\infty$ and $L^2$ errors, whereas for $p=5$ and $p=7$, the secondary-device runs give smaller errors for both $L^\infty$ and $L^2$. Although the two workstations produce somewhat different numerical results, the errors remain of comparable magnitude across the two devices, indicating that the selected FF--PINN configuration gives similar accuracy on both systems while the secondary-device runs require less wall-clock time.
Unlike the higher-order optimization PINN secondary-device runs, the FF--PINN experiments do not show evidence of convergence to a non-ground-state branch;
some numerical variation between the two workstation runs is expected because of different hardware, back-end execution, and other numerical details, but these differences do not substantially change the accuracy of the FF--PINN solutions.

Overall, we find that the secondary-device FF--PINN comparison supports the stability of the selected FF--PINN baseline: when the hyperparameters are held fixed, the GPU+CPU implementation preserves the same accuracy behavior observed on the primary CPU device while achieving shorter run times.

\section{Deep Ritz Variations}\label{S:DEEP_RITZ}

In this section we study the Deep Ritz formulations, which differ substantially from the residual-based PINN and FF--PINN methods studied above. In the PINN-type methods, the neural network is trained by minimizing a least-squares loss built directly from the {\it pointwise ODE} residual together with the relevant boundary residuals. By contrast, the Deep Ritz method is {\it variational}: the neural-network approximation $Q_\theta$ of the profile $Q$ is trained by minimizing a scalar objective derived from the variational structure of the ground state equation, rather than from the equation itself. This changes both the optimization landscape and the interpretation of the reported training objective.

The variational formulation cannot be minimized naively over an unconstrained neural-network class. This is because in the action-based formulation, the negative nonlinear term (see, for example, \eqref{E:action}) can drive the action away from the desired ground-state branch unless additional structure is imposed. Therefore, we use {\it constrained} Deep Ritz formulations, which are designed to select the positive, localized ground-state profile. The constraints that we use include hard boundary conditions, prescribed peak-amplitude information, positivity or branch selection, and, in the Nehari and Weinstein variants, additional integral constraints derived from the structure of the ground state equation.

A common feature of all Deep Ritz variants considered below is that the Adam training stage and the GN/LM refinement optimize different quantities. During the Adam stage, the method minimizes the scalar training objective associated with the selected formulation: the constrained action-based Deep Ritz loss, the Nehari-constrained objective, or the Weinstein quotient objective. These quantities are variational objectives rather than pointwise residual losses, thus, their final values are not expected to vanish even when the computed profile is accurate.

Here, GN/LM is used only as a final residual-based refinement stage. 
Accordingly, the GN/LM steps act on a separate, formulation-dependent least-squares residual vector. We define the corresponding refinement objective by 
\begin{equation}\label{eq:GN-refinement}
\mathcal{J}_{\mathrm{GN}}(\theta)=\frac{1}{2}\|F_{\mathrm{GN}}(\theta)\|_{l^2}^2,
\end{equation}
where $F_{\mathrm{GN}}(\theta)$ is a finite-dimensional residual vector. Note that it has the same least-squares structure as the PINN residual vector \eqref{eq:residual-vector}, but 
its entries depend on the Deep Ritz formulation. For the constrained and Nehari-constrained Deep Ritz methods, it contains the residual components specified in Algorithm~\ref{alg:DeepRitz} and \ref{alg:DeepRitz_nehari}, including, where applicable, the ODE, boundary-condition, non-negativity, and normalized Pokhozhaev residuals. For the Weinstein formulation, $F_{\mathrm{GN}}(\theta)$ instead consists of quadrature-weighted values of the Euler-Lagrange residual of the normalized profile
$$
r_{v,\theta}(x) = -v_\theta''(x)
+bv_\theta(x) -R_{b,L}(\theta)|v_\theta(x)|^{p-1}v_\theta(x),
$$
evaluated at the GN/LM residual points. Thus, the Weinstein refinement acts on the normalized profile $v_\theta$ rather than directly on the ODE residual of the reconstructed ground-state profile $Q_\theta$. By \eqref{E:WDR-EL-residual} the two residuals differ only by the positive constant factor $\left(R_{b,L}(\theta)/\gamma\right)^{1/(p-1)}$, so they have the same zero set, but the choice does affect the scaling of the Jacobian relative to a fixed damping parameter $\lambda$.

Throughout this section the Adam-stage quantity and the GN/LM stage quantity should be interpreted separately: the former measures the variational training objective of the selected Deep Ritz formulation, while the latter measures the residual-refinement objective. This distinction is especially important when interpreting the training-objective curves and the reported values of the Deep Ritz variational objectives and the separate GN/LM refinement objective $\mathcal{J}_{\mathrm{GN}}$.

{ We point out that the comparisons below concern complete, separately tuned training procedures (not a controlled comparison in which only the variational objective is changed). Every reported Deep Ritz configuration includes GN/LM refinement, and the methods differ in trial parameterization, architecture, penalties, and optimization budget. Thus, the final errors and runtimes cannot be attributed to the variational pretraining alone. Instead we wanted to show a more complete comparison of how well the methods could perform.}

\subsection{Constrained Deep Ritz Formulation}\label{S:actionDeepRitz}
To find the (unique) positive ground state branch, we use an action-based Deep Ritz formulation. Instead of training directly on the pointwise ODE residual, we use the action associated with the ground state profile equation \eqref{2ndODE-again}: for a profile $Q$ on $[0,L]$, the action is defined as
\begin{equation}\label{E:action}
\mathcal S_b(Q) = \int_0^L
\left( \frac12 |\partial_x Q|^2
+\frac{b}{2} Q^2 -\frac{\gamma}{p+1} |Q|^{p+1}
\right) \,dx.
\end{equation}
The Deep Ritz stage { performs local action-based optimization over a neural-network trial, satisfying the prescribed boundary values,  
we supplement the action with Pokhozhaev identity penalties.

{
The boundary conditions $Q_\theta(0)=Q_0$ and $Q_\theta(L)=0$
are used to follow the positive ground state branch and are enforced
exactly by the trial function below. 
Although these conditions alone do not make $\mathcal S_b$ lower bounded, we show below that the prescribed height together with the Pokhozhaev penalties makes the desired ground state a strict local minimizer of the corresponding half-line objective. }

For the Dirichlet formulation, we use the hard-constrained trial function 
{
\begin{equation}\label{E:constrained_ansatz_linear}
Q_\theta(x)
=\underbrace{Q_0\left(1-\frac{x}{L}\right)}_{\text{first term}}
+\underbrace{\frac{x}{L}\left(1-\frac{x}{L}\right)
q_\theta\left(\frac{x}{L}\right)}_{\text{trainable correction}}.
\end{equation}
This ansatz enforces the boundary conditions at $0$ and at $L$ 
exactly for all trainable parameters. Indeed, the first term satisfies the prescribed endpoint values ($Q_0$ or $0$), while the second term is a trainable correction as we indicated in \eqref{E:constrained_ansatz_linear} that vanishes at both endpoints, thus, putting them together we get $Q_\theta(0)=Q_0+0=Q_0$ and
$Q_\theta(L)=0+0=0$. Therefore, for the hard Dirichlet trial function in \eqref{E:constrained_ansatz_linear}, the boundary residual is identically zero (i.e., it leaves these boundary values unchanged).

One can observe that this trial function does not enforce positivity, monotonicity, or an upper bound $Q_\theta\leq Q_0$, in particular, the network can take $q_\theta \equiv t$ in \eqref{E:constrained_ansatz_linear}, giving
$$
Q_t(x)=Q_0\left(1-\frac{x}{L}\right)+t\,\frac{x}{L}\left(1-\frac{x}{L}\right),
\quad \mbox{and thus,} \quad \mathcal S_b(Q_t) \rightarrow-\infty
\quad (\mbox{as}~t\to\infty),
$$
and furthermore, the normalized Pokhozhaev penalty in Algorithm~\ref{alg:DeepRitz} remains bounded along this family, and thus, for every fixed finite penalty weight, the total objective also tends to $-\infty$ along this family. 

To explain why this constrained formulation works, we first set 
the constrained Deep Ritz objective to be 
\begin{equation}\label{Loss-DR}
\mathcal{L}_{\mathrm{DR}}(\theta)
= \mathcal S_b(Q_\theta)
+\lambda_{\mathrm{Pokh}}\mathcal{L}_{\mathrm{Pokh}}(\theta),
\end{equation}
where $\mathcal{L}_{\mathrm{Pokh}}(\theta)$ penalizes violations of the Pokhozhaev identities \eqref{E:P1}--\eqref{E:P2}, which are evaluated on the domain $[0,L]$. (Note that a boundary penalty $\mathcal{L}_{\mathrm{BC}}(\theta)$ would normally appear in \eqref{Loss-DR} as well, however, for the trial function in \eqref{E:constrained_ansatz_linear} it vanishes identically, so we drop it from the objective \eqref{Loss-DR}}.)

Next, we note that the ground-state branch is searched via the prescribed central value and the Pokhozhaev penalties. Indeed, the condition
$$
Q_\theta(0) \equiv Q_0
=\left(\frac{b(p+1)}{2\gamma}\right)^{1/(p-1)}
$$
excludes the trivial solution and fixes the correct ground state height,
while the trainable correction determines the profile shape. Moreover, the amplitude variation $Q\mapsto tQ$ is no longer admissible.
The two Pokhozhaev penalties encourage the integral equalities satisfied by the ground state, including
$$
\|Q_\theta'\|_{L^2(0,L)}^2
\approx
\frac{p-1}{p+3}\,b\|Q_\theta\|_{L^2(0,L)}^2.
$$
Thus, the boundary value provides the correct height, while the
variational objective and the integral penalties guide the remaining
shape approximation.

We now give a more precise local variational argument.
For this calculation, consider the exact ground state $Q$ on
$(0,\infty)$, and interpret $\mathcal S_b$, $\mathcal P_1$,
$\mathcal P_2$, and $\mathcal L_{\mathrm{DR}}$ from \eqref{E:action}, \eqref{E:P1}, \eqref{E:P2} and \eqref{Loss-DR}, correspondingly, as their counterparts defined on the half line, i.e.,  
with all integrals taken over $(0,\infty)$.
Then perturbations (preserving the prescribed height) have the form $Q+h$, where $h\in H^1(0,\infty)$ is real-valued and satisfies $h(0)=0$.
Since $Q$ solves the profile equation \eqref{2ndODE-again},
the first variation of the action vanishes for every such
perturbation, namely, 
$$
\mathcal S_b'(Q)[h]
=\int_0^\infty \left(Q'h' + b Q h - \gamma Q^ph \right)\,dx=0.
$$
Thus, $Q$ is a critical point of the action among profiles with
the prescribed boundary value at the origin.
Differentiating the profile
equation gives
\begin{equation}\label{E:kerQ'}
\left(-\partial_x^2+b-p\gamma Q^{p-1}\right)Q'=0.
\end{equation}
Recall that for the branch that we are looking for, it holds $Q'(x)<0$ for $x>0$. Then integration by parts yields the following representation of the second variation
\begin{align*}
\mathcal S_b''(Q)[h,h]
=\int_0^\infty \left(|h'|^2+\left(b-p\gamma \,Q^{p-1}\right)h^2
\right)\,dx
=\int_0^\infty (Q')^2
\left|\left(\frac{h}{Q'}\right)'\right|^2\,dx \geq 0.
\end{align*}
(To justify this, we first consider smooth perturbations supported
away from the endpoints and then extend by density to
$h\in H^1(0,\infty)$ with $h(0)=0$.) Note that the equality holds precisely when $h$ is a multiple of $Q'$.
Thus, prescribing the central value removes the negative quadratic
directions of the action, leaving only this one-dimensional
degeneracy (i.e., $\operatorname{span} \{Q' \}$).

To control this remaining zero direction, we invoke the first Pokhozhaev penalty: using the profile equation, we obtain
\begin{align*}
\mathcal P_1'(Q)[h]
= 2\int_0^\infty(Q'h'+bQh)\,dx
-(p+1)\gamma\int_0^\infty Q^p h\,dx =
-(p-1)\gamma\int_0^\infty Q^p h\,dx,
\end{align*}
and hence,
$$
\mathcal P_1'(Q)[Q']
 \frac{\gamma(p-1)}{p+1}Q_0^{p+1}>0.
$$
Since $\mathcal P_1(Q)=\mathcal P_2(Q)=0$, the second variation of
the penalized objective is
\begin{align}\label{E:L''}
\mathcal L_{\mathrm{DR}}''(Q)[h,h]
={}\mathcal S_b''(Q)[h,h] +
\frac{2\lambda_{\mathrm{Pokh}}}{s(Q)^2}
\left(|\mathcal P_1'(Q)[h]|^2 + |\mathcal P_2'(Q)[h]|^2 \right),
\end{align}
where
$$
s(Q) = b\|Q\|_{L^2(0,\infty)}^2
+\|Q'\|_{L^2(0,\infty)}^2
+\gamma\|Q\|_{L^{p+1}(0,\infty)}^{p+1}
+\varepsilon > 0
$$
is the analog of the normalization scale
$s_\theta=bM_\theta+K_\theta+\gamma N_\theta+\varepsilon$ for the half line, which we use in Algorithm~\ref{alg:DeepRitz}, evaluated at the exact ground state. (In the implementation, we use $\varepsilon=10^{-12}$, and the scale is computed from the current neural approximation.)

The first term in \eqref{E:L''} is nonnegative, and its only zero direction is detected by the first penalty. 
Consider the linearized operator
$-\partial_x^2+b-p\gamma Q^{p-1}$ on $L^2(0,\infty)$,
with domain
$\{ h \in H^2(0,\infty): h(0)=0 \}$.
The boundary condition $h(0)=0$ corresponds to perturbations that 
preserve the height $Q(0)=Q_0$.
The identity \eqref{E:kerQ'} shows that this operator
is nonnegative and has kernel $\operatorname{span}\{Q'\}$.
Since $Q^{p-1}(x) \to 0$ as $x \to \infty$, its essential spectrum
is $[b,\infty)$. Therefore, zero is a simple isolated eigenvalue,
and the quadratic form is bounded below by a positive multiple
of $\|h\|_{L^2(0,\infty)}^2$ on the $L^2$-orthogonal complement
of $Q'$.
Thus, for every fixed $\lambda_{\mathrm{Pokh}}>0$, there
exists $c_\lambda>0$ such that
$$
\mathcal L_{\mathrm{DR}}''(Q)[h,h]
\geq c_\lambda\|h\|_{H^1(0,\infty)}^2
\quad \mbox{with} \quad h(0)=0.
$$
Then using the Taylor expansion, for sufficiently small
$\|h\|_{H^1(0,\infty)}$, we obtain
$$
\mathcal L_{\mathrm{DR}}(Q+h)-\mathcal L_{\mathrm{DR}}(Q)
\geq \frac{c_\lambda}{4}\|h\|_{H^1(0,\infty)}^2.
$$
In particular, $Q$ is an isolated strict local minimizer of the
penalized half-line objective $\mathcal L_{\mathrm{DR}}$ among profiles with the height $Q(0)=Q_0$. (In our implementation, the normalization denominator is recomputed from the current neural approximation at each evaluation but treated as constant when computing derivatives. Although this omits derivatives of the denominator, the resulting gradient has the same linearization at $Q$ as the gradient of the fully differentiated objective, since both Pokhozhaev residuals vanish there.)

In summary, the given height $Q_0$ together with the Pokhozhaev
penalties is sufficient to select the desired ground state locally,
for every $p>1$ (including the cubic, quintic, and septic cases that we consider here).
For the finite-interval problem, functions satisfying $u(L)=0$ can
be extended by zero to the half-line, and a smooth cutoff of $Q$
near $L$ gives an admissible approximation, whose $H^1$ error is
exponentially small as the length $L \to \infty$. The local coercivity then gives a nearby local minimizer of the continuous truncated
objective for sufficiently large $L$. This gives an analytical justification for the local action-based training. In our
computations, Adam pretraining produces a profile near the desired
ground state, and GN/LM subsequently refines it by reducing the ODE
residual, the normalized Pokhozhaev residuals, and a nonnegativity
penalty while retaining the prescribed boundary values.
The agreement with the ground state is verified by
Figure~\ref{fig:DR_converge} and the profile errors in
Table~\ref{tab:Deep_ritz_psweeps}.
} 

This method can be referred to as a {\it constrained Deep Ritz formulation with ground-state branch selection} (rather than a standard residual-based PINN). Our implementation of this DR method, which we use in the numerical experiments, is summarized in Algorithm~\ref{alg:DeepRitz}.

{ We mention that we use the whole line expressions $P_{1,\theta}$ and $P_{2,\theta}$ in Algorithm~\ref{alg:DeepRitz} as ground state regularizers, evaluated by quadrature on $[0,L]$ as ground state regularizers rather than exact finite interval identities. For the trial function \eqref{E:constrained_ansatz_linear}, the omitted boundary terms are $Q_0 Q_\theta'(0)$ and $L |Q_\theta'(L)|^2$, neither of which is forced to vanish by the given boundary values. For comparison, the exact benchmark satisfies $Q_{\mathrm{exact}}'(0)=0$ and $Q_{\mathrm{exact}}'(x)=-\sqrt b\,Q_{\mathrm{exact}}(x)\tanh(\beta x)$. At $b=\gamma=1$, $L=30$ (and $p=3,5,7$), we have $Q_{\mathrm{exact}}(L)<3\times10^{-13}$ and $L|Q_{\mathrm{exact}}'(L)|^2<3\times10^{-24}$, therefore, these estimates motivate the use of the uncorrected expressions as regularizers for the intended whole line profiles.
We note that, for example, our Nehari-constrained approach in Section~\ref{S:NehariDeepRitz}, in contrast, retains the exact finite interval second Pokhozhaev identity correction.
}

{\scriptsize
\begin{algorithm}[!htbp]
\caption{Constrained Deep Ritz Method with  Pokhozhaev Penalties}
\label{alg:DeepRitz}
\begin{algorithmic}[1]
\State Choose $[0,L]$, quadrature points $\{x_j\}_{j=1}^{N_q}$,
$b>0$, $\gamma>0$, $p>1$, and
a Pokhozhaev penalty weight $\lambda_{\mathrm{Pokh}}>0$.

\State  Compute the peak value $Q_0$ from the condition
$$
\frac{b}{2}Q_0^2-\frac{\gamma}{p+1}Q_0^{p+1}=0,
$$
whose positive root is
$Q_0=\left(\frac{b(p+1)}{2\gamma}\right)^{\frac{1}{p-1}}$,
in agreement with \eqref{E:coeff}.

\State Define a neural network $q_\theta$ and { the hard linear trial $Q_\theta$ in \eqref{E:constrained_ansatz_linear}, only $Q_\theta(0)=Q_0$ and $Q_\theta(L)=0$ are imposed exactly.} 
\State Choose
$\mathsf{Opt}_{\mathrm{DR}}\in
\{\text{Adam},\text{L--BFGS},\text{Adam}\to\text{L--BFGS}\}$.
\For{$n=1,2,\ldots,N_{\mathrm{DR}}$}
\State Use automatic differentiation and quadrature to compute
$$
M_\theta=\int_0^LQ_\theta^2dx,\quad
K_\theta=\int_0^L|Q_\theta'|^2dx,\quad
N_{\theta}=\int_0^L|Q_\theta|^{p+1}dx.
$$
\State Compute
$$
\mathcal S_b(Q_\theta)=\frac12 K_\theta+\frac{b}{2} M_\theta -\frac{\gamma}{p+1} N_{\theta}.
$$
\State Compute the { whole-line Pokhozhaev expressions, truncated to $[0,L]$, used as regularizers}
$$
P_{1,\theta} =b M_\theta+K_\theta-\gamma N_{\theta},
\qquad
P_{2,\theta} =b M_\theta-K_\theta-\frac{2\gamma}{p+1}N_\theta.
$$
\State Define
$$
\mathcal L_{\mathrm{Pokh}}
 =\left|\frac{ P_{1,\theta}}{s_{\theta}}\right|^2
 +\left|\frac{ P_{2,\theta}}{s_{\theta}}\right|^2, \qquad 
s_\theta=b M_\theta + K_\theta + \gamma N_\theta + \varepsilon.
$$
where $s_\theta$ nondimensionalizes the residuals by the magnitude of the individual bulk terms and $\varepsilon>0$ is a numerical safeguard (its precise definition is formulation-dependent).

\State Set
$$
\mathcal L_{\mathrm{DR}} =\mathcal S_b(Q_\theta)
+\lambda_{\mathrm{Pokh}}\mathcal L_{\mathrm{Pokh}}
$$
The hard trial makes the boundary loss identically zero.
\State { Take a local optimizer step with  $\mathsf{Opt}_{\mathrm{DR}}$ on the fixed quadrature grid and store checkpoint diagnostics.}
\EndFor
\If{GN/LM refinement is selected}
 \State On at most $N_{\mathrm{GN}}$ uniform points, form $F_\mathrm{GN}(\theta)$ from the quadrature
ODE residual, the two normalized Pokhozhaev residuals, 
and non-negativity residual.
 \For{$m=1,2,\ldots,N_{\mathrm{GN/LM}}$}
  \State Form $J_F(\theta)=\partial F_{\mathrm{GN}}/\partial\theta$ and solve the residual-space system
  $(J_FJ_F^T+\mu I)\alpha= F_{\mathrm{GN}}$, then set $\delta=-J_F^T\alpha$, with step clipping.

\State Accept $\theta\leftarrow\theta+\delta$ (after backtracking) when it decreases
$\mathcal J_{\mathrm{GN}}(\theta) = {\tfrac12}\|F_{\mathrm{GN}}(\theta)\|_{l^2}^2$, shrink $\mu$ on acceptance and grow it otherwise.
  
\EndFor
\EndIf
\State Return the trained approximation $Q_\theta$.
\end{algorithmic}
\end{algorithm}
}

The terms ``soft'', ``hard'', and ``decaying'' in the Trial function line in Table~\ref{tab:deepritz_hyperparams} describe how the boundary conditions are imposed. In the soft formulation, the network output is unconstrained and the boundary conditions are enforced through penalty terms. In the hard formulation, the trial function is constructed so that the boundary conditions are satisfied exactly for every choice of the trainable parameters, as in \eqref{E:constrained_ansatz_linear}. In the decay formulation, the trial function also satisfies the boundary conditions exactly, but uses the stronger built-in decay structure in \eqref{E:constrained_ansatz} instead of the linear boundary interpolator. We tested all three formulations; however, the hard formulation with the linear constraint \eqref{E:constrained_ansatz_linear} produced the smallest errors, therefore, in the actual algorithm (see Algorithm ~\ref{alg:DeepRitz}) we drop the boundary condition term $\mathcal \lambda_{\mathrm{BC}}\mathcal{L}_{\mathrm{BC}} (\theta)$ from the loss $\mathcal{L}_{\mathrm{DR}}(\theta)$.

As discussed above, the Adam stage minimizes the action-based Deep Ritz
objective, while the optional GN/LM stage is used only as a residual-based
refinement of the pretrained profile. Therefore, the GN/LM parameters in
Table~\ref{tab:deepritz_hyperparams} should be interpreted as refinement
parameters for the least-squares residual system, not as parameters of the
scalar action minimization itself.

{\scriptsize
\begin{table}[!htbp]
\centering
\begin{tabular}{ll}
\toprule
\textbf{Parameter} & \textbf{Values tested / Description} \\
\midrule
Domain length $L$ & $30$ \\
Boundary condition & Dirichlet\\

Trial function 
& Soft penalty trial, hard or decaying Dirichlet ansatz \\

Optimizer
& Adam, L--BFGS, or Adam--L--BFGS\\

Refinement strategy 
& GN/LM or Adam--GN/LM \\

Activation functions 
& $\tanh(x)$, $\mathrm{SiLU}(x)$, $\mathrm{SiLUPLUS}(x)$, 
$e^{-x^2}$, $e^{-x^4}$, and $\sech(x)$ \\

Number of hidden layers & $3$, $4$, or $5$ \\

Neurons per layer & $32$, $64$, or $96$ \\

Number of quadrature points $N_q$ & $2^8$, $2^9,\,2^{10},$ $2^{11}$, or $2^{12}$ \\

Adam warm-up length 
& $2{,}500$, $5{,}000$, $7{,}500$, $10{,}000$, $12{,}500$, and $15{,}000$ epochs \\

Adam learning rate $lr$ &  $10^{-3}$ initial; StepLR decay \\

Pokhozhaev penalty weight $\lambda_{\mathrm{Pokh}}$ 
& $0.1$, $1$, $5$, $10$, or $25$ \\

L--BFGS refinement steps 
& $25$, $50$, $75$, $100$, $125$, $150$, or $200$ \\

L--BFGS history size & $20$, $50$, or $100$ \\

L--BFGS max iterations per step & $1$, $5$, $10$, $20$, or $30$ \\

L--BFGS learning rate & $0.5$ or $1.0 $\\

GN/LM refinement steps 
& $25$, $50$, $75$, $100$, $125$, $150$, or $200$ \\

GN/LM initial damping 
& $5\times 10^{-4}$, $10^{-3}$, or $2\times 10^{-3}$ \\

GN/LM Maximum backtracking steps 
& $2$, $4$, 6, 8, and 10 \\

GN/LM Damping shrink factor 
& Values in $[0.4,1.1]$ \\

GN/LM Damping growth factor 
& $5$, $10$, and $15$ \\

Selection criterion 
& Smallest $L^2$ error, with runtime used to break near-ties \\
\bottomrule
\end{tabular}
\caption{Deep Ritz architecture, constraint, and optimization hyperparameters tested for solving \eqref{2ndODE} on $[0,30]$.}
\label{tab:deepritz_hyperparams}
\end{table}
}

The hyperparameter sweep in Table~\ref{tab:deepritz_hyperparams} was used to select the final Deep Ritz configurations. The collocation sweeps gave similar errors for $2^8$, $2^9$, $2^{10}$, $2^{11}$, and $2^{12}$ action quadrature points, so $N_q = 2^9$, was
retained. The GN/LM residual refinement uses an initial damping of $10^{-3}$, at most $8$ backtracking steps, and a fixed maximum of $128$ GN residual points to control memory usage. 
Increasing the GN/LM residual or quotient grid did not materially change the final reported $L^\infty$ or $L^2$ errors, but increasing either grid increased memory usage and runtime. Therefore, we retained the reduced 128-point GN/LM grid as a computationally efficient refinement grid, while verifying the final solution on an independent uniform evolution grid of $4096$ points on $[0,L]$, with the reported $L^2$ error evaluated by composite trapezoid quadrature.
All reported Deep Ritz results therefore use the selected configurations in Table~\ref{tab:deep_ritz_baselines}.

\begin{table}[!htbp]
\centering
\setlength{\tabcolsep}{3pt}
\renewcommand{\arraystretch}{1.05}
\resizebox{\textwidth}{!}{%
\begin{tabular}{cccccccccccccccc}
\hline
 \textbf{Constraint} 
& \textbf{Optimization} 
& \textbf{Refinement} 
& \textbf{Activation} 
& \shortstack[c]{\textbf{Hidden}\\\textbf{Layers}}
& \shortstack[c]{\textbf{Neurons}\\\textbf{per Layer}}
& $N_q$
& \shortstack[c]{\textbf{Adam}\\ \textbf{Warm-up}}
& \shortstack[c]{\textbf{GN/LM}\\ \textbf{Steps}}
& \shortstack[c]{\textbf{Max}\\ \textbf{Backtracks}}
& \shortstack[c]{\textbf{Damping}\\ \textbf{Shrink}}
& \shortstack[c]{\textbf{Damping}\\ \textbf{Growth}}
& \shortstack[c]{\textbf{Initial}\\ \textbf{Damping}}
& \shortstack[c]{\textbf{Residual}\\\textbf{Points}}
& \shortstack[c]{\textbf{Quotient}\\\textbf{Points}}
& $\lambda_{\mathrm{Pokh}}$ 
\\
\hline
 Hard  
& Adam
& GN/LM 
& $\mathrm{SiLU}(x)$ 
& $5$ 
& $32$ 
& $2^9$
& $2{,}500$ 
& $200$ 
& $8$
& $0.7$
& $10$
& $10^{-3}$
& $128$
& $128$
& $25$ 
\\

\hline
\end{tabular}
}
\caption{Selected Deep Ritz configurations for Dirichlet boundary conditions \eqref{E:Dirichlet}.}
\label{tab:deep_ritz_baselines}
\end{table}

Table~\ref{tab:Deep_ritz_psweeps} reports the performance of the constrained Deep Ritz method for $p=3,\,5,$ and 7. The reported $\mathcal{L}_{\mathrm{DR}}(\theta)$ values correspond to the final action-based Deep Ritz objective and are not expected to converge to zero. The residual quality after refinement is instead measured by $\mathcal{J}_{\mathrm{GN}}(\theta)$, which is the least-squares objective used in the GN/LM refinement stage.

We find that the constrained Deep Ritz method gives very accurate results for the cubic and quintic cases. For $p=3$, the method produces $\mathcal{J}_{\mathrm{GN}}(\theta)$ of order $10^{-11}$, with $L^\infty$ and $L^2$ errors both of order $10^{-6}$. For $p=5$ the residual-refinement objective is reduced further to order $10^{-14}$, while the $L^\infty$ and $L^2$ errors are of $10^{-5}$ and $10^{-6}$, respectively. The septic case $p=7$ is slightly more challenging: although $\mathcal{J}_{\mathrm{GN}}$ reaches order $10^{-16}$, the final errors increase to order $10^{-3}$. We interpret this as indicating that the sharper higher-power profile is harder to approximate accurately, even when the residual-based refinement objective is very small. In particular, these results show that a smaller value of $\mathcal{J}_{\mathrm{GN}}(\theta)$ does not necessarily imply a smaller profile error across different values of $p$.

\begin{table}[!htbp]
\centering
\begin{tabular}{|c||c|c|c|c|c|}
\hline
$p$
& $\mathcal{L}_{\mathrm{DR}}(\theta)$ 
& $\mathcal{J}_{\mathrm{GN}}(\theta)$
& CPU time
& $L^\infty$ error
& $L^2$ error \\
\hline\hline
3 & 6.666667e-01 & 9.620502e-11 & 264.092212 & 4.302695e-06 & 9.264227e-06 \\ \hline
5 & 6.801746e-01 & 2.375216e-14 & 252.688689 & 1.224911e-05 & 6.049316e-06 \\ \hline
7 & 6.679167e-01 & 4.063798e-16 & 246.570857 & 3.198027e-03 & 2.417198e-03 \\ \hline
\end{tabular}
\caption{Performance of the constrained Deep Ritz method with Dirichlet boundary conditions \eqref{E:Dirichlet} for $Q^{\prime\prime}=Q-Q^p$.}
\label{tab:Deep_ritz_psweeps}
\end{table}

The CPU runtimes are relatively consistent across the three powers, ranging from around $246$--$264$ seconds. Thus, we find that the main cost is determined primarily by the fixed network architecture, quadrature size, and GN/LM refinement settings, rather than by the particular value of $p$.

\begin{figure}[!htbp]
\centering
\includegraphics[width=0.49\textwidth, height=0.32\textwidth]{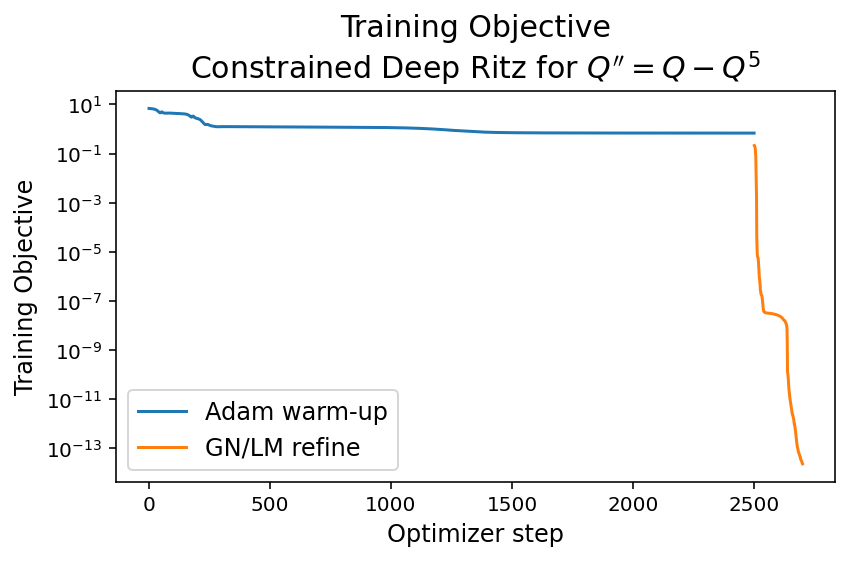}
\includegraphics[width=0.49\textwidth, height=0.32\textwidth]{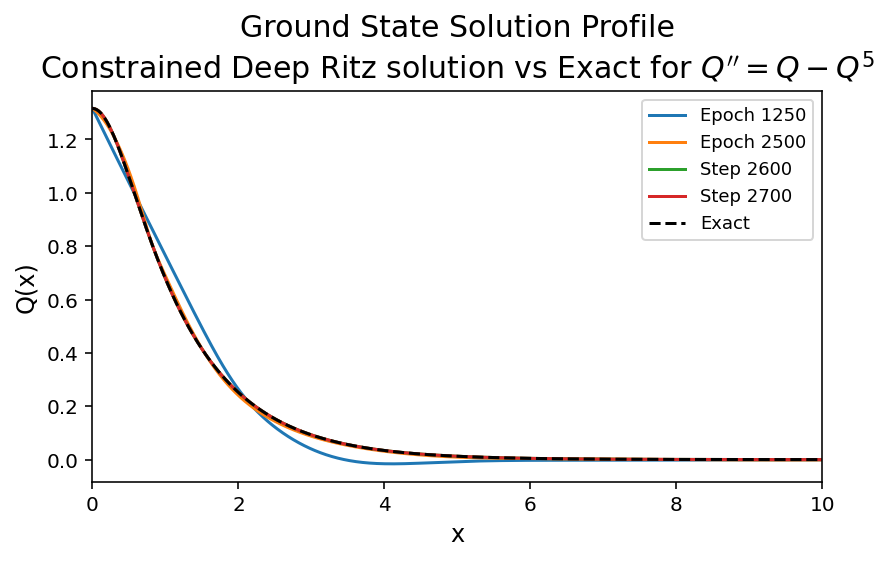}
\caption{Training objective (left) and ground state solution profile convergence (right) for constrained Deep Ritz method (with parameters as seen in Table~\ref{tab:deep_ritz_baselines}) for $Q^{\prime\prime}=Q-Q^5$. The (black dashed) curve labeled \emph{Exact} in the right panel is the explicit solution given by \eqref{Qexact}. The final learned profile corresponds to $2{,}500$ Adam steps followed by $200$ GN/LM steps.}
\label{fig:DR_converge}
\end{figure}

Figure~\ref{fig:DR_converge} illustrates the two-stage behavior of the method. In the left panel, we observed that, during the Adam warm-up, the action-based Deep Ritz objective decreases and stabilizes near its final variational value. Once GN/LM refinement begins, the residual-based objective decreases by several additional orders of magnitude. In the right panel, we compare the 
learned profiles with the exact ground state given by \eqref{Qexact}. We observe that the Adam-trained approximation already visually captures the main ground-state shape, while the GN/LM stage visually improves the agreement with the exact solution. The final profile after $2{,}500$ Adam steps and $200$ GN/LM steps aligns closely with the exact curve throughout the plotted domain. The remaining discrepancy is consistent with the $p=5$ errors in Table~\ref{tab:Deep_ritz_psweeps}: the final $L^\infty$ error is $1.224911\times 10^{-5}$, while the $L^2$ error is $6.049316\times 10^{-6}.$

\subsubsection{Comparison with the ``secondary” computational device for constrained Deep Ritz}\label{S:secondaryDeepRitz}

Table~\ref{tab:Deep_ritz_psweeps_secondary} reports the same constrained Deep Ritz configuration on the secondary GPU+CPU workstation. No additional hyperparameter tuning was performed on this device. The purpose of this test is to determine whether the selected constrained Deep Ritz baseline transfers reliably to a different computational environment and to quantify the corresponding runtime change. We observe that the secondary device substantially reduces the elapsed runtime for all three powers. The resulting primary-to-secondary runtime ratios correspond to speed-ups of approximately $6.4$, $6.5$, and $5.8$ for $p=3,\,5$, and $7$, respectively.

\begin{table}[!htbp]
\centering
\begin{tabular}{|c||c|c|c|c|c|}
\hline
$p$
& $\mathcal{L}_{\mathrm{DR}}(\theta)$ 
& $\mathcal{J}_{\mathrm{GN}}(\theta)$ 
& GPU+CPU time
& $L^\infty$ error
& $L^2$ error \\ \hline\hline
3 & 6.666667e-01 & 7.300457e-11 & 41.048747 & 3.697873e-06 & 7.913867e-06 \\ \hline
5 & 6.801749e-01 & 1.826449e-12  & 38.984522 & 2.828127e-05 & 1.107741e-05 \\ \hline
7 & 6.677846e-01 & 1.326228e-09 & 42.222312 & 1.397995e-03 & 1.007490e-03 \\ \hline
\end{tabular}
\caption{Secondary-device GPU+CPU performance of the constrained Deep Ritz method with Dirichlet
boundary conditions \eqref{E:Dirichlet} for $Q^{\prime\prime}=Q-Q^p$.}
\label{tab:Deep_ritz_psweeps_secondary}
\end{table}

The final errors on the secondary device remain comparable to those obtained on the primary device. For $p=3$, both the $L^\infty$ and $L^2$ errors improve slightly and remain of order $10^{-6}$, while for $p=5$, the secondary run give slightly larger errors: $L^\infty$ error remains of order $10^{-5}$, while the $L^2$ error increases from order $10^{-6}$ to order $10^{-5}$. For $p=7$, both errors decrease relative to the primary-device values but remain of order $10^{-3}$. Thus, we find that the secondary-device runs preserve the same accuracy patterns as the primary device runs: strong accuracy for $p=3$ and $p=5$, with a more challenging approximation for the sharper $p=7$ profile. 

We also observe that the relationship between the residual-refinement objective and the profile error is not monotone across the two devices. For example, for $p=7$, the secondary-device value of $\mathcal{J}_{\mathrm{GN}}(\theta)$ is of order $10^{-9}$ rather than $10^{-16}$, yet the secondary-device $L^\infty$ and $L^2$ errors are smaller. This again indicates that the magnitude of the sampled residual objective alone does not completely determine the error relative to the exact ground state \eqref{Qexact}.

Overall, we find that the constrained Deep Ritz method transfers well to the secondary GPU+CPU workstation. The main effect of the secondary device is a large runtime reduction, while the final errors remain comparable to the primary-device results. On the secondary device, we also observe that the constrained Deep Ritz method improves substantially over the corresponding FF--PINN results for $p=3$ and $p=5$, while the comparison is less decisive for the sharper $p=7$ case.

\subsection{Nehari-Constrained Deep Ritz Formulation}
\label{S:NehariDeepRitz}

For the ground-state solution which we consider on $[0,L]$, of equation
\eqref{2ndODE-again}, we consider the associated action
\begin{equation}\label{E:action-again}
\mathcal S_b(Q) = \int_0^L \left[\frac12 |Q'(x)|^2 +\frac b2 Q(x)^2
-\frac{\gamma}{p+1}Q(x)^{p+1} \right] dx.
\end{equation}
The ground state profile equation \eqref{2ndODE-again} is the Euler--Lagrange equation for this action,
$\mathcal S_b'(Q)=0$. However, the focusing ($\gamma >0$) action is not bounded from
below on the full trial space. To show that, consider any nonzero $u$ to get
$$
\mathcal S_b(tu)=\frac{t^2}{2}
\int_0^L\bigl(|u'|^2+bu^2\bigr)\,dx
- \gamma \, \frac{ t^{p+1}}{p+1} \int_0^L |u|^{p+1}\,dx
   ~~ \rightarrow -\infty \quad\text{as }t\to\infty.
$$
Thus, the ground state cannot be described as an {\it unconstrained} minimizer of $\mathcal S_b$. In the previous section, { in the constrained DR algorithm, we considered a {\it specific} branch, where the initial condition $Q_0=A$ helped us locate the correct branch, which was possible since we knew the explicit solutions, and hence, the appropriate branch of solutions (although we also showed why fixing the $Q_0$ gives the appropriate branch).  
However, if we did not have explicit solutions, and thus, the peak height explicitly, then we would need to add extra constraints or other conditions.} 

A possible constraint is, for example, to use one of the Pokhozhaev identities. Taking the first identity \eqref{E:P1}, sometimes also called the Nehari identity (and to distinguish it from the second Pokhozhaev identity, we denote it here via $\mathcal N_b$), we define it on a finite interval $[0,L]$
\begin{equation}\label{E:Nehari-DR}
\mathcal N_b(u):=\int_0^L\left(|u'|^2+bu^2-\gamma |u|^{p+1}\right)dx=0.
\end{equation}
The corresponding Nehari manifold is
$$
\mathcal{M}_L=\{u\in H^1(0,L) {{\setminus\{0\}} }: \mathcal N_b(u)=0 ~ \mbox{and} ~ u(L)=0\}.
$$
Then the nontrivial ground state is characterized variationally by minimizing
$\mathcal S_b$ over nonzero profiles on the Nehari manifold $\mathcal M_L$.  With this constraint, we have
\begin{equation}\label{E:action-on-Nehari}
\mathcal S_b(u) =\frac{p-1}{2(p+1)} \int_0^L\bigl(|u'|^2+bu^2\bigr)\,dx >0.
\end{equation}
Thus, the value of the action at the trained ground state is not
expected to converge to zero.

To avoid prescribing the unknown amplitude $Q(0)$, we first construct a
positive shape $\widetilde Q_\theta$ with a unit height. For that we set
$$
z=\left(\frac{x}{L}\right)^2,
\qquad \rho(s)=\log(1+e^s),
$$
and define
\begin{equation}\label{E:DR-shape}
\widetilde Q_\theta(x)
=(1-z)\frac{\rho\bigl(q_\theta(z)\bigr)}{\rho\bigl(q_\theta(0)\bigr)}.
\end{equation}
Then we get
$$
\widetilde Q_\theta(0)=1,
\qquad \widetilde Q_\theta'(0)=0,
\qquad \widetilde Q_\theta(L)=0,
\qquad \widetilde Q_\theta(x)>0\quad (0\leq x<L).
$$
We note that the symmetry condition at the origin and the truncated Dirichlet
condition at $x=L$ are now imposed exactly, while the amplitude remains to be
determined by the variational problem.

Define
\begin{equation}\label{E:DR-shape-integrals}
\begin{aligned}
\widetilde M_\theta
&=\int_0^L \widetilde Q_\theta^2\,dx,
 &
\widetilde K_\theta
&=\int_0^L |\widetilde Q_\theta'|^2\,dx,
&
\widetilde N_\theta
    &=\int_0^L \widetilde Q_\theta^{p+1}\,dx.
\end{aligned}
\end{equation}
Then the Nehari projection is obtained by setting
\begin{equation}\label{E:DR-amplitude}
    A_\theta
    =
    \left(
        \frac{\widetilde K_\theta+b\widetilde M_\theta}
             {\gamma\widetilde N_\theta}
    \right)^{\!1/(p-1)},
    \qquad
    Q_\theta=A_\theta\widetilde Q_\theta.
\end{equation}
From scaling it follows that 
$$
K_\theta=A_\theta^2\widetilde K_\theta,
\qquad
M_\theta=A_\theta^2\widetilde M_\theta,
\qquad
N_\theta=A_\theta^{p+1}\widetilde N_\theta,
$$
and thus, \eqref{E:DR-amplitude} gives
\begin{equation*}\label{E:DR-Nehari-exact}
K_\theta+bM_\theta-\gamma N_\theta=0.
\end{equation*}
Therefore, the first Pokhozhaev
identity is enforced by the above Nehari projection, and we do not need to include it as a penalty in constrained Deep Ritz approach.
The projected action can be evaluated as 
\begin{equation}\label{E:DR-projected-action}
\mathcal S_b(Q_\theta)
=\frac{p-1}{2(p+1)} A_\theta^2
\bigl(\widetilde K_\theta+b\widetilde M_\theta\bigr).
\end{equation}

For the truncated Dirichlet problem, the second Pokhozhaev identity contains a boundary correction as we discussed in \S \ref{S:Pokhozhaev}, thus, it becomes 
$P_{2,L}(Q):= bM(Q)-K(Q)-\frac{2\gamma}{p+1}N(Q) +L|Q'(L)|^2 =0.$
Note that the last term vanishes completely only in the infinite-domain. If we omit it, it may introduce a small truncation error. We therefore use
\begin{equation}\label{E:DR-Pokh-loss}
\mathcal L_{\mathrm{Pokh}}(\theta)
=\left| \frac{P_{2,L}(Q_\theta)}{s_\theta} \right|^2,
\quad \mbox{where} \quad s_\theta = bM_\theta+K_\theta
    +\frac{2\gamma}{p+1}N_\theta +L|Q_\theta'(L)|^2+\varepsilon.
\end{equation}
The Nehari-constrained Deep Ritz loss is
\begin{equation}\label{E:DR-corrected-loss}
\mathcal L_{\mathrm{NDR}}(\theta) =
    \mathcal S_b(Q_\theta)
    +\lambda_{\mathrm{Pokh}}
     \mathcal L_{\mathrm{Pokh}}(\theta).
\end{equation}
We note that no boundary penalty is required for the hard-constrained ansatz
\eqref{E:DR-shape}.

{\scriptsize
\begin{algorithm}[!htbp]
\caption{Nehari-Constrained Deep Ritz Method for
$Q''=bQ-\gamma Q^p$}
\label{alg:DeepRitz_nehari}
\begin{algorithmic}[1]
\State Choose $b>0$, $\gamma>0$, $p>1$, the domain $[0,L]$,
quadrature points $\{x_j\}_{j=1}^{N_q}$, and
$\lambda_{\mathrm{Pokh}}\geq0$.
\State Define a neural network $q_\theta(z)$ for $z\in[0,1]$ and set
$\rho(s)=\log(1+e^s)$.
\State Define the positive hard-constrained shape
$$
    z=\left(\frac{x}{L}\right)^2,
    \qquad
    \widetilde Q_\theta(x)
    =
    (1-z)
    \frac{\rho(q_\theta(z))}{\rho(q_\theta(0))}.
$$
\State Choose the variational optimizer
$$
    \mathsf{Opt}_{\mathrm{DR}}
    \in
    \{\text{Adam},\ \text{L--BFGS},\
      \text{Adam}\to\text{L--BFGS}\}.
$$
\For{$n=1,2,\ldots,N_{\mathrm{epochs}}$}
    \State Use automatic differentiation and quadrature to compute
    $$
        \widetilde M_\theta
        =\int_0^L\widetilde Q_\theta^2\,dx,
        \quad
        \widetilde K_\theta
        =\int_0^L|\widetilde Q_\theta'|^2\,dx,
        \quad
        \widetilde N_\theta
        =\int_0^L\widetilde Q_\theta^{p+1}\,dx.
    $$
    \State Compute the differentiable Nehari amplitude
    $$
        A_\theta
        =
        \left(
            \frac{\widetilde K_\theta+b\widetilde M_\theta}
                 {\gamma\widetilde N_\theta}
        \right)^{1/(p-1)}
    $$
    and set $Q_\theta=A_\theta\widetilde Q_\theta$.
    \State Keep $A_\theta$ in the automatic-differentiation graph; do not
    detach it from the trainable parameters.
    \State Compute
    $$
        M_\theta=A_\theta^2\widetilde M_\theta,
        \qquad
        K_\theta=A_\theta^2\widetilde K_\theta,
        \qquad
        N_\theta=A_\theta^{p+1}\widetilde N_\theta.
    $$
    \State Compute the projected action
    $$
        \mathcal S_b(Q_\theta)
        =
        \frac{p-1}{2(p+1)}
        A_\theta^2
        \bigl(\widetilde K_\theta+b\widetilde M_\theta\bigr).
    $$
    \State Compute the finite-interval Pokhozhaev residual
    $$
        P_{2,L,\theta}
        =
        bM_\theta-K_\theta
        -\frac{2\gamma}{p+1}N_\theta
        +L\left(Q_\theta'(L)\right)^2
    $$
    and the normalized penalty
    $$
        \mathcal L_{\mathrm{Pokh}}(\theta)
        =
        \left|\frac{P_{2,L,\theta}}{s_{\theta,L}}\right|^2,
    $$
    where
    $$
    s_{\theta,L}
    =
    bM_\theta+K_\theta
    +\frac{2\gamma}{p+1}N_\theta
    +L\left(Q_\theta'(L)\right)^2+\varepsilon,
    $$
    and $\varepsilon>0$ is a small numerical safeguard (the precise definition
    of the normalization scale $s_{\theta,L}$ is formulation-dependent).
    \State Form
    $$
        \mathcal L_{\mathrm{NDR}}(\theta)
        =
        \mathcal S_b(Q_\theta)
        +\lambda_{\mathrm{Pokh}}
         \mathcal L_{\mathrm{Pokh}}(\theta).
    $$
    \State Update $\theta$ using $\mathsf{Opt}_{\mathrm{DR}}$.
   
    \State Using the fixed quadrature grid, store
    $\mathcal S_b(Q_\theta)$, $\mathcal L_{\mathrm{Pokh}}(\theta)$, and
    $|K_\theta+bM_\theta-\gamma N_\theta|$.
\EndFor
\State Evaluate the pointwise ODE diagnostic
$$
    \mathcal L_{\mathrm{ODE}}(\theta)
    =
    \frac1{N_r}\sum_{j=1}^{N_r}
    \left|
        Q_\theta''(x_j)-bQ_\theta(x_j)
        +\gamma Q_\theta(x_j)^p
    \right|^2.
$$
\State Optionally apply a GN/LM residual refinement while retaining the
hard boundary conditions and the Nehari projection.
\State Return $Q_\theta$, together with
$\mathcal S_b(Q_\theta)$, $\mathcal L_{\mathrm{ODE}}(\theta)$,
$\mathcal L_{\mathrm{Pokh}}(\theta)$, and the Nehari residual.
\end{algorithmic}
\end{algorithm}
}
After the variational training stage, we evaluate the pointwise equation
residual
\begin{equation}\label{E:DR-ODE-diagnostic}
\mathcal{R}_\theta(x)
=Q_\theta''(x)-bQ_\theta(x)+\gamma Q_\theta(x)^p,
\qquad \mathcal L_{\mathrm{ODE}}(\theta)
=\frac1{N_r}\sum_{j=1}^{N_r}|\mathcal{R}_\theta(x_j)|^2.
\end{equation}

In this formulation, the network first represents a positive shape with a unit height, satisfying the hard boundary and symmetry constraints, and the physical amplitude is then determined through the Nehari projection. Thus, unlike the constrained Deep Ritz formulation above, the peak value $Q(0)$ is not prescribed from the exact solution. Instead, the amplitude is selected variationally through \eqref{E:DR-amplitude}, which enforces the Nehari identity up to quadrature and floating-point error.

\begin{remark}
The Nehari-projected action is directly related to the Weinstein quotient. Eliminating $A_\theta$ from \eqref{E:DR-projected-action} gives
\begin{equation}\label{E:DR-Weinstein-relation}
\mathcal S_b(Q_\theta)
=
\frac{p-1}{2(p+1)}\gamma^{-\frac2{p-1}}
\left[\frac{\widetilde K_\theta+b\widetilde M_\theta}{(\widetilde N_\theta)^{\frac2{p+1}} }
\right]^{\frac{p+1}{p-1}}.
\end{equation}
Hence, the Nehari-constrained DR method and the Weinstein DR method produce equivalent minimizers, in the sense that 
{ when $\lambda_{\mathrm{Pokh}}=0$ and the same raw shape and quadrature are used, the Nehari-projected action is a strictly increasing function of the Weinstein quotient and they have the same minimizing shapes. (We mention that our reported NDR and WDR procedures nevertheless differ through the additional second-identity penalty in NDR, the residuals used in GN/LM refinement, and their separately selected parameterizations and hyperparameters.)}
\end{remark}

{\scriptsize
\begin{table}[!htbp]
\centering
\begin{tabular}{ll}
\toprule
\textbf{Parameter} & \textbf{Values tested / Description} \\
\midrule
Domain length $L$ & $30$ \\
Boundary condition & mixed\\

Trial function 
& Soft penalty trial, hard or decaying Dirichlet ansatz \\

Optimizer
& Adam, L--BFGS, or Adam--L--BFGS\\

Refinement strategy 
& GN/LM or Adam--GN/LM \\

Activation functions 
& $\tanh(x)$, $\mathrm{SiLU}(x)$, $\mathrm{SiLUPLUS}(x)$, 
$e^{-x^2}$, $e^{-x^4}$, and $\sech(x)$ \\

Number of hidden layers & $3$, $4$, or $5$ \\

Neurons per layer & $32$, $64$, or $96$ \\

Number of quadrature points $N_q$ & $2^8$, $2^9,\,2^{10},$ $2^{11}$, or $2^{12}$ \\

Adam warm-up length 
& $2{,}500$, $5{,}000$, $7{,}500$, $10{,}000$,  $12{,}500$, or $15{,}000$ steps \\

Adam learning rate $lr$ &  $10^{-3}$ initial; StepLR decay \\

Pokhozhaev penalty weight $\lambda_{\mathrm{Pokh}}$ 
& $0.1$, $1$, $5$, $10$, or $25$ \\

L--BFGS refinement steps 
& $25$, $50$, $75$, $100$, $125$, $150$, or $200$ \\

L--BFGS history size & $20$, $50$, or $100$ \\

L--BFGS max iterations per step & $1$, $5$, $10$, $20$, or $30$ \\

L--BFGS learning rate & $0.5$ or $1.0 $\\

GN/LM refinement steps 
& $25$, $50$, $75$, $100$, $125$, $150$, or $200$ \\

GN/LM initial damping 
& $5\times 10^{-4}$, $10^{-3}$, or $2\times 10^{-3}$ \\

GN/LM Maximum backtracking steps 
& $2$, $4$, 6, 8, and 10 \\

GN/LM Damping shrink factor 
& Values in $[0.4,1.1]$ \\

GN/LM Damping growth factor 
& $5$, $10$, and $15$ \\

Selection criterion 
& Smallest $L^2$ error, with runtime used to break near-ties \\
\bottomrule
\end{tabular}
\caption{Nehari-constrained Deep Ritz architecture, constraint, and optimization hyperparameters tested for solving \eqref{2ndODE} on $[0,30]$.}
\label{tab:nehari_deepritz_hyperparams}
\end{table}
}

Table~\ref{tab:nehari_deepritz_hyperparams} lists the architecture, quadrature, constraint, and optimizer parameters tested in the hyperparameter sweep. The sweep varied the activation function, network depth and width, number of quadrature points, Adam warm-up length, L--BFGS settings, GN/LM refinement parameters, and the Pokhozhaev penalty weight $\lambda_{\mathrm{Pokh}}$. Since the Nehari projection enforces the first identity directly and the hard trial function enforces the boundary conditions, the remaining penalty is used only to control the finite-interval second Pokhozhaev residual.

\begin{table}[!htbp]
\centering
\setlength{\tabcolsep}{3pt}
\renewcommand{\arraystretch}{1.05}
\resizebox{\textwidth}{!}{%
\begin{tabular}{ccccccccccccccc}
\hline
\textbf{Optimization} 
& \textbf{Refinement} 
& \textbf{Activation} 
& \shortstack[c]{\textbf{Hidden}\\\textbf{Layers}}
& \shortstack[c]{\textbf{Neurons}\\\textbf{per Layer}}
& $N_q$
& \shortstack[c]{\textbf{Adam}\\ \textbf{Warm-up}}
& \shortstack[c]{\textbf{GN/LM}\\ \textbf{Steps}}
& \shortstack[c]{\textbf{Max}\\ \textbf{Backtracks}}
& \shortstack[c]{\textbf{Damping}\\ \textbf{Shrink}}
& \shortstack[c]{\textbf{Damping}\\ \textbf{Growth}}
& \shortstack[c]{\textbf{Initial}\\ \textbf{Damping}}
& \shortstack[c]{\textbf{Residual}\\\textbf{Points}}
& \shortstack[c]{\textbf{Quotient}\\\textbf{Points}}
& $\lambda_{\mathrm{Pokh}}$ 
\\
\hline
Adam
& GN/LM 
& $\tanh(x)$ 
& $5$ 
& $64$ 
& $2^9$
& $10{,}000$ 
& $200$ 
& $6$
& $0.7$
& $15$
& $10^{-3}$
& $128$
& $128$
& $0.1$ 
\\

\hline
\end{tabular}
}
\caption{Selected Nehari-constrained Deep Ritz configuration for mixed boundary conditions \eqref{E:mixed}.}
\label{tab:nehari_deep_ritz_baselines}
\end{table}

The final configuration selected from this sweep is reported in Table~\ref{tab:nehari_deep_ritz_baselines}. This baseline uses Nehari-constrained Adam training followed by GN/LM residual refinement, with the quadrature and residual-point counts fixed to balance accuracy, runtime, and memory usage. All Nehari-constrained Deep Ritz results below use this fixed configuration unless otherwise stated, so the following tables measure the performance of the selected Nehari baseline across the nonlinear powers rather than a separate hyperparameter search for each value of $p$.

Table~\ref{tab:Nehari_Deep_ritz_psweeps} reports the primary-device results for the Nehari-constrained Deep Ritz method. Unlike the constrained Deep Ritz formulation above, this method does not prescribe the exact peak amplitude $Q(0)=A$. Instead, the amplitude is determined dynamically through the Nehari projection \eqref{E:DR-amplitude}, while the hard trial function enforces positivity, symmetry at the origin, and the endpoint condition at $x=L$. Thus, this formulation tests whether the ground-state branch can be selected variationally without inserting the known exact amplitude.

\begin{table}[!htbp]
\centering
\begin{tabular}{|c||c|c|c|c|c|}
\hline
$p$
& $\mathcal{L}_{\mathrm{NDR}}(\theta)$ 
& $\mathcal{J}_{\mathrm{GN}}(\theta)$
& CPU time
& $L^\infty$ error
& $L^2$ error \\
\hline\hline
3 & 6.666667e-01 & 2.141298e-13 & 313.008810 & 1.690701e-07 & 2.733932e-07 \\ \hline
5 & 6.801748e-01 & 1.853534e-11 & 315.139802 & 5.223740e-06 & 4.561398e-06 \\ \hline
7 & 6.677478e-01 & 1.803535e-09 & 300.718944 & 2.275145e-04 & 1.384767e-04 \\ \hline
\end{tabular}
\caption{Performance of the Nehari-constrained Deep Ritz method with mixed boundary conditions \eqref{E:mixed} for $Q^{\prime\prime}=Q-Q^p$.}
\label{tab:Nehari_Deep_ritz_psweeps}
\end{table}

We observe that the Nehari projection gives a more stable approximation across the nonlinear powers tested. The method achieves $L^2$ errors of orders $10^{-7}$, $10^{-6}$, and $10^{-4}$ for $p=3,\,5,$ and $7$, respectively. The higher-power case remains the most difficult, but the deterioration from $p=5$ to $p=7$ is much less severe than that observed for the constrained Deep Ritz method in Table~\ref{tab:Deep_ritz_psweeps}. In particular, for $p=7$, the Nehari-constrained method reduces the $L^2$ error from order $10^{-3}$ to $10^{-4}$. We interpret this as evidence that the Nehari projection provides a more effective mechanism for selecting the correct branch when the profile becomes sharper.

More specifically, Table~\ref{tab:Nehari_Deep_ritz_psweeps} shows that, for $p=3$, both the $L^\infty$ and $L^2$ errors are of order $10^{-7}$. For $p=5$, both errors are of order $10^{-6}$, while for $p=7$ they are of order $10^{-4}$. Thus, we obtain a gradual deterioration in profile accuracy as $p$ increases, rather than the larger loss of accuracy observed for the constrained Deep Ritz formulation at $p=7$.

\begin{figure}[!htbp]
\centering
\includegraphics[width=0.49\textwidth, height=0.32\textwidth]{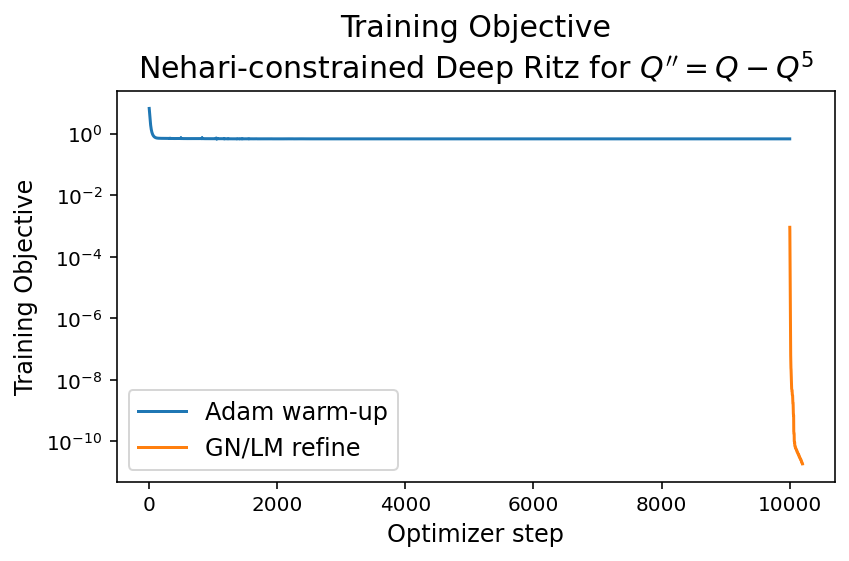}
\includegraphics[width=0.49\textwidth, height=0.32\textwidth]{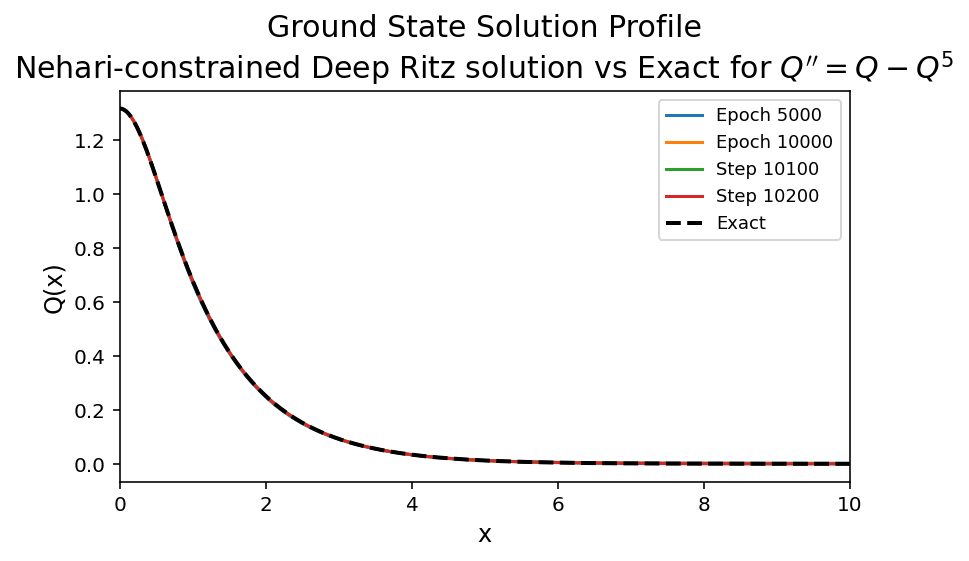}
\caption{Training objective (left) and ground state solution profile convergence (right) for the Nehari-constrained Deep Ritz method (with parameters as seen in Table~\ref{tab:nehari_deep_ritz_baselines}) for $Q^{\prime\prime}=Q-Q^5$. The (black dashed) curve labeled \emph{Exact} in the right panel is the explicit solution given by \eqref{Qexact}.  The final learned profile corresponds to $10{,}000$ Adam steps followed by $200$ GN/LM steps.}
\label{fig:NDR_converge}
\end{figure}

The improvement in robustness comes with a moderate increase in runtime. On the primary device, the Nehari-constrained runs require approximately $301$--$315$ seconds, compared with approximately $247$--$264$ seconds for the constrained Deep Ritz method. Thus, the Nehari formulation is not the fastest Deep Ritz variant, but we find that it gives a better balance of accuracy across all three powers.

In the left panel of Figure \ref{fig:NDR_converge}, one can observe that the Adam stage stabilizes the Nehari-constrained variational objective, after which the GN/LM stage further reduces the residual-refinement objective. In the right panel of the same figure, we compare the learned profiles with the exact ground state given by \eqref{Qexact}. The learned profiles are approaching the exact curve throughout the plotted interval, including near the peak at the left endpoint, even though the exact peak amplitude is not prescribed by the formulation. The final profile after $10{,}000$ Adam steps and $200$ GN/LM steps is not only visually  indistinguishable from the exact solution at the plotting scale, but also its remaining numerical discrepancy is consistent with the $p=5$ errors reported in Table~\ref{tab:Nehari_Deep_ritz_psweeps}: the final $L^\infty$ and $L^2$ errors are $5.223740\times 10^{-6}$ and $4.561398\times 10^{-6}$, respectively, and are therefore, both of order $10^{-6}$.

\FloatBarrier
\subsubsection{Comparison with the ``secondary” computational device for Nehari-constrained Deep Ritz}\label{S:secondary_nehari}

Table~\ref{tab:Nehari_Deep_ritz_psweeps_secondary} shows that the Nehari-constrained Deep Ritz method benefits substantially from the secondary GPU+CPU workstation. On the primary device, the three runs require approximately $301$--$315$ seconds, while on the secondary device the runtime is reduced to approximately $95$--$97$ seconds. We find that this corresponds to a speedup of about three across the powers tested, with very little variation between $p=3$, $p=5$, and $p=7$. Thus, we observe that the secondary device provides a consistent acceleration rather than an improvement limited to a single exponent.

\begin{table}[!htbp]
\centering
\begin{tabular}{|c||c|c|c|c|c|}
\hline
$p$
& $\mathcal{L}_{\mathrm{NDR}}(\theta)$ 
& $\mathcal{J}_{\mathrm{GN}}(\theta)$
& GPU+CPU time
& $L^\infty$ error
& $L^2$ error \\ \hline
\hline
3 & 6.666667e-01 & 5.934479e-13 & 96.631522 & 3.255530e-07 & 5.289161e-07 \\ \hline
5 & 6.801748e-01 & 1.871966e-11 & 94.922791 & 5.221092e-06 & 4.573099e-06 \\ \hline
7 & 6.677480e-01 & 2.051331e-09 & 96.848960 & 2.730569e-04 & 1.542854e-04 \\ \hline
\end{tabular}
\caption{Secondary-device GPU+CPU performance of the Nehari-constrained Deep Ritz method with mixed boundary conditions \eqref{E:mixed} for $Q^{\prime\prime}=Q-Q^p$.}
\label{tab:Nehari_Deep_ritz_psweeps_secondary}
\end{table}

Importantly, this reduction in runtime does not substantially change the accuracy pattern. The final errors remain close to the primary-device values and preserve the same ordering across the nonlinear powers: the method is most accurate for $p=3$, remains highly accurate for $p=5$, and is least accurate for the sharper $p=7$ profile. 
For $p=3$, both the secondary-device errors remain of order $10^{-7}$; for $p=5$, both remain of order $10^{-6}$; and for $p=7$, both remain of order $10^{-4}$. Thus, we observe exactly the same error orders as on the primary device. We find that the main effect of moving the Nehari-constrained Deep Ritz method to the secondary device is a roughly threefold reduction in wall-clock time, while the numerical accuracy remains essentially comparable.

The close agreement between the two devices is also reflected in the residual-refinement objective. The values of $\mathcal{J}_{\mathrm{GN}}(\theta)$ remain of the same order on both devices for every power tested: $10^{-13}$ for $p=3$, $10^{-11}$ for $p=5$, and $10^{-9}$ for $p=7$, differing only in the leading constant. Thus, unlike pure GN/LM in Section~\ref{S:secondaryPINN}, the Nehari-constrained formulation does not exhibit a device-dependent change of solution branch.

\subsection{Weinstein Deep Ritz Formulation}\label{S:WeinsteinDeepRitz}

We next introduce a variational neural-network method based on the Weinstein
quotient rather than on the action. 
Following the classical variational characterization \cite{Weinstein1983}, we use the Weinstein quotient to define the best constant as
\begin{equation}\label{E:R-b-WDR}
R_b(p) := \inf_{u\in H^1(\R)\setminus\{0\}}
  \frac{q_b(u)} {\|u\|_{L^{p+1}(\R)}^2},
\end{equation}
where the quadratic form $q$ is defined as
\begin{equation}\label{E:q-b-WDR}
q_b(u)
:=
\|u'\|_{L^2(\R)}^2+b\|u\|_{L^2(\R)}^2.
\end{equation}
For brevity, we drop the dependence of $R_b$ on $p$, and simply write $R_b$.
If $v$ is a positive optimizer normalized by
$\|v\|_{L^{p+1}(\R)}=1$,
then $q_b(v)=R_b$, and the corresponding Euler--Lagrange equation is
\begin{equation}\label{E:v-Euler-WDR}
  -v^{\prime\prime}+bv=R_bv^p.
\end{equation}
Consequently,
\begin{equation}\label{E:Q-from-v-WDR}
  Q=\left(\frac{R_b}{\gamma}\right)^{\frac{1}{p-1}}v
\end{equation}
solves the initial ground state equation \eqref{2ndODE-again}.

This characterization is especially convenient for a Deep Ritz
implementation. The quotient is invariant under multiplication of $u$ by a
nonzero constant, so it does not select the trivial solution. We remove the
remaining scale invariance by normalizing the neural-network trial function
exactly. The optimization determines the normalized shape $v$ and the value
$R_b$; the physical amplitude of $Q$ is then recovered from
\eqref{E:Q-from-v-WDR}. In particular, for fixed $b$ and $p$, the variational
training is independent of $\gamma$, which enters only through the final
amplitude rescaling.

Since the ground state is even, we work on $[0,L]$ and interpret each trial
function by its even extension to $[-L,L]$. Let $q_\theta$ be a neural network
and let $\rho:\R\to(0,\infty)$ be a smooth positivity map, for example,
$\rho(z)=\log(1+e^z)$. We use the hard-constrained raw trial function
\begin{equation}\label{E:WDR-raw-ansatz}
w_\theta(x)
= \left(1-\frac{x^2}{L^2}\right)
  \rho\!\left(q_\theta\!\left(\frac{x^2}{L^2}\right)\right),
  \qquad 0\le x\le L,
\end{equation}
which gives
$$
w_\theta'(0)=0,\qquad w_\theta(L)=0,
\qquad w_\theta(x)>0\quad\text{for }0\le x<L.
$$

We normalize the even extension of $w_\theta$ by defining 
\begin{equation}\label{E:WDR-normalization}
v_\theta(x):=A_\theta^{-\frac{1}{p+1}}w_\theta(x), ~~\mbox{where} ~~A_\theta =2\int_0^L |w_\theta(x)|^{p+1}\,dx.
\end{equation}
We then have  $2\int_0^L |v_\theta(x)|^{p+1}\,dx=1$.

The truncated Weinstein quotient, evaluated on the even extension, is
\begin{align}
R_{b,L}(\theta) :=
2\int_0^L \left(|v_\theta'(x)|^2+b|v_\theta(x)|^2\right)\,dx
=\frac{ 2 \int_0^L \left(|w_\theta'(x)|^2 + b|w_\theta(x)|^2\right) \,dx }{\left( 2 \int_0^L|w_\theta(x)|^{p+1} \,dx\right)^{\frac{2}{p+1}} }.
\label{E:WDR-Rayleigh-quotient}
\end{align}
The Weinstein Deep Ritz loss is then defined as
\begin{equation}\label{E:WDR-loss}
\mathcal L_{\mathrm{WDR}}(\theta) = R_{b,L}(\theta).
\end{equation}
We point out that no other normalization penalty, amplitude condition, or Dirichlet boundary penalty is needed because we built these constraints into
\eqref{E:WDR-raw-ansatz}--\eqref{E:WDR-normalization}.

After minimizing \eqref{E:WDR-loss}, following the definition of $Q$ in \eqref{E:Q-from-v-WDR}, we set
\begin{equation}\label{E:WDR-reconstruction}
Q_\theta(x) = \left(\frac{R_{b,L}(\theta)}{\gamma}\right)^{\frac{1}{p-1}}
v_\theta(x).
\end{equation}
At an exact constrained minimizer {(i.e., a minimizer of the underlying function space variational problem, rather than a minimizer over a fixed neural-network class)}, the normalized profile satisfies
\begin{equation}\label{E:WDR-Euler-finite_wrong}
-v_\theta^{\prime\prime}+bv_\theta =R_{b,L}(\theta)v_\theta^p,
\qquad v_\theta'(0)=0,\qquad v_\theta(L)=0,
\end{equation}
and also $-Q_\theta^{\prime\prime}+bQ_\theta-\gamma Q_\theta^p=0.$

During and after training, we monitor the Euler--Lagrange residual
\begin{equation}\label{E:WDR-EL-residual}
r_{v,\theta}(x) = - v_\theta^{\prime\prime}(x)+ b v_\theta(x)
- R_{b,L}(\theta) (v_\theta(x))^{p}.
\end{equation}
The corresponding profile residual satisfies
\begin{equation}\label{E:res-ODE}
-Q_\theta^{\prime\prime}+bQ_\theta-\gamma Q_\theta^p
=\left(\frac{R_{b,L}(\theta)}{\gamma}\right)^{\frac{1}{p-1}} r_{v,\theta}.
\end{equation}
Thus, either residual \eqref{E:WDR-EL-residual} or \eqref{E:res-ODE} may be used as a post-training diagnostic. If a GN/LM stage is used, it is applied only after the Weinstein-quotient pretraining to a least-squares vector assembled from $r_{v,\theta}$ at residual points. The
normalization and the hard boundary structure are retained during this
refinement, so that the optimization cannot collapse to the zero solution.

Unlike the action-based Deep Ritz formulation, the Weinstein Deep Ritz method does not require the value $Q(0)=Q_0$ as input. It obtains the normalized profile shape from the quotient minimization and then recovers the correct amplitude through \eqref{E:WDR-reconstruction}. 
{ Thus, forming the WDR loss and reconstruction does not require an explicit profile, which makes this approach promising for computing ground states with no explicit solutions, which, for example, we investigate in the follow-up work \cite{HRR-cNLS}.} 
The Weinstein Deep Ritz implementation is summarized in Algorithm~\ref{alg:WeinsteinDeepRitz}.

{\scriptsize
\begin{algorithm}[!htbp]
\caption{Weinstein Deep Ritz Method for
$Q^{\prime\prime}=bQ-\gamma Q^p$}
\begin{algorithmic}[1]
\State Choose $b>0$, $\gamma>0$, $p>1$, a truncated domain $[0,L]$,
quadrature points $\{x_j\}_{j=1}^{N_q}$, and optional residual-refinement
points $\{\widehat x_j\}_{j=1}^{N_r}$.
\State Define a neural network $q_\theta$ and a smooth positivity map
$\rho:\R\to(0,\infty)$  such that 
$$\rho(s)=\log(1+e^s).$$
\State Define the hard-constrained raw trial function
$$
 w_\theta(x)
 =\left(1-\frac{x^2}{L^2}\right)
  \rho\!\left(q_\theta\!\left(\frac{x^2}{L^2}\right)\right),
$$
which satisfies $w^\prime_\theta(0)=0$ and $w_\theta(L)=0$.
\State Use automatic differentiation to compute $w_\theta'(x_j)$.
\State Compute the normalization
$$
A_\theta=2\int_0^L |w_\theta|^{p+1}\,dx, \qquad
 v_\theta=A_\theta^{-1/(p+1)}w_\theta.
$$
\State Compute the truncated Weinstein quotient
$$
 R_{b,L}(\theta)
 =2\int_0^L\left(|v_\theta'|^2+b|v_\theta|^2\right)\,dx.
$$
\State Define the variational loss
$$
 \mathcal L_{\mathrm{WDR}}(\theta)=R_{b,L}(\theta).
$$
\State Choose a quotient-minimization strategy
$$
 \mathsf{Opt}_{\mathrm{WDR}}
 \in\{\text{Adam},\ \text{L--BFGS},\
       \text{Adam}\to\text{L--BFGS}\}.
$$
\For{$n=1,2,\dots,N_{\mathrm{epochs}}$}

\State Evaluate the quadrature approximations of $A_\theta$ and
$R_{b,L}(\theta)$ on the fixed uniform quadrature grid.

\State Differentiate $\mathcal L_{\mathrm{WDR}}$ with respect to $\theta$
and update $\theta$ using the variational stage of
$\mathsf{Opt}_{\mathrm{WDR}}$.

\State Store $R_{b,L}(\theta)$ together with the normalization and
boundary diagnostics.

\EndFor

\If{GN/LM refinement is selected}

\State Recompute $A_\theta$, $v_\theta$, and $R_{b,L}(\theta)$ and form the
Euler--Lagrange residual
$$
r_{v,\theta}(x)
=-v_\theta^{\prime\prime}(x)
+bv_\theta(x)-R_{b,L}(\theta)
\left(v_\theta(x)\right)^p.
$$

\State Form $F_{\mathrm{GN}}(\theta)$ from the quadrature-weighted values
of $r_{v,\theta}$ at the GN/LM residual points.

\State Apply GN/LM to minimize
$$
\mathcal J_{\mathrm{GN}}(\theta)
=\frac12\|F_{\mathrm{GN}}(\theta)\|_2^2,
$$
while retaining the normalized, hard-constrained representation of
$v_\theta$.

\EndIf

\State Reconstruct the ground-state profile
$$
Q_\theta
=
\left(
\frac{R_{b,L}(\theta)}{\gamma}
\right)^{1/(p-1)}
v_\theta.
$$

\State Evaluate the pointwise ground-state residual and the mass, energy,
and finite-interval Pokhozhaev diagnostics for $Q_\theta$.

\State Return $v_\theta$, $R_{b,L}(\theta)$, and the reconstructed profile
$Q_\theta$
\end{algorithmic}
\label{alg:WeinsteinDeepRitz}
\end{algorithm}
}

The method first trains a normalized profile by minimizing the truncated Weinstein quotient and then reconstructs the physical ground-state amplitude from the computed quotient value.

For the GN/LM refinement, the residual points are uniformly spaced on $[0,L]$ and the residual vector is formed using the square-root composite-trapezoidal weights,
$$
F_{{\rm GN},j}(\theta) =\sqrt{\omega_j}\,r_{v,\theta}(\widehat x_j),
\qquad j=1,\ldots,N_r,
$$
where $\omega_1=\omega_{N_r}=\Delta x/2$, $\omega_j=\Delta x$ for $2\leq j\leq N_r-1$, and $\Delta x=L/(N_r-1)$. No additional $N_r^{-1/2}$ normalization is used, unlike in the PINN formulation. The general GN/LM damping and backtracking procedure is the one described above. During the variational stage, $N_q$ denotes the fixed quadrature grid used to evaluate the Weinstein quotient. During GN/LM refinement, the residual grid is used to form $F_{\rm GN}$, while a separate quotient grid is used to recompute the normalization and $R_{b,L}(\theta)$ at each trial parameter vector.

Thus, the hyperparameter sweep in Table~\ref{tab:weinstein_deepritz_hyperparams} varies both the neural-network approximation parameters and the optimization parameters that affect the quotient minimization and subsequent residual refinements. In particular, we varied the activation function, network depth and width, number of quadrature points used to evaluate the quotient, Adam warm-up length, L--BFGS settings, and GN/LM refinement parameters. Since the normalization and boundary structure are built into the Weinstein trial function, no separate amplitude penalty, normalization penalty, or boundary penalty is required.

{\scriptsize
\begin{table}[!htbp]
\centering
\begin{tabular}{ll}
\toprule
\textbf{Parameter} & \textbf{Values tested / Description} \\
\midrule
Domain length $L$ & $30$ \\
Boundary condition & mixed \\

Trial function 
& Soft penalty trial, hard or decaying Dirichlet ansatz \\

Optimizer
& Adam, L--BFGS, or Adam--L--BFGS\\

Refinement strategy 
& GN/LM or Adam--GN/LM \\

Activation functions 
& $\tanh(x)$, $\mathrm{SiLU}(x)$, $\mathrm{SiLUPLUS}(x)$, 
$e^{-x^2}$, $e^{-x^4}$, and $\sech(x)$ \\

Number of hidden layers & $3$, $4$, or $5$ \\

Neurons per layer & $32$, $64$, or $96$ \\

Number of quadrature points $N_q$ & $2^8$, $2^9,\,2^{10},$ $2^{11}$, or $2^{12}$ \\

Adam warm-up length 
& $2{,}500$, $5{,}000$, $7{,}500$, $10{,}000$,  $12{,}500$, or $15{,}000$ steps \\

Adam learning rate $lr$ &  $10^{-3}$ initial; StepLR decay \\

L--BFGS refinement steps 
& $25$, $50$, $75$, $100$, $125$, $150$, or $200$ \\

L--BFGS history size & $20$, $50$, or $100$ \\

L--BFGS max iterations per step & $1$, $5$, $10$, $20$, or $30$ \\

L--BFGS learning rate & $0.5$ or $1.0 $\\

GN/LM refinement steps 
& $25$, $50$, $75$, $100$, $125$, $150$, or $200$ \\

GN/LM initial damping 
& $5\times 10^{-4}$, $10^{-3}$, or $2\times 10^{-3}$ \\

GN/LM Maximum backtracking steps 
& $2$, $4$, 6, 8, and 10 \\

GN/LM Damping shrink factor 
& Values in $[0.4,1.1]$ \\

GN/LM Damping growth factor 
& $5$, $10$, and $15$ \\

Selection criterion 
& Smallest $L^2$ error, with runtime used to break near-ties \\
\bottomrule
\end{tabular}
\caption{Weinstein Deep Ritz architecture, constraint, and optimization hyperparameters tested for solving \eqref{2ndODE} on $[0,30]$.}
\label{tab:weinstein_deepritz_hyperparams}
\end{table}
}

\begin{table}[!htbp]
\centering
\setlength{\tabcolsep}{3pt}
\renewcommand{\arraystretch}{1.05}
\resizebox{\textwidth}{!}{%
\begin{tabular}{cccccccccccccc}
\hline
\textbf{Optimization} 
& \textbf{Refinement} 
& \textbf{Activation} 
& \shortstack[c]{\textbf{Hidden}\\\textbf{Layers}}
& \shortstack[c]{\textbf{Neurons}\\\textbf{per Layer}}
& $N_q$
& \shortstack[c]{\textbf{Adam}\\ \textbf{Warm-up}}
& \shortstack[c]{\textbf{GN/LM}\\ \textbf{Steps}}
& \shortstack[c]{\textbf{GN/LM Max}\\ \textbf{Backtracks}}
& \shortstack[c]{\textbf{Damping}\\ \textbf{Shrink}}
& \shortstack[c]{\textbf{Damping}\\ \textbf{Growth}}
& \shortstack[c]{\textbf{Initial}\\ \textbf{Damping}}
& \shortstack[c]{\textbf{Residual}\\\textbf{Points}}
& \shortstack[c]{\textbf{Quotient}\\\textbf{Points}}
\\
\hline
Adam
& GN/LM 
& $\mathrm{SiLU}(x)$
& $5$ 
& $32$ 
& $2^9$
& $5{,}000$ 
& $200$ 
& $6$ 
& $0.5$
& $15$
& $10^{-3}$
& $128$
& $128$
\\
\hline
\end{tabular}
}
\caption{Selected Weinstein Deep Ritz configurations for mixed boundary conditions \ref{E:mixed}.}
\label{tab:weinstein_deep_ritz_baselines}
\end{table}

The selected configuration is reported in Table~\ref{tab:weinstein_deep_ritz_baselines}. This baseline uses an Adam-trained Weinstein quotient stage followed by GN/LM residual refinement, with the quadrature and residual-point counts fixed in order to balance accuracy with memory usage and runtime. All Weinstein Deep Ritz results reported below use this fixed configuration unless otherwise stated, so the subsequent tables measure the performance of the selected Weinstein baseline across the nonlinear powers rather than a new hyperparameter search for each case.

\begin{table}[!htbp]
\centering
\begin{tabular}{|c||c|c|c|c|c|}
\hline
$p$
& $\mathcal{L}_{\mathrm{WDR}}(\theta)$ 
& $\mathcal{J}_{\mathrm{GN}}(\theta)$
& CPU time
& $L^\infty$ error
& $L^2$ error \\
\hline\hline
3 & 2.309401e+00 & 5.781202e-15 & 148.970724 & 4.832268e-08 & 6.092994e-08 \\ \hline
5 & 2.553766e+00 & 3.953484e-14 & 150.616751 & 5.235757e-06 & 2.322385e-06 \\ \hline
7 & 2.592439e+00 & 1.687038e-09 & 145.797130 & 2.802301e-04 & 1.451584e-04 \\ \hline
\end{tabular}
\caption{Performance of the Weinstein Deep Ritz method with mixed boundary conditions \eqref{E:mixed} for $Q^{\prime\prime}=Q-Q^p$.}
\label{tab:Weinstein_Deep_ritz_psweeps}
\end{table}

Table~\ref{tab:Weinstein_Deep_ritz_psweeps} reports the performance of the Weinstein Deep Ritz method for $Q''=Q-Q^p$, with $p=3,5,$ and $7$. Unlike the action-based Deep Ritz and Nehari-constrained formulations, this method does not minimize the action directly. Instead, the Adam stage minimizes the Weinstein quotient $\mathcal{L}_{\mathrm{WDR}}(\theta)$, which selects a normalized profile shape. The physical amplitude is then recovered from the computed quotient value using the reconstruction formula \eqref{E:WDR-reconstruction}. As before, $\mathcal{J}_{\mathrm{GN}}(\theta)$ denotes the separate residual-based GN/LM refinement objective, not the quotient itself.

We find that the selected Weinstein Deep Ritz configuration gives the most favorable error--runtime balance among the Deep Ritz variants tested. For the cubic case, the method achieves errors of order $10^{-8}$, while for the quintic case the errors are of order $10^{-6}$. The septic case remains the most challenging, with errors increasing to the $10^{-4}$ level. Thus, we observe a gradual deterioration in accuracy as the profile becomes sharper, but the Weinstein formulation still remains accurate across all three nonlinear powers.

More precisely, Table~\ref{tab:Weinstein_Deep_ritz_psweeps} shows that both errors are of order $10^{-8}$ for $p=3$ and order $10^{-6}$ for $p=5$. For $p=7$, both errors increase to order $10^{-4}$. Thus, we obtain the same gradual loss of accuracy with increasing $p$ as in the Nehari-constrained formulation, but without the order $10^{-3}$ error observed for the action-based constrained Deep Ritz method.

The runtime is also substantially lower than for the other Deep Ritz variants. The Weinstein runs require approximately $146$--$151$ seconds, whereas the constrained Deep Ritz method requires approximately $247$--$264$ seconds and the Nehari-constrained method requires approximately $301$--$315$ seconds. Hence, among the selected configurations, the Weinstein formulation is the most accurate Deep Ritz variant for $p=3$ and the fastest for all three powers. We observe that for $p=5$, the Nehari-constrained formulation gives the slightly smaller $L^\infty$ error, while the Weinstein formulation gives the smaller $L^2$. For $p=7$, its errors are comparable to those of the Nehari-constrained method, with both methods producing errors of order $10^{-4}$, while the Weinstein method retains the shorter runtime. Hence, we find that the Weinstein formulation gives the best overall balance between profile accuracy and computational cost among the three selected Deep Ritz variants we consider in this paper.

\begin{figure}[!htbp]
\centering
\includegraphics[width=0.49\textwidth, height=0.32\textwidth]{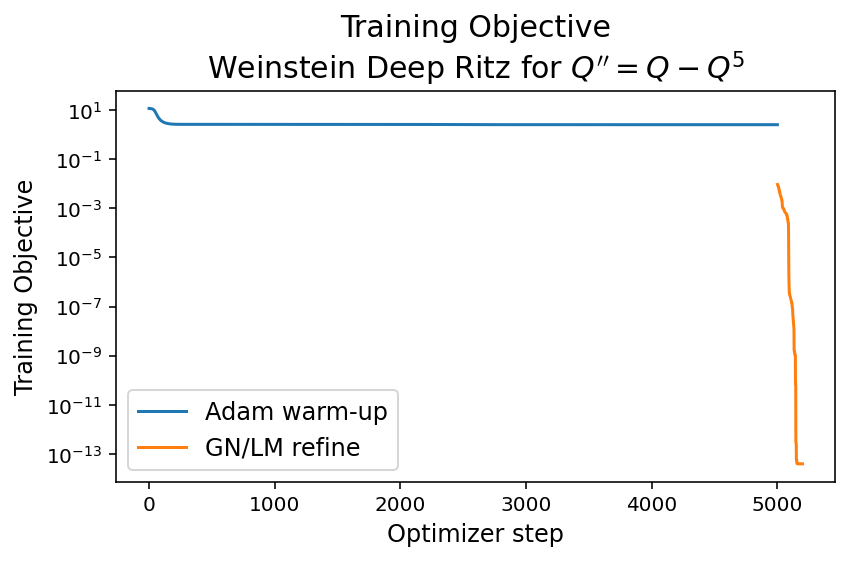}
\includegraphics[width=0.49\textwidth, height=0.32\textwidth]{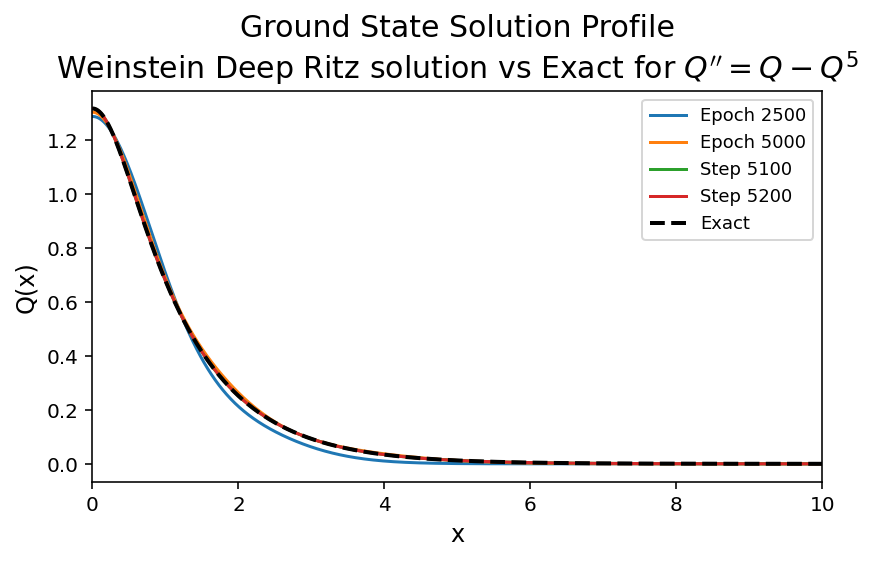}
\caption{Training objective (left) and ground state solution profile convergence (right) for the Weinstein Deep Ritz method (with parameters as seen in Table~\ref{tab:weinstein_deep_ritz_baselines}) for $Q^{\prime\prime}=Q-Q^5$. 
The (black dashed) curve labeled \emph{Exact} in the right panel is the explicit solution given by \eqref{Qexact}. The final learned profile corresponds to $5{,}000$ Adam steps followed by $200$ GN/LM steps.}
\label{fig:WDR_converge}
\end{figure}

Figure~\ref{fig:WDR_converge} illustrates the two-stage behavior of the method. In the left panel, we observe that the Adam warm-up reduces the Weinstein quotient, thereby improving the normalized profile shape and quotient value. The GN/LM stage then acts as a residual-refinement step and sharply reduces the residual-based objective. In the right panel of Figure ~\ref{fig:WDR_converge}, we compare the 
reconstructed profiles with the exact ground state given by \eqref{Qexact}. We observe that the Adam-trained approximation already captures the main shape and amplitude of the ground state, while the GN/LM stage further improves the agreement with the exact profile. This behavior is consistent with the structure of the method: the variational stage determines the normalized shape, while the final amplitude is obtained only after the quotient minimization. The final profile after $5{,}000$ Adam steps and $200$ GN/LM steps is visually indistinguishable from the exact solution at the plotting scale. Its remaining discrepancy agrees with the $p=5$ values in Table~\ref{tab:Weinstein_Deep_ritz_psweeps}: the final $L^\infty$ and $L^2$ errors are $5.235757\times 10^{-6}$ and $2.322385\times 10^{-6}$, respectively, and are hence, both of order $10^{-6}$. Although the computations are performed on the full interval $[0,30]$, the solution profiles are displayed on $[0,10]$ to emphasize the left side of the profile and the peak-amplitude region.

\FloatBarrier
\subsubsection{Comparison with the ``secondary” computational device for Weinstein Deep Ritz}\label{S:secondary_weinstein}

Table~\ref{tab:Weinstein_Deep_ritz_psweeps_secondary} reports the same Weinstein Deep Ritz baseline on the secondary GPU+CPU workstation.  We observe that the main effect of the secondary device is a reduction in runtime from approximately $146$--$151$ seconds to approximately $40$--$44$ seconds.  This corresponds to a speedup of about $3.3$--$3.7$ across the three nonlinear powers.  The acceleration is consistent across $p=3$, $p=5$, and $p=7$, indicating that the secondary device improves the practical wall-clock cost without changing the qualitative behavior of the method.

\begin{table}[!htbp]
\centering
\begin{tabular}{|c||c|c|c|c|c|}
\hline
$p$
& $\mathcal{L}_{\mathrm{WDR}}(\theta)$ 
& $\mathcal{J}_{\mathrm{GN}}(\theta)$
& GPU+CPU time
& $L^\infty$ error
& $L^2$ error \\ \hline
\hline
3 & 2.309401e+00 & 1.934108e-14 & 40.420787 & 9.883258e-08 & 1.186556e-07 \\ \hline
5 & 2.553766e+00 & 5.054238e-14 & 44.468209 & 5.204065e-06 & 2.321890e-06 \\ \hline
7 & 2.592439e+00 & 1.258231e-09 & 44.377295 & 3.016073e-04 & 1.542014e-04 \\ \hline
\end{tabular}
\caption{Secondary-device GPU+CPU performance of the Weinstein Deep Ritz method with mixed boundary conditions \eqref{E:mixed} for $Q^{\prime\prime}=Q-Q^p$.}
\label{tab:Weinstein_Deep_ritz_psweeps_secondary}
\end{table}

The secondary-device errors remain of the same orders as the corresponding primary-device results. For $p=3$, the $L^\infty$ error remains of order $10^{-8}$, while the $L^2$ error changes from order $10^{-8}$ to order $10^{-7}$. For $p=5$, both errors remain of order $10^{-6}$, and for $p=7$, both remain of order $10^{-4}$. Thus, we find that the GPU+CPU implementation mainly improves runtime while preserving the same accuracy pattern.  In particular, the higher-power case remains the most difficult, but the secondary-device result stays comparable to the original Weinstein run.

The quotient values are also essentially unchanged between the two devices for all three powers. Similarly, the primary- and secondary-device residual-refinement objectives remain of the same orders for $p=5$ and $p=7$. For $p=3$, $\mathcal{J}_{\mathrm{GN}}(\theta)$ changes from order $10^{-15}$ to order $10^{-14}$, while the corresponding profile errors remain very small. Thus, we observe no evidence of a device-dependent change in the selected solution branch.

Overall, we find that the Weinstein Deep Ritz method gives the best error--runtime combination among the Deep Ritz variants considered  here for the selected runs.  It is the fastest formulation and gives the smallest errors for $p=3$, while remaining comparable to the Nehari-constrained formulation for $p=5$ and $p=7$.  The Nehari-constrained Deep Ritz method, however, remains competitive for the sharper $p=7$ profile and selects the same solution branch in the reported runs without prescribing the exact peak amplitude.  Thus, we regard the Weinstein formulation as the most efficient and often the most accurate variant, among the selected runs, while Nehari provides a robust constrained alternative whose accuracy remains comparable for the
higher-power case.

\FloatBarrier
\section{Conclusions and Discussion}\label{S:Conclusion}

We have developed and compared enhanced neural-network-based methods for computing {\it one-dimensional} ground-state profiles for several classical models in nonlinear dispersive equations. The profile equation $Q''=bQ-\gamma Q^p$, arising from the NLS, gKdV, and NLKG models, admits explicit $\sech$-type solutions for all integer power 
$p>1$, thus, providing exact benchmarks for systematic error analysis. Building on the comparison in~\cite{HRW2026}, we investigated three directions of improvement over the standard Adam-trained PINN baseline: (i) higher-order optimization (GN/LM refinement), (ii) Fourier-feature PINNs, and (iii) constrained Deep Ritz formulations.

The main findings of this study are as follows.

{\it Higher-order optimization improves PINN accuracy.}  Replacing or supplementing standard Adam training with Gauss--Newton/Levenberg--Marquardt refinement consistently reduces the pointwise and $L^2$ errors of standard PINNs. An Adam warm-up followed by GN/LM refinement (Adam--GN/LM) is generally more robust than pure GN/LM, which is sensitive to the damping parameter and initialization and may converge to an incorrect solution branch, such as a double-bump profile. GN/LM refinement generally gives the best results for the standard PINN and Deep Ritz formulations, whereas L--BFGS performs better for the FF--PINN formulation.

{\it Fourier-feature PINNs provide modest additional gains.} Embedding the spatial coordinate into a trigonometric feature space before passing it through the multilayer perceptron is designed to reduce the spectral bias of coordinate-based networks and, in some configurations, we observe improved 
convergence. However, the gains from Fourier features are smaller than those from higher-order optimization, and the method introduces additional hyperparameters including the scale of the frequency matrix 
and whether the raw spatial coordinate is included.

{\it Deep Ritz variants outperform residual-based methods.} Among all neural-network methods tested, the variational Deep Ritz formulations produce the strongest overall performance. By minimizing the action or related variational objectives instead of a pointwise residual, these methods exploit the variational structure of the ground-state equation and yield lower errors for comparable or shorter runtimes.

{\it The Weinstein Deep Ritz method gives the best error--runtime balance overall.} The Weinstein Deep Ritz formulation, based on minimizing the Weinstein quotient~\eqref{E:Weinstein-Rb} over a normalized, hard-constrained trial class, achieves the smallest errors among the tested neural-network methods for $p=3$, while also requiring the shortest runtimes, among the Deep Ritz variants (approximately 146--151 seconds on the primary device versus 247--264 seconds for the constrained Deep Ritz and 301--315 seconds for the Nehari-constrained variant). For $p=5$ and $p=7$, the Weinstein and Nehari-constrained formulations produce errors on the same order ($\sim 10^{-6}$ and $\sim 10^{-4}$, respectively). A key advantage of the Weinstein method is that it does not require prescribing the exact peak amplitude $Q(0)$; the normalized profile shape is determined by the quotient minimization, and the physical amplitude is recovered through the reconstruction formula~\eqref{E:Weinstein-recovery}. This makes the method promising to problems where no explicit solution is available.

{\it The Nehari-constrained Deep Ritz method provides a robust alternative.} The Nehari-constrained formulation, which projects the amplitude onto the Nehari manifold and supplements the action with a Pokhozhaev penalty, provides a stable branch-selection mechanism without prescribing the exact peak amplitude. Its accuracy remains competitive, particularly for the sharper $p=7$ profile.

{\it{The dependence of accuracy on the nonlinear power is method-dependent.}} The numerical results do not show a uniform deterioration in accuracy as $p$ increases from $3$ to $7$. For the Deep Ritz formulations, the higher-power cases generally produce larger errors. For the PINN formulations, the trend depends on the optimizer: the pure L--BFGS errors are smaller at $p=7$ than at $p=3$, the pure GN/LM errors decrease from $p=3$ to $p=5$ before increasing at $p=7$, and the Adam--L--BFGS and Adam--GN/LM configurations exhibit different trends across the tested powers. The FF--PINN errors remain of comparable order across the three powers. Thus, the observed dependence on the nonlinear power varies with the selected formulation and optimizer.

{\it GPU acceleration preserves accuracy while reducing runtime.} Moving from the primary CPU-based device to the secondary GPU+CPU workstation consistently reduces wall-clock runtimes, with primary-to-secondary runtime ratios ranging from approximately $1.6$ to $9.7$ across the reported configurations, while the profile errors remain of comparable order for most of the reported configurations.

{\it Neural-network methods do not yet match classical Petviashvili solvers.} Although the enhanced neural-network methods substantially improve upon the Adam-trained PINN baseline, they do not yet achieve in 1D the accuracy or computational efficiency of the Petviashvili iteration reported in~\cite{HRW2026}. For the one-dimensional benchmark problem, the Petviashvili method achieves machine-precision accuracy, on the order of $ 10^{-14}$, in shorter runtimes. Further studies would be needed, for example, for higher dimensions and more complex nonlinearities. 

The methods developed here provide a one-dimensional proof of concept, they { are intended to} establish a foundation for extending neural-network-based ground-state computation to settings where classical methods face greater challenges: higher-dimensional problems subject to the curse of dimensionality, equations with combined or non-standard nonlinearities where explicit solutions are unavailable, and problems involving fractional or higher-order differential operators. The Weinstein Deep Ritz method, which we introduced in this paper, does not rely on knowledge of the exact solution and is therefore  { prospectively} well-suited for such extensions.

%
\appendix \section{Newton direction and the strong Wolfe conditions}
\label{S:appendix-optimizers}

For reader's convenience in this appendix we give a standard derivation of the quasi-Newton search direction \eqref{eq:search-direction} and the strong Wolfe conditions used in the L--BFGS refinement stage (though the material is by now classical, see, e.g., \cite{NocedalWright2006} for a full treatment).

\paragraph{\bf Newton search direction.}
First, note that the Hessian $\nabla_\theta^2\mathcal{L}(\theta_k)$ is the matrix of second derivatives of the loss at the given current iterate $\theta_k$, thus, it describes its local second-order behavior. The second-order Taylor expansion about $\theta_k$
gives the quadratic expansion
\begin{equation}\label{eq:newton-quadratic-model}
\mathcal{L}(\theta_k+p)\approx\mathcal{L}(\theta_k)
+\nabla_\theta\mathcal{L}(\theta_k)^Tp
+\tfrac{1}{2}p^T\nabla_\theta^2\mathcal{L}(\theta_k)p,
\end{equation}
whose stationary point $p_k$ solves
\begin{equation}\label{eq:newton-system}
\nabla_\theta^2\mathcal{L}(\theta_k)\,p_k=-\nabla_\theta\mathcal{L}(\theta_k).
\end{equation}

If $\nabla_\theta^2\mathcal{L}(\theta_k)$ is positive definite, this stationary point is the unique minimizer of the quadratic model. When $\nabla_\theta^2\mathcal{L}(\theta_k)$ is non-singular, then 
$p_k=-\bigl(\nabla_\theta^2\mathcal{L}(\theta_k)\bigr)^{-1}\nabla_\theta\mathcal{L}(\theta_k)$, see 
\cite[Sec.~2.2, Eqs.~(2.14)--(2.15)]{NocedalWright2006}. 
In practice, however, the search direction is obtained by solving \eqref{eq:newton-system} without forming the inverse. Forming and storing the full $N_\theta\times N_\theta$ Hessian is
impractical for large networks, so L--BFGS instead uses a limited number of
recent parameter and gradient differences
$$
s_j=\theta_{j+1}-\theta_j, \qquad
y_j=\nabla_\theta\mathcal{L}(\theta_{j+1})-\nabla_\theta\mathcal{L}(\theta_j).
$$

A first-order expansion of the gradient gives
$y_j\approx\nabla_\theta^2\mathcal{L}(\theta_j)s_j$, so gradient changes carry
information about the action of the Hessian. 
Then the corresponding inverse secant
relation is $H_{j+1}y_j\approx s_j$ with
$H_{j+1}\approx\bigl(\nabla_\theta^2\mathcal{L}(\theta_{j+1})\bigr)^{-1}$, see
\cite[Sec.~2]{LiuNocedal1989}, \cite[Sec.~7.2, Algs.~7.4--7.5]{NocedalWright2006} for more details.
This yields the search direction \eqref{eq:search-direction}.
\smallskip

\paragraph{\bf Strong Wolfe conditions.}
Once $p_k$ is fixed, a line search minimizes
$\phi_k(\alpha)=\mathcal{L}(\theta_k+\alpha p_k)$ over $\alpha>0$. Writing
$g_k=\nabla_\theta\mathcal{L}(\theta_k)$ and $\phi_k'(0)=g_k^Tp_k$, and assuming
$p_k$ is a descent direction ($g_k^Tp_k<0$), the strong Wolfe conditions require
\begin{align}
\mathcal{L}(\theta_k+\alpha_k p_k)&\leq \mathcal{L}(\theta_k)+c_1\alpha_k g_k^Tp_k,\label{eq:strong-wolfe-decrease}\\
\bigl|\nabla_\theta\mathcal{L}(\theta_k+\alpha_k p_k)^Tp_k\bigr|&\leq c_2\,|g_k^Tp_k|,\qquad 0<c_1<c_2<1.\label{eq:strong-wolfe-curvature}
\end{align}
Condition \eqref{eq:strong-wolfe-decrease} is for sufficient decrease, and
the one in \eqref{eq:strong-wolfe-curvature} is for the strong curvature,  
which rules out two things: short steps where the loss is still decreasing steeply, and long steps where the directional derivative becomes large and positive, see \cite[Sec.~1, Eqs.~(1.1)--(1.2)]{MoreThuente1994}. Setting
$s_k=\alpha_k p_k$ and
$y_k=\nabla_\theta\mathcal{L}(\theta_{k+1})-\nabla_\theta\mathcal{L}(\theta_k)$,
the curvature condition implies $s_k^Ty_k>0$. This keeps the BFGS
inverse-Hessian approximation positive definite and preserves a descent
direction in the next iteration, as shown in \cite[Sec.~6.1, Eq.~(6.8)]{NocedalWright2006}.
Fixed steps or backtracking either depend more strongly on the learning
rate or do not by themselves enforce $s_k^Ty_k>0$, which is why the
strong Wolfe search is used.


\bibliographystyle{acm} 
\bibliography{main}

\end{document}